\documentclass{amsart}
\usepackage[english]{babel}
\usepackage{graphicx}
\usepackage{amssymb}
\usepackage{amsmath}
\usepackage{amsfonts}
\usepackage[table]{xcolor}
\usepackage{mathtools}
\usepackage{algpseudocode}
\usepackage{algorithm}
\usepackage[hidelinks]{hyperref}
\usepackage[foot]{amsaddr}
\usepackage{tikz}
\usepackage[top=1in, bottom=1.25in, left=1.25in, right=1.25in]{geometry}

\graphicspath{{images/}}
\usepackage{float}
\usepackage{subcaption}
\usepackage{multirow}
\usepackage{ulem}
\usepackage{siunitx}
\usepackage{soul}
\usepackage{booktabs}
\usepackage{overpic}

\usepackage{pifont}
\newcommand{\cmark}{\ding{51}}
\newcommand{\xmark}{\ding{55}}

\newcommand{\bmu}{{\boldsymbol{\mu}}}

\newcommand{\N}{{\mathcal{N}}}
\newcommand{\R}{{\mathbb{R}}}

\newcommand{\norm}[1]{\left\lVert#1\right\rVert}

\begin{document}

\title[A GCA approach to MOR for parametrized PDEs]{A graph convolutional autoencoder approach to model order reduction for parametrized PDEs}

\author{Federico Pichi$^{1,*}$, Beatriz Moya$^{2,*}$, Jan S. Hesthaven$^{1}$}
\address{$^1$ Chair of Computational Mathematics and Simulation Science, École Polytechnique Fédérale de Lausanne, 1015 Lausanne, Switzerland}
\address{$^2$ Aragon Institute in Engineering Research (I3A), University of Zaragoza, 50018 Zaragoza, Spain}
\address{$^*$ \textup{Both authors contributed equally to this work.}}

\begin{abstract}
    The present work proposes a framework for nonlinear model order reduction based on a Graph Convolutional Autoencoder (GCA-ROM). In the reduced order modeling (ROM) context, one is interested in obtaining real-time and many-query evaluations of parametric Partial Differential Equations (PDEs). Linear techniques such as Proper Orthogonal Decomposition (POD) and Greedy algorithms have been analyzed thoroughly, but they are more suitable when dealing with linear and affine models showing a fast decay of the Kolmogorov n-width. On one hand, the autoencoder architecture represents a nonlinear generalization of the POD compression procedure, allowing one to encode the main information in a latent set of variables while extracting their main features. On the other hand, Graph Neural Networks (GNNs) constitute a natural framework for studying PDE solutions defined on unstructured meshes. Here, we develop a non-intrusive and data-driven nonlinear reduction approach, exploiting GNNs to encode the reduced manifold and enable fast evaluations of parametrized PDEs.
    We show the capabilities of the methodology for several models: linear/nonlinear and scalar/vector problems with fast/slow decay in the physically and geometrically parametrized setting. The main properties of our approach consist of (i) high generalizability in the low-data regime even for complex behaviors, (ii) physical compliance with general unstructured grids, and (iii) exploitation of pooling and un-pooling operations to learn from scattered data.
    
    \begin{center}
    \textbf{Code availability:} \url{https://github.com/fpichi/gca-rom}
    \end{center}
\end{abstract}

\maketitle
\tableofcontents

\section{Introduction}
The numerical approximation of parameterized partial differential equations (PDEs) using standard high-fidelity techniques (i.e.\ Finite Element, Finite Volume, Spectral Element Methods) is often unfeasible in many-query and real-time contexts. The exploitation of Reduced Order Models (ROMs) \cite{BennerModelOrderReduction2020,benner2017model} can reduce the computational resources (CPU time and storage) required for their analysis and simulation. 
During the last decades, due to increasing interest and efforts, ROMs have become a well-established class of methodologies based on solid mathematical foundations. 

Among them, the Reduced Basis method \cite{hesthavenCertifiedReducedBasis2015,quarteroniReducedBasisMethods2016} enables fast and reliable evaluations of the solution for new parameter values.
To build the reduced space, which allows for these efficient computations, one usually exploits Proper Orthogonal Decomposition (POD), an SVD-based method to extract the principal components, or a Greedy algorithm, iteratively augmenting the space with a basis corresponding to the worst approximation in the parametric space.

These techniques allow decoupling the computation in two stages, offline and online, achieving good and predictable accuracy at a small computational cost in the online phase.
Nevertheless, these are usually linear approaches, which are known for losing efficiency when the models are not easily reducible, or they are characterized by non-affine and nonlinear terms.

Exploring fast and efficient non-intrusive model order reduction techniques from a deep learning perspective \cite{lee2020model, fresca2021comprehensive, vinuesa2022enhancing} allows users to overcome some of the limitations of traditional approaches. Exploiting nonlinear machine learning methodologies helps in reaching a low-dimensional representation of the latent subspace, and to capture the correlations among the features due to its optimal capacity in learning patterns. 

The autoencoder architecture, which has recently been the subject of many studies, is particularly well-suited in the ROM context. Indeed, it generalizes linear compression procedures, such as the POD. In the machine learning context, an autoencoder is constituted by nonlinear encoding and decoding structures connected through a bottleneck. The bottleneck identifies the latent dimension and plays the role of the reduced space, in which we compress the information of the high order system with the encoder. Then we ``project'' the reduced representation back to its original dimension thanks to the decoder. This is an example of an unsupervised learning task, where, instead of having input-output pairs, one seeks to approximate the identity operator by consecutively compressing and decompressing the features.

Although past works take advantage of Fully Connected Neural Networks (FCNNs) \cite{MilanoNeuralNetworkModeling2002}, studies have progressively adopted more optimized structures to exploit spatial and temporal correlations to perform more efficient training procedures \cite{hernandez2021deep,hesthaven2018non,PichiArtificialNeuralNetwork2023}. An architecture that has been thoroughly investigated is Convolutional Neural Networks (CNNs), which exploit the spatial information of the data provided, hence detecting and learning patterns. Originally meant for image classification, they have recently been applied extensively to dynamical modeling and parametrized PDEs \cite{maulik2021reduced,peng2020unsteady,wu2021reduced,lee2020model,fresca2021comprehensive}, showing great performance in several fluid dynamics benchmarks \cite{pawar2020data}. However, CNNs usually work with structured datasets, resembling the shape of an image composed by pixels. Although this is a common structure in the field of computer vision, it is a rare condition during the investigation of physical problems. One could attempt to employ Cartesian meshes to simulate dynamic behaviors \cite{morimoto2021convolutional}, but irregular meshes are often required to adapt to the domain under investigation. Indeed, PDEs are often posed on complex and possibly parametrized domains, which involve unstructured meshes.
A possible approach consists in reshaping these grids into image-like objects \cite{fresca2021comprehensive}. Despite its good performance, there is a lack of consistency in the interpretation of locality while convoluting nearby pixels. Indeed, it is hard to find an efficient ordering and reshaping which is physically consistent with the problem. As an alternative, interpolation and level set approximations reshape the data as structured meshes. However, this approach requires dense meshes not to lose information in the reshaping process, hence increasing the computational cost. For instance, in \cite{sharma2019deepinsight}, the authors propose to perform a preprocessing step to find a 2D representation of the data with $k$-PCA. Nevertheless, algorithms progressively have evolved toward their adaptation to the geometry of the domain. Works such as mesh-informed neural networks \cite{francoLearningOperatorsMeshInformed2022} and continuous convolutional filters \cite{coscia2022continuous,DohertyQuadConvQuadratureBasedConvolutions2023} suggest to re-think the neural network structures to incorporate geometrical information, and keep an inductive bias which is physically consistent. 

Geometric deep learning \cite{bronstein2021geometric} arises as a unifying theory to treat data exploiting information on geometry \cite{battaglia2018relational}.
For example, an active line of research is devoted to generalizing convolutions to graphs and manifolds \cite{kipf2016semi,monti2017geometric}.
Here, we seek to preserve the underlying geometric structure of the data, known as \textit{geometric priors}. The imposition of these provides a deeper understanding and a more accurate treatment of the dataset. Moreover, the additional knowledge included in the algorithm can compensate when dealing with scarce data, without drastically compromising the generalization properties.

This work is dedicated to the study of deep learning in model order reduction for parametric PDEs, developing and analyzing a data-driven and non-intrusive framework to deal with models defined on unstructured grids, capturing their geometric features through a combination of Convolutional Autoencoders and Graph Neural Networks (GNNs).

We take advantage of the straightforward interpretation of a field defined on a general mesh as a simple, undirected and connected graph, with the field as a feature. 

\begin{figure}
\begin{center}
\includegraphics[width=\textwidth]{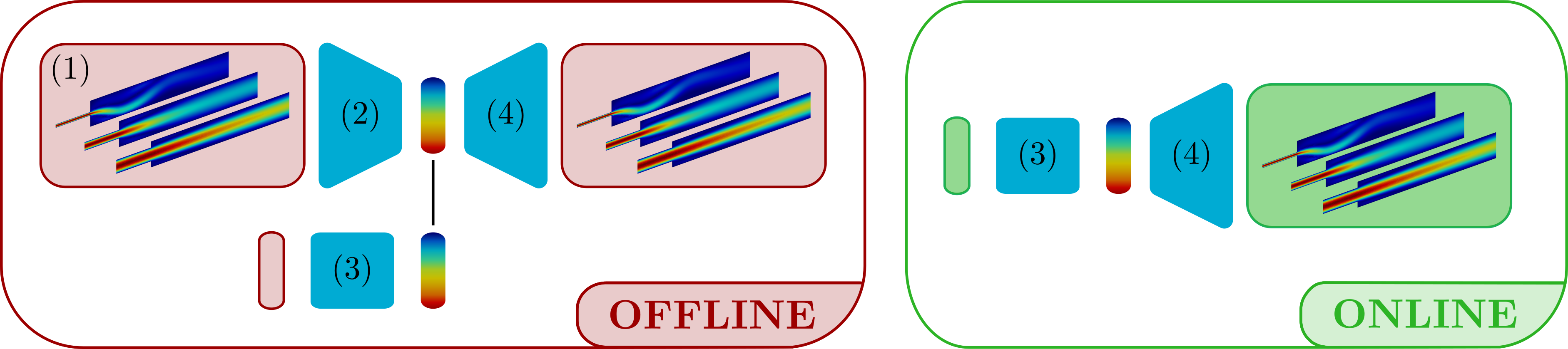}
\end{center}
\caption{{Schematic representation of the GCA-ROM architecture. In the offline phase, a graph autoencoder is trained via unsupervised learning, while a multilayer perceptron learns the encoded representation through a mapping from the parametric space to the low dimensional space. In the online phase, the bottleneck map is evaluated, and the decoder is used to reconstruct the field. The numbers in the figure correspond to the indices referenced in its description.}}
\label{fig:scheme}
\end{figure}

Inspired by \cite{lee2020model,fresca2021comprehensive}, we propose a modular architecture, namely Graph Convolutional Autoencoder for Reduced Order Modelling (GCA-ROM), see Figure \ref{fig:scheme}, which subsequently exploits:
\begin{itemize} 
    \item[(1)] a graph-based layer to express an unstructured dataset;
    \item[(2)] an encoder module compressing the information through:
    \begin{itemize} 
        \item[(i)] spatial convolutional layers based on MoNet \cite{monti2017geometric} to identify patterns between geometrically close regions;
        \item[(ii)] skip-connection operation, to keep track of the original information and help the learning procedure;
        \item[(iii)] a pooling operation, to down-sample the data to obtain smaller networks;
    \end{itemize}
    \item[(3)] a bottleneck map, connected to the encoder's output by means of a dense layer, which learns the latent behavior in a vector;
    \item[(4)] a decoder module,  {connected to the latent vector through a dense layer}, recovering the original data by applying the same operations as in the encoder, but in reverse order.
\end{itemize}

Using this approach, we study problems covering a wide spectrum of characteristics. Mainly, we focus on both linear and nonlinear parametrized PDEs with scalar and vector unknown fields, and analyze different configurations of the proposed architecture, to evaluate their advantages and drawbacks for several benchmarks in model order reduction. The methodology has been tested on a complex phenomenon exhibiting a bifurcating behavior (linked to the non-uniqueness of the solution). In this case, we exploit the bottleneck and a clustering technique to classify bifurcating regimes, and thus detect the bifurcation points. 

We acknowledge the use of the open source software FEniCS \cite{fenics} and RBniCS \cite{rbnics} to generate the dataset, and PyTorch Geometric \cite{FeyLenssen2019} for the implementation and training of the graph convolutional autoencoder. The codes to reproduce the results and create a new test dataset from high-fidelity data are publicly available in the GitHub repository: \url{https://github.com/fpichi/gca-rom}.

This paper is structured as follows. After a brief overview of the GNNs literature, Section 2 is dedicated to the description of the model order reduction setting, differentiating intrusive projection-based methodologies and non-intrusive data-driven techniques. In Section 3, we detail the basic concepts of GNNs, and present our GCA-ROM architecture detailing all its components. Section 4 discusses the results obtained for different test cases. Finally, we state the conclusions and future research lines opening from the proposed methodology in Section 5.

\subsection{Literature review on {Graph Neural Networks}}
In recent years, Graph Neural Networks \cite{ZhouGraphNeuralNetworks2020,HamiltonGraphRepresentation2020} have been an active area of research with rapid growth. GNNs are a specific kind of deep learning architectures capable of extracting information from a graph-structured dataset, including geometric configurations, relationships between nodes and connectivity, and features' behavior.
They are used in different fields such as recommendation systems, social networks and to express molecular structures, for their ability to process non-euclidean data. Several tasks can be performed through GNNs 
\begin{itemize}
    \item[$\circ$] at  graph-level: where the goal is the property/label for the whole structure,
    \item[$\circ$] at node-level: to infer the features at some specific location,
    \item[$\circ$] at edge-level: to predict the existence of a given link/connection.
\end{itemize}
The information is usually propagated between nodes through two of the basic operations of GNNs: aggregation and message passing.  

A key property of GNNs is that they preserve discrete symmetries, i.e.\ nodes do not have to be in any specific order when they are provided as input to the network. Consequently, the functions are permutation-invariant. This is a great advantage that allows generalizing common architectures such as Convolutional Neural Networks that require a fixed template/filter. Graph networks have been widely employed in recent research \cite{wu2020comprehensive}. Notably, it has been especially successful in large deformation settings and complex simulations, drug discovery, and in general in scenarios where  data is unstructured. 

Focusing on our area of interest, graphs have been applied to the study of computational mechanics for the reconstruction of solution fields \cite{xu2021convolutional}, the integration of complex systems in time \cite{sanchez2020learning}, and combined with learning biases to preserve fundamental physical laws \cite{hernandez2022thermodynamics}.

In the model order reduction context, we can discern between two classes. On the one hand, we find those that reduce the mesh to alleviate graph operations \cite{fortunato2022multiscale}, and its edge-augmented version \cite{GladstoneGNNbasedPhysicsSolver2023}. The work by Han et al.\ \cite{han2022predicting} fall into first class and focus on preserving the graph information by finding a coarser representation by collapsing nodes to speed up the simulation. 

On the other hand, the second class of reduction methods is related to the classical ROM setting, with which a vector low-dimensional representation of the data is found. For example, in \cite{wang2019dynamic} the main goal is to deal with 3D representations for reconstruction and tracking. Graph U-nets \cite{gao2019graph} present an architecture of graph convolutions combined with down-sampling layers that lead to a low dimensional manifold. Ranjan et al.\ \cite{ranjan2018generating} introduce the so-called COMA network consisting of a combination of spectral Chebyshev convolutions \cite{defferrard2016convolutional} and down-sampling operations to generate new faces. In \cite{zhou2020fully}, the authors define the pooling operations based on varying kernels formulated at the mesh points. Also, the work of Kashefi et al.\ \cite{kashefi2021point} presents the reconstruction of solution fields for a given geometry. 
In the context of parametrized PDEs, we acknowledge the recent work of Gruber et al.\ \cite{gruber2022comparison}, which compared different ML approaches for data-driven ROMs, considering spectral convolutions without exploiting down-sampling procedures.

Indeed, the definition of down-sampling and up-sampling operations to coarsen and reconstruct domains is one of the most critical steps for an autoencoder-based architecture. These methods are widely studied for classification problems \cite{cangea2018towards}. However, pooling, un-pooling, and inverse convolutions are not defined for non-euclidean domains. Although there are functions that perform an approximation of these operations \cite{qin2020uniform}, their application is still ongoing research. Currently, the most widely used reconstruction function for up-sampling is the one exploiting $k$-nearest neighbors \cite{qi2017pointnet++}. Finally, we remark that since GNNs are still in their early stage of development, they offer much room for further investigations and developments.

\section{{ Reduced Order Models for parametrized Partial Differential Equations}}

In this section we review the standard setting of model order reduction in the context of parametrized partial differential equations. 
In particular, we will outline the main differences between classical approaches based on intrusive projections, as opposed to non-intrusive data-driven techniques. 

\subsection{Intrusive projection-based methodologies}
In the field of computational mechanics, the numerical approximation of physical models described by PDEs is of utmost importance. Such models are usually characterized by the domain in which they are posed, geometric and/or physical parameters, and possibly nonlinear terms.
The so-called high-fidelity techniques, such as the Finite Element method, involve a large number of degrees of freedom, preventing a fast resolution of the model under investigation. This issue is especially harmful in the parametrized context, when one is interested in recovering the solution for different physical/geometrical configurations, for many-query or real time investigations.

Reduced Order Models have been developed during the last decade to address this computational bottleneck.
In particular, the Reduced Basis (RB) method \cite{hesthavenCertifiedReducedBasis2015,quarteroniReducedBasisMethods2016} aims at constructing a low-dimensional space which captures the key features of the system. This way, one uses a Galerkin projection to project the model w.r.t.\ the basis spanning the low-dimensional manifold, to obtain a reduced system. 

Given a suitable Hilbert space $\mathbb{X}_{\bmu} \doteq \mathbb{X}(\Omega(\bmu))$, we consider the nonlinear, steady, and parametrized problem: given $\boldsymbol{\mu} \in \mathcal{P} \subset \mathbb{R}^P$, seek $u(\bmu) \in \mathbb{X}_{\bmu}$ such that
\begin{equation*}
\label{eq:strong_form}
F(u(\boldsymbol{\mu}); \boldsymbol{\mu}) = 0 \qquad \text{in}\quad \mathbb{X}'_{\bmu}
\end{equation*}
where $\boldsymbol{\mu} \in \mathcal{P}$ is the parameter, $\Omega(\bmu)$ the computational domain (possibly $\bmu$-dependent), and $F$ expresses the PDE operator.
A possible high-fidelity approach consists in a Galerkin projection of the associated weak-formulation over some finite dimensional subspace $\mathbb{X}_{\bmu}^\mathcal{N} \subset \mathbb{X}_{\bmu}$ of dimension $\mathcal{N}$.

Thanks to the Newton-Kantorovich method \cite{ciarlet2013linear}, the formulation of the linearized system {at the $k$-th iteration} reads: given $\bmu \in \mathcal{P}$, seek $\delta \textbf{u}_\mathcal{N} \in \mathbb{R}^{\mathcal{N}}$ such that 
\begin{equation}
\label{eq:newton_hf}
\mathsf{J}_\N(\textbf{u}_{\mathcal{N}}^{k}; \bmu)\delta \textbf{u}_{\mathcal{N}} = -\mathsf{F}_{\mathcal{N}}(\textbf{u}_{\mathcal{N}}^{k}; \bmu) \qquad \text{in}\ \mathbb{R}^{{\mathcal{N}}}
\end{equation}
where $\mathsf{F}_\N \in \R^\N$ and $\mathsf{J}_\N \in \R^{\N\times\N}$ are the high-fidelity residual vector and the Jacobian matrix, respectively. Finally, we update the iteration as $\textbf{u}_{\mathcal{N}}^{k} = \textbf{u}_{\mathcal{N}}^{k-1} + \delta \textbf{u}_{\mathcal{N}}$ until convergence.

Exploiting the same concept in the ROM context, one builds a basis $\{\boldsymbol{\zeta}_k\}_{k=1}^{N} \subset \mathbb{X}_{\bmu}^\mathcal{N}$ for the reduced space $\mathbb{X}_{\bmu}^N$, e.g.\ by means of the POD technique \cite{volkwein2011model}, to project the system once again, obtaining the following formulation: given $\bmu \in \mathcal{P}$, seek $\delta \textbf{u}_N \in \mathbb{R}^{N}$ such that 
\begin{equation}
\label{eq:newton_rb}
\mathsf{V}^T\mathsf{J}_\N(\mathsf{V} \textbf{u}_N^{k}; \boldsymbol{\mu})\mathsf{V} \delta \textbf{u}_N = -\mathsf{V}^T\mathsf{F}_\N(\mathsf{V} \textbf{u}_N^{k}; \boldsymbol{\mu}) \qquad \text{in}\quad \mathbb{R}^{N}
\end{equation}
where the basis matrix is encoded in $\mathsf{V} = [\boldsymbol{\zeta}_1 | \dots | \boldsymbol{\zeta}_N] \in \mathbb{R}^{\N \times N}$, and we update the iteration as $\textbf{u}_N^{k} = \textbf{u}_N^{k-1} + \delta \textbf{u}_N$.
The reduced basis approximation is finally given by, $u_N(\bmu) = \sum_{k=1}^{N} u_N^{(k)}(\bmu) \boldsymbol{\zeta}_k$, where $\textbf{u}_N(\bmu) = \{u_N^{(k)}(\bmu)\}_{k=1}^{N}$ is the reduced coefficient vector.

The POD basis can be extracted by exploiting the FE data, i.e.\ a set of solutions corresponding to different parameters $\{\bmu^{i}\}_{i=1}^{N_{\text{S}}}$, the so-called snapshots $\{u_\N(\bmu^{i})\}_{i=1}^{N_{\text{S}}}$.
The efficiency of the classical approach relies on the affine decomposition assumption, and on the independence of the online computations from the degrees of freedom $\N$.
Unfortunately, these assumptions are not fulfilled in the general nonlinear context, since \eqref{eq:newton_rb} requires repeated assembly of the system involving the basis matrix $\mathsf{V}$, drastically compromising the performance. The Empirical Interpolation Method and its variants \cite{barrault04:_empir_inter_method,Chaturantabut2010} are affine-recovery techniques which interpolate the non-affine or nonlinear terms, suggesting a considerable speed-up. Such methodologies entail an intrusive approach, requiring the a priori knowledge of the physics and its high-fidelity operators.

\subsection{Non-intrusive and data-driven approaches}
As detailed above, the intrusive nature of projection-based model order reduction can jeopardize the goal of traditional reduced order models, i.e.\ high computational efficiency. Non-intrusive and data-driven ROMs arise as an alternative to these techniques, with the scope of reducing the computational complexity without the need to project the governing equation. The starting point to build such ROMs is the dataset constituted by the aforementioned snapshots.
Several linear and nonlinear methodologies have been developed and analyzed, occasionally benefiting by machine learning approaches.

The key reference of the first class is represented by the RB method, expressing the reduced approximation as the linear expansion over the basis functions as $u_{\N}(\bmu) \approx u_N(\bmu) = \mathsf{V} \textbf{u}_N(\bmu)$.
While traditional methods require the solution of a reduced (nonlinear) system to obtain the reduced coefficient vector $\textbf{u}_N(\bmu)$, non-intrusive approaches exploit different workarounds to replace this step.

Proper Orthogonal Decomposition with interpolation (PODI) \cite{Bui-ThanhProperOrthogonalDecomposition2003,DemoNonintrusiveApproachReconstruction2019} projects the snapshots onto the POD basis and then interpolates the reduced coefficients. When interested in regression rather than interpolation, machine learning techniques become relevant. With POD-NN, one recovers these by exploiting a feedforward neural network for a supervised learning task, where the dataset is given by the projected snapshots \cite{hesthaven2018non,PichiArtificialNeuralNetwork2023,BarnettNeuralNetworkAugmentedProjectionBasedModel2022}. Physics-reinforced neural network (PRNN) addresses the low-data regime, augmenting the POD-NN loss with the reduced equations encoding the physics \cite{ChenPhysicsinformedMachineLearning2021}. 

The aforementioned methodologies help to recover an efficient reconstruction of the solution, but they are still exploiting a linear basis expansion. Unfortunately, many problems cannot be effectively tackled by linear ROMs, such as (i) advection-dominated problems for which a large number of modes is required \cite{GreifDecayKolmogorovNwidth2019}, and (ii) bifurcation problems, innately of nonlinear and ill-posed nature \cite{pichi_phd}.

Nonlinear reduction techniques have been developed to cope with these issues. These involve an approximation of the form $\textbf{u}_{\N}(\bmu) \approx \psi(\textbf{u}_N(\bmu))$, where $\psi$ is a suitable nonlinear map. The linear approach corresponds to $\psi(\textbf{u}_N(\bmu)) = \mathsf{V} \textbf{u}_N(\bmu)$. As alternatives, we highlight: the registration method \cite{TaddeiRegistrationMethodModel2020}, kernel principal component analysis \cite{DiezNonlinearDimensionalityReduction2021}, the method of freezing \cite{OhlbergerNonlinearReducedBasis2013}, the shifted POD \cite{ReissShiftedProperOrthogonal2018}, localized models \cite{AmsallemNonlinearModelOrder2012} and the use of the Wasserstein space \cite{EhrlacherNonlinearModelReduction2020}.

These methodologies allow to obtain much better accuracy, but are usually characterized by a less trivial online phase.
This is where machine learning techniques come into play. As discussed in the introduction, the autoencoder architecture, with its nonlinear compression and decoding phases, represents the ideal generalization of the POD technique.

Recently, they have been widely used to develop nonlinear reduced order models. In \cite{lee2020model}, the authors present a nonlinear manifold least-squares Petrov-Galerkin method based on convolutional autoencoders. A hyper-reduced extension exploiting the physics have been introduced in \cite{RomorNonlinearManifoldReducedOrder2023}.
The DL-ROM approach developed in \cite{fresca2021comprehensive} augments the training phase with a supervised task, in which the bottleneck is learned through a feedforward neural network. A POD preprocessing have been considered to reduce the size of the network in \cite{fresca2022pod}. An interesting extension entails the reconstruction of the hyper-reduction operators \cite{CicciDeepHyROMnetDeepLearningBased2022}. 
Dynamic integration based on thermodynamics priors exploits autoencoders to learn mappings to a low-dimensional manifold \cite{hernandez2021deep, moya2022physics}.
Many other non-intrusive approaches have been proposed recently, e.g.\ based on Gaussian Process Regression (GPR) \cite{GuoReducedOrderModeling2018} and Operator Inference \cite{PeherstorferDatadrivenOperatorInference2016}. 
Finally, Neural Operators are successful in approximating mappings between function spaces, rather than functions themselves. These methodologies usually benefit from both the discretization-invariant and the universal approximation properties \cite{kovachki2021neural,LuLearningNonlinearOperators2021a,DemoDeepONetMultiFidelityApproach2023}.

Following the path set by these works, graph convolutional autoencoders are the natural next step for the development of nonlinear model order reduction techniques, combining the classical PDE framework with solutions defined on unstructured grids and machine learning structures encoding geometric features. 

\section{Graph convolutional autoencoder for nonlinear model order reduction}

Computational problems derived from physical models are naturally formulated over discretized domains for their analysis and resolution. In many cases, the geometry is not simple, suggesting the use of unstructured mesh. Deep learning techniques typically operate with structured datasets, e.g.\ standard CNNs, that depend strongly on the ordering of data points. Treating complex domain using such approaches often results in inefficient trainings, as they do not take the geometrical information into consideration. 

Geometric Deep Learning arose as a mathematical framework to augment neural networks capabilities with geometric priors. 
This results in data driven approaches, characterized by both the spatial domain and the physical behavior of the system.

Graph Neural Networks are a branch of geometric deep learning. In this context, GNNs interpret the computational mesh as a peculiar graph structure, where the nodes are the vertices, and the edges constitute the boundary of each element of the triangulation. 
This representation makes the use of graph neural networks suitable for the development of computational problems defined on unstructured grids.

\subsection{Brief introduction to graph neural networks}
 
\begin{figure}
	\begin{minipage}{0.3\textwidth}
		\includegraphics[width=\textwidth]{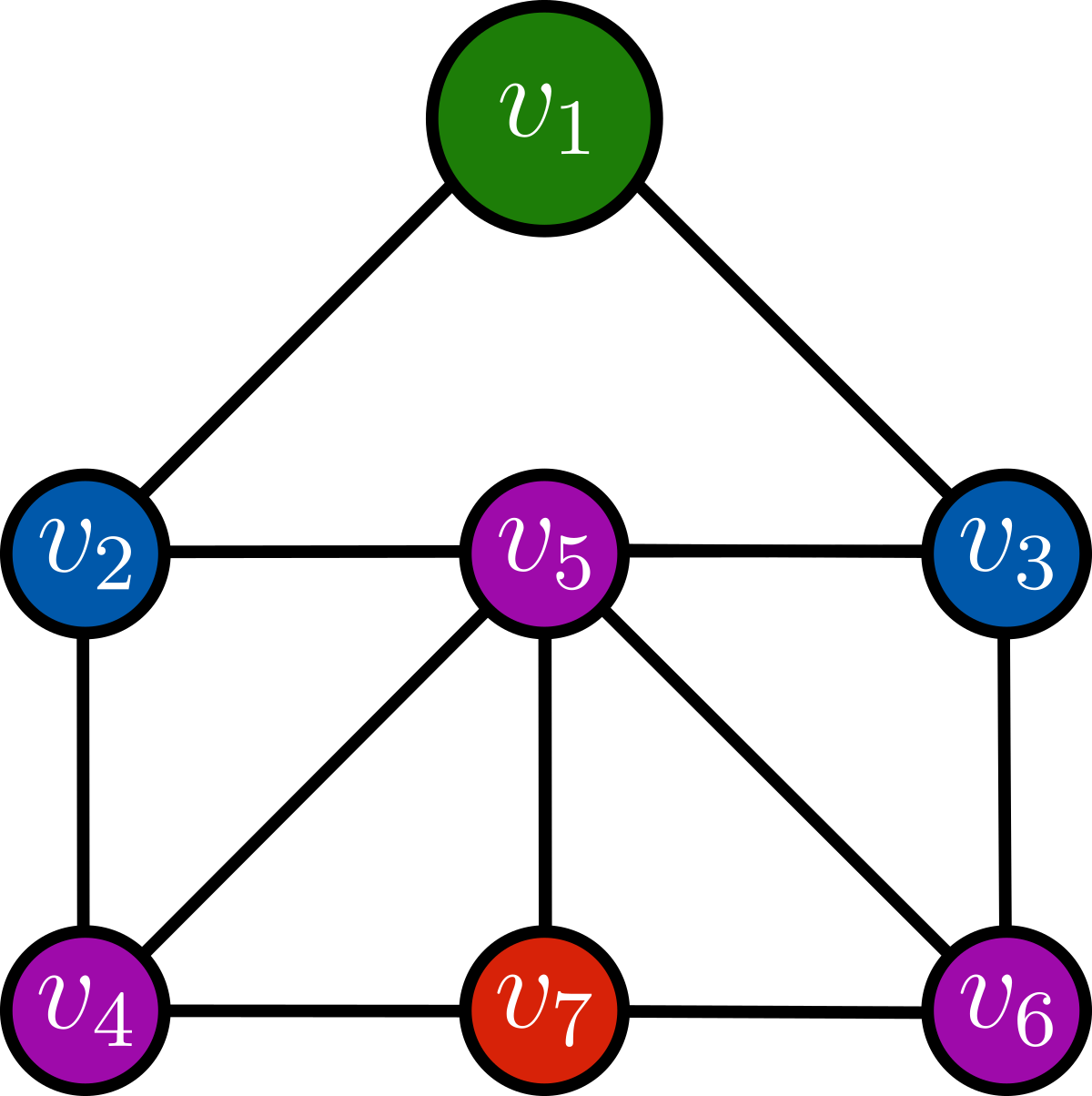}
	\end{minipage}\hfill
	\begin{minipage}{0.65\textwidth}
		\includegraphics[width=\textwidth]{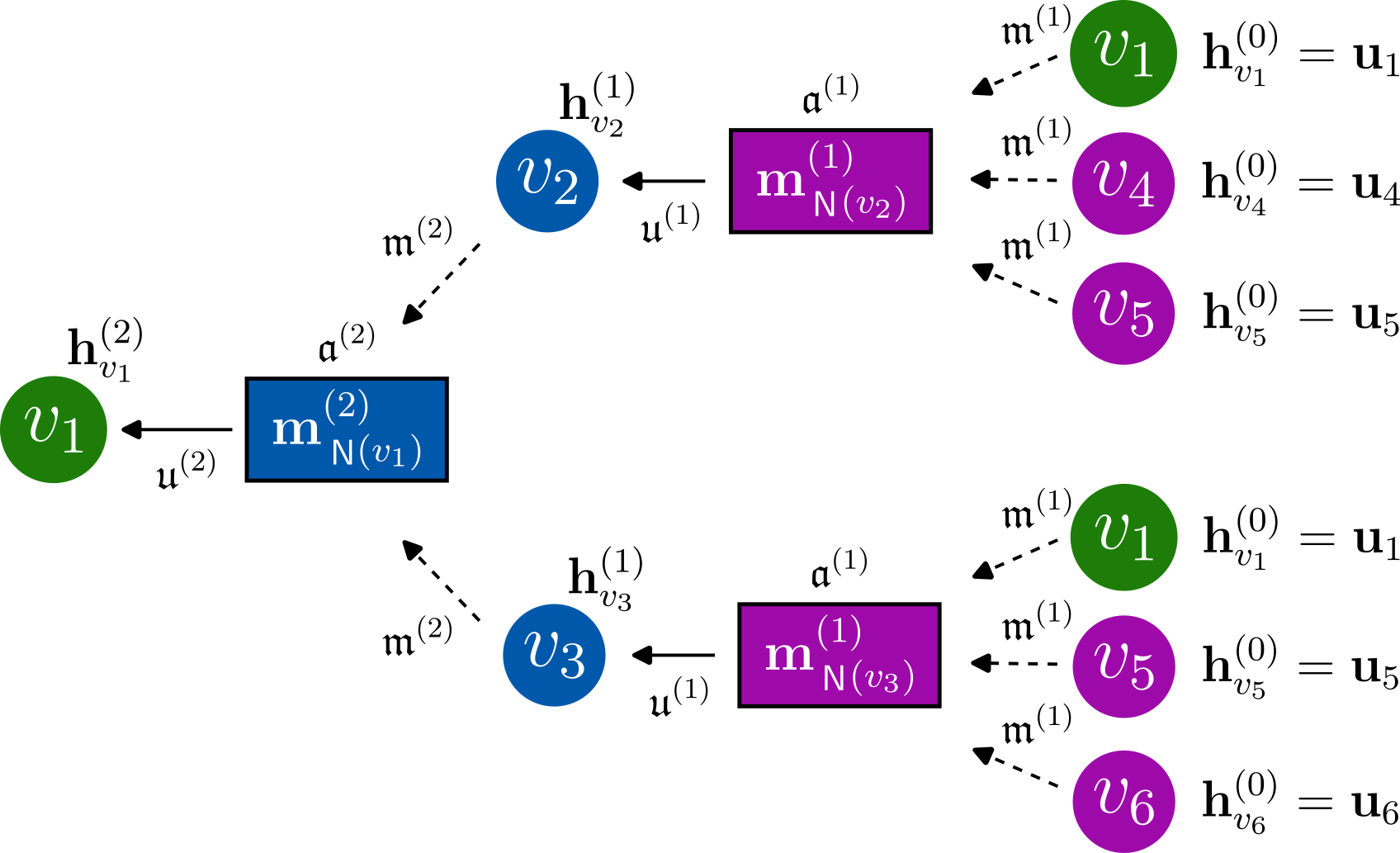}
	\end{minipage}
	\caption{{Example of a graph with message passing and aggregation procedures, showing the flow of information to compute a node's hidden embedding, where colors denote the connectivity-based distance from node $v_1$ to be updated.}}
	\label{fig:graph_basic_operations}
\end{figure}

Let us consider the mesh $T(\mathcal{V},\mathcal{E},\mathcal{F})$, which can be seen as a simple, undirected, and connected graph $\mathcal{G}(\mathcal{V} , \mathcal{E})$ with nodes $\mathcal{V}$ and edges $\mathcal{E}$  (see Figure \ref{fig:graph_basic_operations}), including additional information about the elements of the mesh $\mathcal{F}$.

The nodes have associated features $\mathbf{u}$, that represent the evaluation of a set of state variables at the vertices of the mesh. This information can be embedded in the matrix $\mathbf{U} = \left[\mathbf{u}_1 | \dots | \mathbf{u}_{N_h}\right]^T \in \mathbb{R}^{N_h \times d}$, for all the $N_h$ nodes of the graph\footnote{The dataset is obtained by extracting the values of the high-fidelity solutions at the corresponding vertexes of the mesh, regardless the polynomial order of the Lagrange elements.}, which have no particular order, and $d$ features, i.e.\ a scalar or vector field. The adjacency matrix $\mathbf{A} \in \mathbb{R}^{N_h\times N_h}$ records of the connections among nodes. A convenient way of expressing it is by means of the adjacency list, where the $k$-th entry, corresponding to $e_k \in \mathcal{E}$, is the pair $(i, j)$ denoting the existence of a link between $v_i, v_j \in \mathcal{V}$.
This is crucial to efficiently define operations which are permutation-invariant and equivariant. Thus, one can randomly permute the original ordering of the nodes induced by the mesh labelling without changing the final output. This is a key difference with standard CNNs, where a fixed filter produces different outputs if two pixels are swapped.

Starting from the graph structure as defined above, one obtains a graph neural network by defining a set of optimizable operations acting on all attributes of the graph. Despite their recent introduction, many works have focused on how to process graph-data through a neural network. 
Thus, let us review in more detail the basic operations which constitutes the backbone of many GNN architectures: (i) message passing framework, (ii) convolutional layers, and (iii) down-sampling and up-sampling procedures.

\subsubsection{Message passing framework}
The main idea is to propagate the information to local neighborhoods of each node $v$, which is commonly denoted by $\mathsf{N}(v)$ and defines its computation graph with degree $|\mathsf{N}(v)|$. Such messages are exchanged between nodes at different $k$-hops (or layers) of the graph, and are updated by exploiting neural networks.
At the $k$-th hop, one computes the hidden embedding $\mathbf{h}^{(k)}_u \in \mathbb{R}^{d(k)}$ of the node $u$, a vector representing a transformation of its original features by means of differentiable aggregate and update computations.
In practice, at a node $u \in \mathcal{V}$, one assembles the messages to be sent through the operation $\mathfrak{m}^{(k)}$ as $$\mathbf{m}_v^{(k)} = \mathfrak{m}^{(k)}(\mathbf{h}_v^{(k-1)}), \quad \text{from each node} \ v \in \mathsf{N}(u)$$ aggregates them with $\mathfrak{a}^{(k)}$ in $$\mathbf{m}_{\mathsf{N}(u)}^{(k)} = \mathfrak{a}^{(k)}(\{\mathbf{m}_v^{(k)},\ \forall v \in \mathsf{N}(u)\}),$$ and finally updates the hidden embedding by means of the function $\mathfrak{u}^{(k)}$  $$\mathbf{h}_u^{(k)} = \mathfrak{u}^{(k)}(\mathbf{m}_{\mathsf{N}(u)}^{(k)}).$$
For each node $v_j \in \mathcal{V}$, the initialization of the hidden embedding are defined as its own input features, i.e.\ $\mathbf{h}_{v_j}^{(0)} = \mathbf{u}_j$. Since we seek to update the hidden embedding without erasing previous information, we consider ghost self-edges, such that $\mathsf{N}(u)$ contains the node $u$ itself.

The whole procedure can be performed in several ways. The simplest example of message function $\mathfrak{m}^{(k)}$ is the multiplication of the hidden embedding with a weight matrix $\mathbf{W}^{(k)} \in \mathbb{R}^{d^{(k)} \times d^{(k-1)}}$, i.e.\  $\mathbf{m}_v^{(k)} = \mathbf{W}^{(k)} \mathbf{h}_v^{(k-1)}$.
For the aggregate function $\mathfrak{a}^{(k)}$, one usually considers the normalized sum of the neighbor embeddings given by $\mathbf{m}_{\mathsf{N}(u)}^{(k)} = \sum_{v \in \mathsf{N}(u)}\frac{\mathbf{m}_v^{(k)}}{|\mathsf{N}(u)|}$, to better balance nodes with much higher degree. Finally, the update operation $\mathfrak{u}^{(k)}$ can be seen as the nonlinear activation function $\mathbf{h}_u^{(k)} = \sigma(\mathbf{m}_{\mathsf{N}(u)}^{(k)})$ in the neural network context, e.g.\  with $\sigma(x) = \text{ReLU}(x)$ or $\sigma(x) = \tanh(x)$. 

To summarize, a basic GNN with $K$ layers can be described by the iteration to update the node embeddings\footnote{We remark that with the self-edges notation, in this explanatory example one cannot distinguish between information coming from the node itself and its neighbors. Sharing the same weight matrix $\mathbf{W}^{(k)}$ possibly compromising the capability of the GNN. A common approach is to consider different weight matrices for the two types of information.} 
\begin{equation}
	\mathbf{h}_u^{(k)} = \sigma\left(\frac{1}{{|\mathsf{N}(u)|}}\sum_{v \in \mathsf{N}(u)}\mathbf{W}^{(k)} \mathbf{h}_v^{(k-1)}\right) \quad  \text{for } k = 1, \dots, K, \ \text{and } u \in \mathcal{V}.
	\label{eq:node_embeddings}
\end{equation}
A bias term $\mathbf{b} \in \mathbb{R}^{d(k)}$, omitted here for simplicity, improves the performance.
Mimicking the standard neural networks' framework, one can also rewrite the message passing iteration \eqref{eq:node_embeddings} with the following compact graph-level notation 

\begin{equation}
	\mathbf{H}^{(k)} = \sigma\left(\mathbf{D}^{-1}(\mathbf{A} + \mathbf{I})\mathbf{H}^{(k-1)}{\mathbf{W}^{(k)}}^T\right) \quad  \text{for } k = 1, \dots, K,
	\label{eq:graph_node_embeddings}
\end{equation}
where $\mathbf{H}^{(k)} = \left[\mathbf{h}_{v_1}^{(k)} | \dots | \mathbf{h}_{v_{N_h}}^{(k)}\right]^T \in \mathbb{R}^{N_h \times d^{(k)}}$ expresses the matrix with the hidden embedding for each node taken as row, $\mathbf{D}^{-1}$ is the diagonal matrix with the inverse of the degree of each node $|\mathsf{N}(v)|$, and $\mathbf{A}, \mathbf{I}$ are the adjacency and the identity matrices (representing the self-edges), respectively.

These operations are permutation invariant/equivariant, thus allowing for a message passing framework to produce the same results independent of the initial ordering of the nodes, and based solely on their connections.

\subsubsection{Convolutional layers in the non-Euclidean setting}

A graph convolutional network (GCN) is a specific type of GNN aiming at extending the concept of CNNs, operating in regular Euclidean domains, to handle non-grid data. The GCN learns to combine the hidden embeddings by defining convolutional operations able to capture the relationships between the nodes, and optimize a given loss function. 
This process is similar to the standard convolutional layers in CNNs, where a fixed filter slides over the pixels of an image to produce gathered information. The main issue with CNNs is that the operations are not invariant w.r.t.\ the nodes' order. Given also the different cardinality of each node's neighbors, such extension to graphs is not trivial \cite{kipf2016semi,defferrard2016convolutional}.

One can distinguish between two classes of convolutional layers: spectral and spatial ones.
Spectral approaches are characterized by the spectral analysis of the graph Laplacian matrix to obtain a suitable basis. These methods have strong mathematical foundations from signal processing, but they depend on the graph's cardinality $|N_h|$, requiring an eigen-decomposition to define the filters, which compromises their efficiency. Such procedures define a global filter that can be applied to the whole graph.

Spectral convolutional methods seek to build the filters as $g_\mathbf{w} * \mathbf{H} \doteq \mathbf{U} g_\mathbf{w} \mathbf{U}^T \mathbf{H}$, where $g_\mathbf{w}$ is the filter, and $\mathbf{U}$ the eigenvector matrix of the graph Laplacian $\mathbf{L} = \mathbf{D} - \mathbf{A}$. While ChebNet exploits a Chebyshev polynomial of degree $K$ in $\mathbf{L}$ to approximate the convolution  $g_\mathbf{W} * \mathbf{H}$, GCN simplifies this by considering only the first term \cite{kipf2016semi}.

Conversely, spatial convolutions reflect an ordinary CNN approach, where each node is related to its local neighbors. Being local rather than global, they do not depend on the graph's dimensionality, resulting in faster and more efficient methods.

Spatial methods are mainly based on the message-passing operation, which can be defined in several ways. In GraphSage \cite{hamilton2017inductive} features are sampled and then aggregated, while in GINN \cite{BerroneGraphInformedNeuralNetworks2022} the adjacency matrix is exploited for regression tasks.
In this work, we have decided to exploit the MoNet spatial framework introduced in \cite{monti2017geometric}, which can be interpreted as a Gaussian Mixture Model (GMM), and thus a general class for convolutions in non-Euclidean domains. 
MoNet builds a set of pseudo-coordinates $\mathbf{e}$ used to define the weights of an optimizable Gaussian kernel with $Q$ filters, through the iteration procedure
\begin{equation}
	\mathbf{h}_u = \frac{1}{{|\mathsf{N}(u)|}}\sum_{v \in \mathsf{N}(u)} \frac{1}{Q}\sum_{q=1}^{Q}\boldsymbol{\omega}^{q}(\mathbf{e}_u) \odot {\mathbf{W}^q} \mathbf{h}_v,
	\label{eq:node_monet}
\end{equation}
{where $\odot$ is the element-wise multiplication, and $\boldsymbol{\omega}^{q}$ is the weighting function defined in terms of a trainable mean vector $\boldsymbol{\mu}_q$ and a diagonal covariance matrix $\boldsymbol{\Sigma}_q$ as} $$\boldsymbol{\omega}^{q}(\mathbf{e}_u)=\text{exp}\left(-\frac{1}{2}(\mathbf{e}_u-\boldsymbol{\mu}_q)^T \boldsymbol{\Sigma}_q^{-1} (\mathbf{e}_u-\boldsymbol{\mu}_q)\right).$$
In practice, MoNet considers as pseudo-coordinate the edge attributes given by the distance between two connected nodes, introducing a geometric bias in the learning process. These geometrically-informed convolutions could increase the cost of the learning task w.r.t.\ other convolutional layers, but they also improve the generalizability. 

To conclude, we recall that GCN can be seen as a spatial method, and thus a particular instance of MoNet with $Q = 1$ and $\boldsymbol{\omega}^{1} = 1/|\mathsf{N}(u)|$. In this more general setting the filters $\boldsymbol{\omega}^{q}$ are optimizable.

\begin{figure}
	\includegraphics[width=0.4\textwidth]{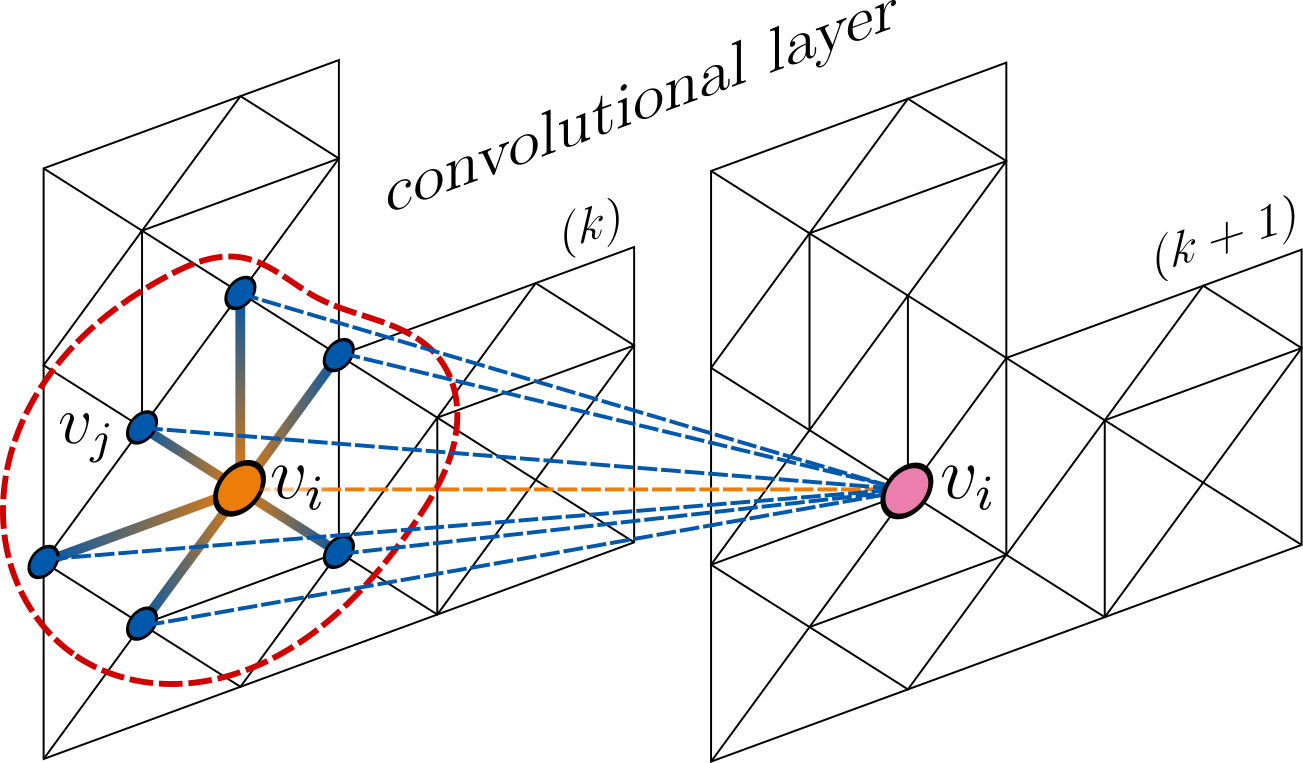}\hfill
	\includegraphics[width=0.4\textwidth]{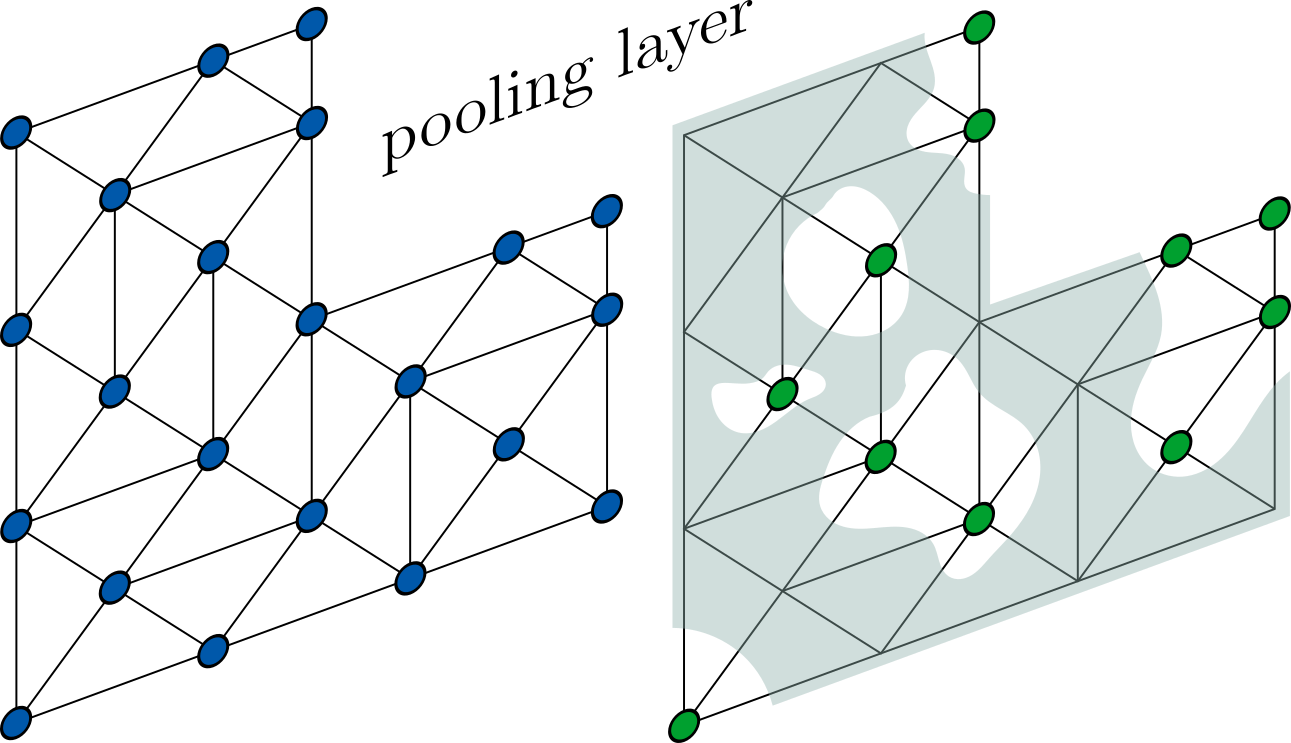}
	\caption{Convolutional and pooling layers.}
	\label{fig:gca_conv_pool}
\end{figure}

\subsubsection{Down-sampling and up-sampling procedures}
Another key difference between CNNs and GNNs consists in the property of the former to naturally reduce the spatial dimension of the input. As we have seen in the previous section, the convolutional layers are defined for each node of the original graph/mesh. For this reason, in some cases, it can be useful to rely on down-sampling and up-sampling.
These methodologies are needed in the encoder and the decoder, to obtain representations of the graph which are, respectively, coarser and finer.

Pooling is the most widely used technique to down-sample the size of the input by aggregating information from multiple nodes and edges. This results in a smaller and more manageable graph, improving generalization and performance. 
Although it is a well-defined operation for CNNs, it cannot be directly applied to graphs, as there is no natural hierarchy in nodes' importance. 
The selection of subsets of nodes can be performed in several ways, from random mask selection to more sophisticated operations including clustering and attention mechanisms \cite{ranjan2018generating,bianchi2020spectral,liu2022graph}. 
Applying down-sampling techniques multiple times in a GNN, one can construct a multiscale hierarchy of coarser domains in the spirit of algebraic multigrid methods \cite{xu2017algebraic}.
However, when misused, these approaches result in information loss, and a deteriorating model performance. 

Un-pooling can be considered the inverse operation of pooling and belongs to the class of up-sampling techniques. While it is relatively simple to collect information to reduce the number of nodes, the reverse is still an open research challenge. Depending on the ML application, up-scaling can be achieved by storing the geometric information obtained during the coarsening step, e.g.\ for classification, segmentation, or noise filtering purposes. COMA stores the barycentric information of the eliminated points \cite{ranjan2018generating}, while Mincut requires the adjacency matrix \cite{bianchi2020spectral}. To our knowledge, there are only a few functions that perform up-sampling without information obtained during the down-sampling step. This is crucial to obtain an online reconstruction fully decoupled from the original solution field, i.e.\ encoder-free. 
PointNet++ proposes a k-NN interpolation of the points to up-sample by considering the position and the features of the nodes in the down-sampled coarser configuration \cite{qi2017pointnet++}.
In particular, given a node at position $\mathbf{x}_i$, we define its feature vector $\mathbf{u}_i$ as the weighted interpolation w.r.t.\ its $k$ neighbors given by
$$\mathbf{u}_i=\frac{\sum_{j=1}^{k} \xi(\mathbf{x}_j)\mathbf{u}_j}{\sum_{j=1}^{k}\xi(\mathbf{x}_j)}, \quad\quad \text{where} \quad \xi(\mathbf{x}_j)=\frac{1}{d(\mathbf{x}_i, \mathbf{x}_j)^2} .$$
We remark that, in this preliminary investigation, the goal of the pooling procedures is restricted to  the study of the learning capabilities from a coarser representation of the mesh, obtained e.g.\ from a fixed set of sensors. In the following we will show that, if the loss of information is not excessive, our naive approach is able to recover the field with acceptable accuracy.
This difficult task is of crucial importance to drastically reduce the number of trainable parameters and speed-up the convergence.
Performing an optimal pooling remains open, with interesting works about multiscale layers \cite{grattarola2022understanding,BarweyMultiscaleGraphNeural2023}, and will be the subject of future investigations.

\subsection{{ A Graph Convolutional Autoencoder approach for Reduced Order Modelling}}
We are now ready to propose the graph convolutional autoencoder for model order reduction applications. The main sources of inspiration is CNN-based autoencoder architecture introduced in \cite{fresca2021comprehensive,lee2020model}.
As detailed in the previous sections, such approaches are particularly useful when dealing with structured meshes, that can be seen as images with a fixed number of neighboring pixels. For this reason, we consider a similar architecture, but extend its applicability in a geometrically consistent manner to unstructured meshes defined over complex domains, when a Cartesian representation is no longer possible.

ML approaches mimic the standard ROM setting, having an offline phase in which one forms the dataset and build the reduced model by training the neural network, and an online phase for the real-time evaluation of the field of interest. 

\begin{figure}[th]
	\centering
	\includegraphics[width=\textwidth]{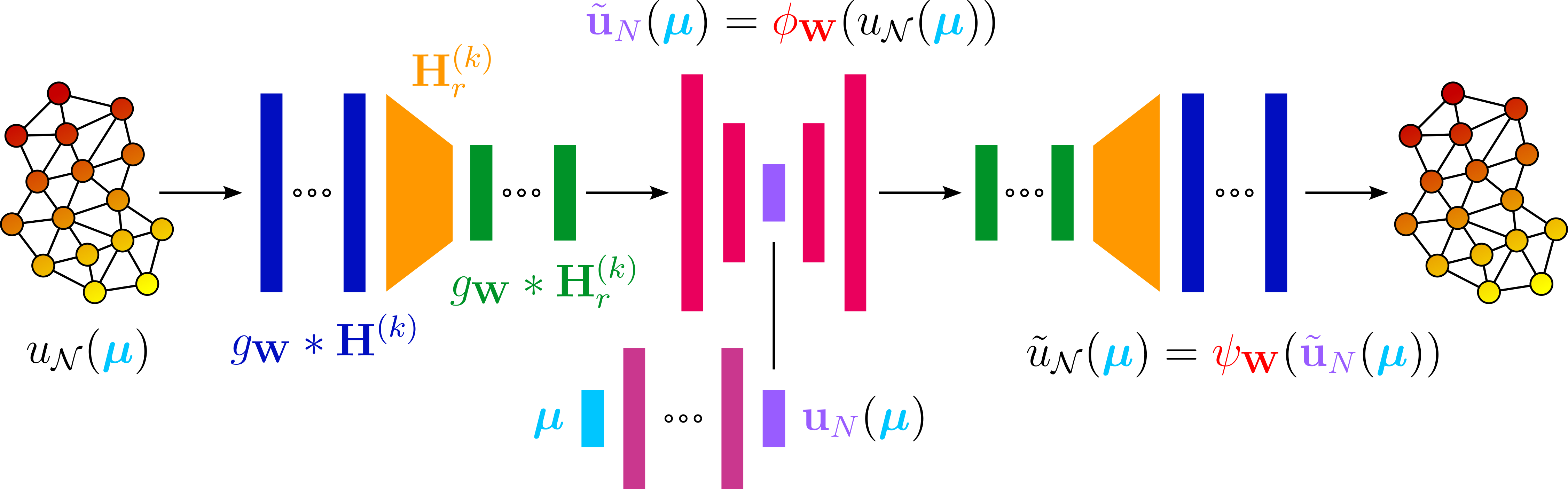}
	\caption{GCA-ROM architecture for the offline phase.}
		\label{fig:gca_rom_offline}
\end{figure}

\begin{figure}[!th]
	\centering
	\includegraphics[width=0.65\textwidth]{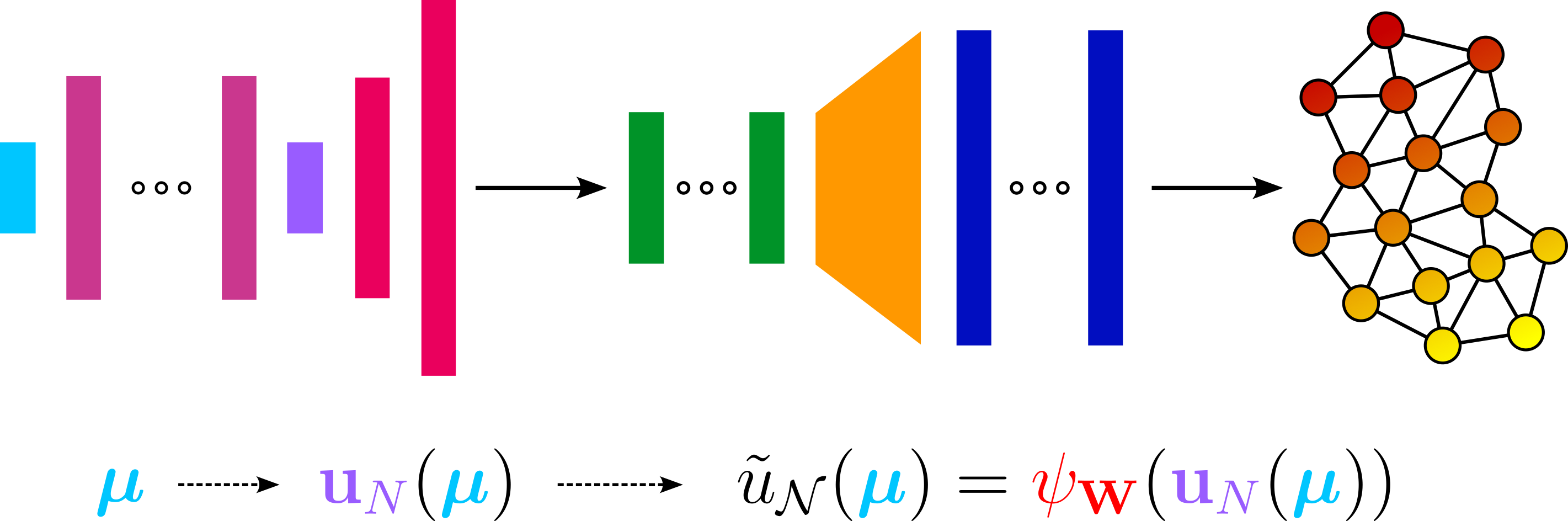}
	\caption{GCA-ROM architecture for the online phase.}
		\label{fig:gca_rom_online}
\end{figure}

We consider the graph dataset $\Xi = \{\textbf{u}_\N(\bmu^{i}), \Omega_{\N}(\bmu^{i})\}_{i=1}^{N_{\text{S}}}$, formed by $N_{\text{S}}$ solutions $\textbf{u}_\N(\bmu^{i})$ of a parametrized PDE defined over unstructured meshes $\Omega_{\N}(\bmu^{i})$, corresponding to the set of parameters $\{\bmu^{i}\}_{i=1}^{N_{\text{S}}}$. 
Following the first DL-ROM approach \cite{fresca2021comprehensive}, we detail the steps of the modular architecture and refer to Figures \ref{fig:gca_rom_offline} and \ref{fig:gca_rom_online}, respectively, for the offline and online phases.

The architecture for the offline training is composed of an autoencoder and a multi-layer perceptron (MLP). The former seeks to approximate the identity map while encoding the information into a low-dimensional space expressed by the bottleneck or latent space.

We feed the encoder module $\phi_\mathbf{W}$ with the graph data $\Xi$, exploiting MoNet \eqref{eq:node_monet} as the message passing algorithm.
In this step, by means of graph convolutions $g_\mathbf{W} * \mathbf{H}^{(k)}$, we extract the most meaningful information between the nodes and their evolution w.r.t.\ the samples.
When interested in reducing the graph dimensionality, one can consider an additional module by applying a fixed random mask for down-sampling, to obtain a reduced graph with $\mathbf{H}_r^{(k)}$ hidden embeddings. A second batch of convolutional layers $g_\mathbf{W} * \mathbf{H}_r^{(k)}$ can be used to further encode the latent space. 
{To improve training, we also augment the network with a residual learning strategy exploiting skip-connections, featuring identity shortcuts inspired by the ResNet framework developed in \cite{he2016deep}. This step, which can only be applied in batches of convolutions sharing the same graph's size, enables the architecture to acquire the residual knowledge between layers, while simultaneously promoting the flow of gradient information throughout the network.}
At this stage, we use fully-connected layers to connect the graph structure with the bottleneck. We thus encode the original information in the latent vector $\tilde{\mathbf{u}}_N(\bmu) \in \mathbb{R}^N$, obtained as the evaluation of the solution field through the encoder structure $\tilde{\mathbf{u}}_N(\bmu) = \phi_\mathbf{W}(\textbf{u}_\N(\bmu))$. The decoding structure is composed by the same operations as the encoding ones, but in reversed order. This way, processing the bottleneck with the decoder $\psi_\mathbf{W}$, we obtain the reconstruction of the original field $\tilde{u}_\N(\bmu) = \psi_\mathbf{W}(\tilde{\mathbf{u}}_N(\bmu))$, and we can assemble the first term in the loss function for the unsupervised learning task as 
\begin{equation*}
	\mathcal{L}_{\text{MSE}} = \frac{1}{N_\text{tr}}\sum_{i=1}^{N_\text{tr}}\left\|\tilde{\textbf{u}}_\N(\bmu^i) - \textbf{u}_\N(\bmu^i) \right\|^2_2.
\end{equation*}
Having encoded the latent variables with the autoencoder, we use this dataset for the supervised learning task with the MLP. Thus, we are non-intrusively building the map $\textbf{u}_N(\bmu) = \text{MLP}(\bmu)$, from which we can recover the reduced coefficient. The MLP, fundamental for the online evaluation phase, defines the loss
\begin{equation*}
	\mathcal{L}_{\text{BTT}} = \frac{1}{N_\text{tr}}\sum_{i=1}^{N_\text{tr}}\left\|\mathbf{u}_N(\bmu^i) - \tilde{\mathbf{u}}_N(\bmu^i) \right\|^2_2.
\end{equation*}
These two Mean Squared Error (MSE) loss terms can be balanced through the hyperparameter $\lambda$, resulting in the loss function
\begin{equation*}
	\mathcal{L} = \mathcal{L}_{\text{MSE}} + \lambda \mathcal{L}_{\text{BTT}},
\end{equation*}
which guides the training procedure through the Adam optimizer.
Finally, during the online phase, the bottleneck can be directly evaluated by exploiting the MLP for a new parameter $\bmu$, and successively decompressed through the graph decoder, to recover the corresponding field defined over its geometry. 

As for the preprocessing of the dataset, we applied a double normalization with respect to the sample and the features, i.e.\ we impose zero mean and unitary variance for each column (i.e.\ node-wise), and subsequently for each row (i.e.\ parameter-wise). This allowed us to treat all the benchmarks within a unique framework, greatly improving the results. Moreover, we split the dataset $\Xi$ into training and testing set w.r.t.\ the percentage rate $r_\text{t}$ in the following manner: $\Xi_\text{tr} = \mathsf{shuffled}({\Xi})[{:}, {:}r_\text{t} N_\text{S} ]$ and $\Xi_\text{te} = \Xi \setminus \Xi_\text{tr}$ of dimension $N_\text{te}$.

Common disadvantages of machine learning approaches are the tuning of the hyperparameter and the learning procedure. As we will show in Section \ref{sec:num_res}, the GCA-ROM architecture is quite robust w.r.t.\ both. Among the possible quantities to be investigated we focus on the following set: the training rate percentage $r_\text{t}$, the pooling rate percentage $r_\text{p}$, the number of nodes in the feedforward network of the GCA $\text{ffn}$, the number of nodes in each MLP layer of the latent map $n_l$, the dimension of the bottleneck $n$, the weight $\lambda$ in the loss function, and the number of a-priori and a-posteriori hidden channels, $\text{hcp}$ and $\text{hcd}$, respectively.

The scope of the hyperparameter analysis (showed in Appendix \ref{sec:arch}) is not to produce the best possible accuracy for the model, but to explore the sensitivity of the architecture to different autoencoder configurations. For this reason, we have not focused on the usual ones, e.g.\  activation functions, learning rates and number of epochs, nor exploited multiple restart approaches with different seeds.
In particular, we considered the combinations of the following configurations for the aforementioned hyperparameters\footnote{In the pooling case, testing the architecture without convolutional layers is equivalent to learning the map only from sensors data.}: $r_\text{t} \in [10, 30, 50]$, $r_\text{p} \in [30, 50, 70]$, $\text{ffn} \in [100, 200, 300]$, $n_l \in [50, 100]$, $n \in [15, 25]$, $\lambda \in [0.1, 1, 10]$, and $\text{hcp}, \text{hcd} \in [1, 2, 3]$.

We remark that, as in CNNs, our GNN architecture shares the weights across the nodes\footnote{The weights multiply the hidden embedding $\mathbf{h}_u^{(k)}$, thus the matrix $\mathbf{W}^{(k)}$ depends only on the number of features and filters, but not on the number of neighbors $|\mathsf{N}(u)|$. This allows for many possible generalizations, from different meshes to unseen nodes.}, making it scalable to large datasets, in contrast to fully-connected autoencoders. This results in a much lower number of trainable parameters, and substantial speed-up in the training phase. The size of the optimization problem is thus dominated by the MLP dimension. For this reason the pooling approach could be fundamental for finer meshes.

This geometry-informed architecture not only reduces the computational complexity, but also fosters learning. Indeed, already with a relatively small training horizon, i.e.\ $N_\text{epochs} = 5000$, the network reaches accurate reconstruction errors for the benchmarks, and still show a decreasing behavior. 
Finally, we note that GCA-ROM is a convolve-then-reshape method, extracting the inductive bias directly from its original framework, while the ones exploiting standard CNNs are reshape-then-convolve approaches.

\section{Numerical results}\label{sec:num_res}
{In this section, we present different applications showcasing the capabilities of GCA-ROM as an efficient nonlinear reduced order model. The main goal of this approach consists in providing an interpretable and effective framework to investigate parametric PDEs defined over parametrized domains, to enhance the performance of linear ROMs, and to decrease the number of expensive simulations during the offline phase.
Toward these goals, we highlight the main properties of the GCA-ROM architecture when facing these issues: (i) its consistency and versatility when dealing with unstructured varying geometries, (ii) its strength in providing better accuracy/speedup for problems with slow Kolmogorov $n$-width decay, and (iii) its expressive power, characterized by high generalizability even in the low-data regime.}

We analyze the performance of the proposed methodology with respect to four benchmark parametrized PDEs: (i) nonlinear Poisson equation, (ii) advection-dominated problem, (iii) Graetz flow, and (iv) Navier-Stokes system. Each of these problems features a different complexity that GCA-ROM is capable of handling. 
We compare the performance of the novel architecture w.r.t.\ standard intrusive techniques such as POD-Galerkin method, and the original data-driven DL-ROM, both in terms of accuracy and speed-up. {
As detailed below, the method presents remarkable improvements w.r.t.\ state-of-the-art methodologies, and paves the way for new investigations bridging standard numerical techniques for PDEs and neural network approaches.}

In particular, we test the architecture, with and without exploiting up- and down-sampling operations. This way, one can distinguish between the generalization property in the plain setting, and the case in which one has only partial access to the high-fidelity solutions, i.e.\ through sensors.

As concerns the error analysis, we denote with $$\varepsilon_{GCA}(\boldsymbol{\mu}) = \frac{\norm{\textbf{u}_\N(\boldsymbol{\mu}) - \tilde{\textbf{u}}_\N(\boldsymbol{\mu})}_2}{\norm{\textbf{u}_\N(\boldsymbol{\mu})}_2} \quad \text{and} \quad \overline{\varepsilon}_{GCA} = \frac{1}{(1-r_\text{t})N_\text{S}}\sum_{\boldsymbol{\mu} \in \Xi_\text{te}}\varepsilon_{GCA}(\boldsymbol{\mu})$$ the GCA-ROM relative error for a given parameter $\boldsymbol{\mu} \in \mathcal{P}$ and its mean over the testing set $\Xi_\text{te}$, respectively.

Notice that, rather than finding the optimal configuration of the architecture, we focus on testing several hyperparameter's configurations to investigate the robustness of the methodology w.r.t.\ the training procedure and the influence of the graphs' architecture. Thus, we will show the best results within our analysis, and refer to Appendix \ref{sec:appendix} for a more detailed description of the test cases and architectures, and performance w.r.t.\ different configurations.

\subsection{Nonlinear Poisson equation}\label{sec:poi_result}
The first benchmark consists of a nonlinear elliptic problem in a two-dimensional hexagonal domain $\Omega$ with a hole, depicted in Figure \ref{fig:nonlinear_elliptic_domain}. With this toy problem, we aim at testing the ability of GCA-ROM to cope with an unstructured mesh with $N_h = 2562$ nodes, where the Euclidean geometry based on Cartesian coordinates does not provide a valid description of the domain. 
{
    \begin{figure}[!ht]
    \centering
\begin{minipage}{0.4\textwidth}
    \centering
    \def\l{2}
    \def\ll{1}
    \def\a{{\l*sqrt(3.)/2}}
    \def\b{{\l*sqrt(3.)}}
    \begin{tikzpicture}[scale=1.]
    \node[below] at (\ll,0) {\normalsize{(1, 0)}};
    \node[below] at (\l+\ll,0) {\normalsize{(3, 0)}};
    \node[right] at (2*\l,\a) {\normalsize{(4, $\sqrt{3}$)}};
    \node[above] at (\l+\ll,\b){\normalsize{(3, 2$\sqrt{3}$)}};
    \node[above] at (\ll,\b) {\normalsize{(1, 2$\sqrt{3}$)}};
    \node[left] at (0,\a) {\normalsize{(0, $\sqrt{3}$)}};
    \draw[fill=green!55!blue] (\ll,0.) -- (\l+\ll,0.) -- (2*\l,\a) --(\l+\ll, \b) -- (\ll,\b) -- (0,\a) -- cycle;
    \draw[fill=white] (\l,\a) circle[radius=0.5];
    \node[white] at (1.5,0.7) {\Large{$\Omega$}};
    \node[above] at (\l+0.15,\a) {\footnotesize{0.5}};
    \draw[fill=green!55!blue] (\l,\a) -- (\l+0.5,\a);
    \end{tikzpicture}
\end{minipage}\qquad\qquad\qquad
\begin{minipage}{0.4\textwidth}
\centering
    \includegraphics[width=0.73\textwidth]{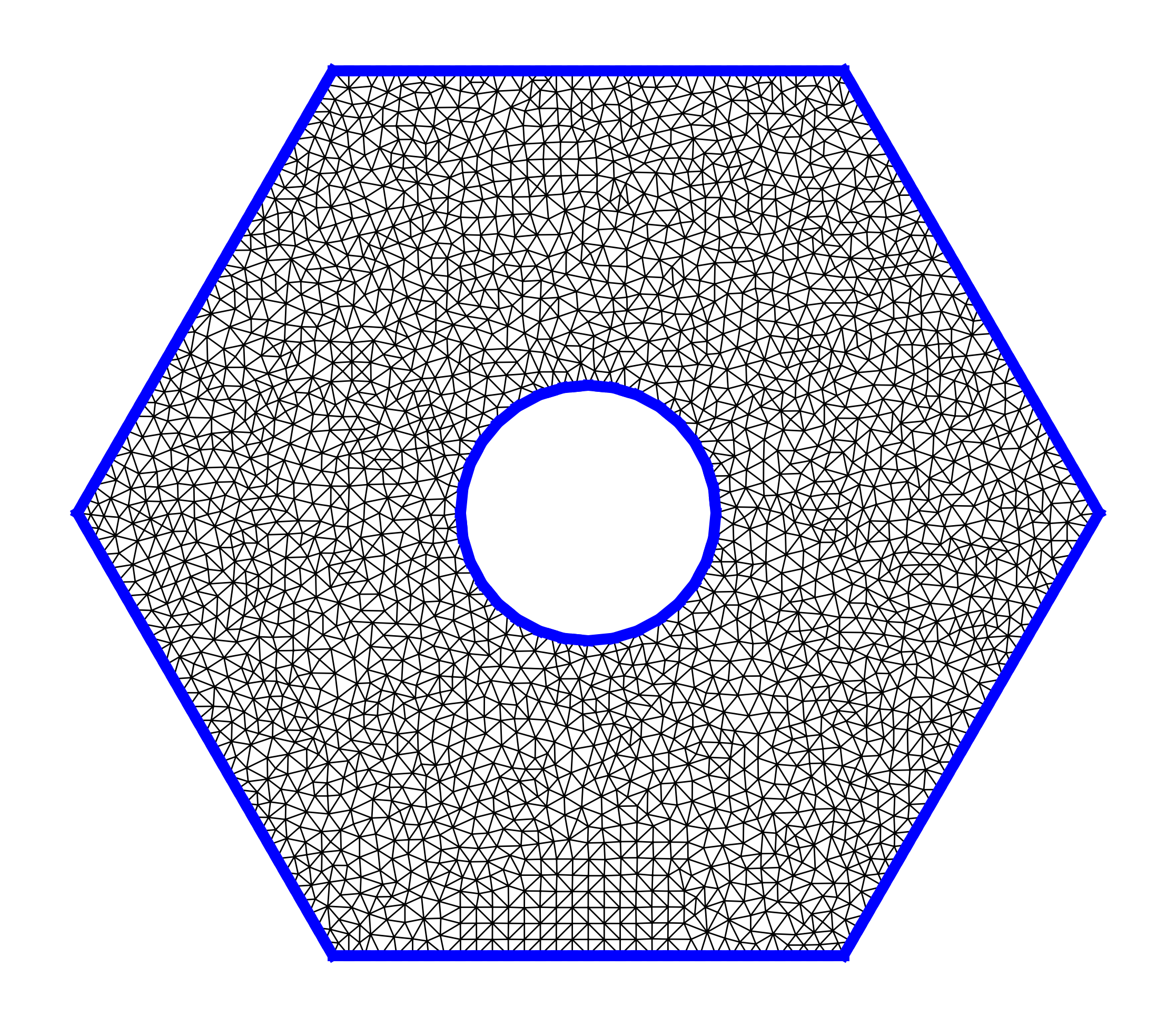}
\end{minipage}
    \caption{{The hexagonal domain with a hole for the nonlinear elliptic problem and its computational mesh.}}
    \label{fig:nonlinear_elliptic_domain}
\end{figure}
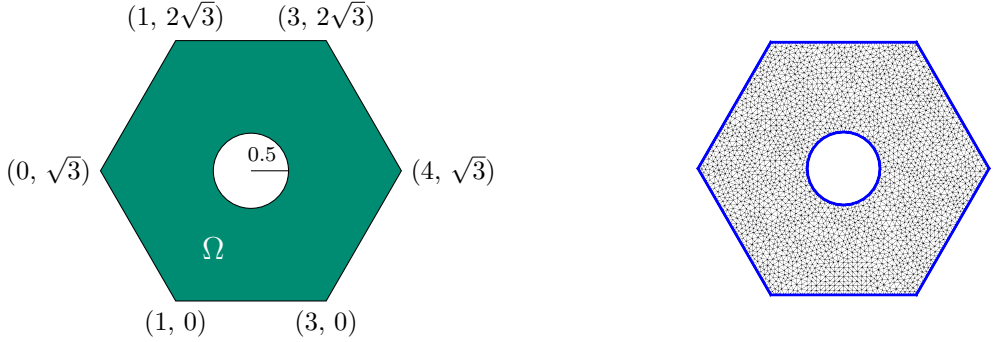}

The continuous formulation of the PDE governing this parametrized problem reads as: given $\boldsymbol{\mu}\in\mathbb{P}$, find $u(\boldsymbol{\mu})$ such that
$$ -\nabla^2u(\boldsymbol{\mu})+\mu_1\frac{e^{\mu_2\,u(\boldsymbol{\mu})}-1}{\mu_2}=g(\boldsymbol{x}),$$
where the physical multi-parameter $\boldsymbol{\mu} = (\mu_1, \mu_2) \in \mathbb{P} = [0.01, 10]^2$ controls the strength of the sink term and of the nonlinearity, and the source term is characterized by the following expression
$$g(\boldsymbol{x}) = 100\sin(2\pi x)\cos(2\pi y), \quad \forall \boldsymbol{x} = (x, y) \in \Omega.$$
The dataset is constituted by $N_\text{S} = 100$ snapshots (10 equispaced parameter samples in each direction) by means of Finite Element approximations over a uniform grid in $\mathbb{P}$. 
The nonlinear term prevents the straightforward applicability of ROMs (hyper-reduction techniques are needed to recover the efficiency), while the lack of a structured mesh leads to a geometrically inconsistent use of standard CNNs.

\begin{figure}[!ht]
	\centering
	\begin{subfigure}{0.32\textwidth}
	\includegraphics[width=1.015\textwidth]{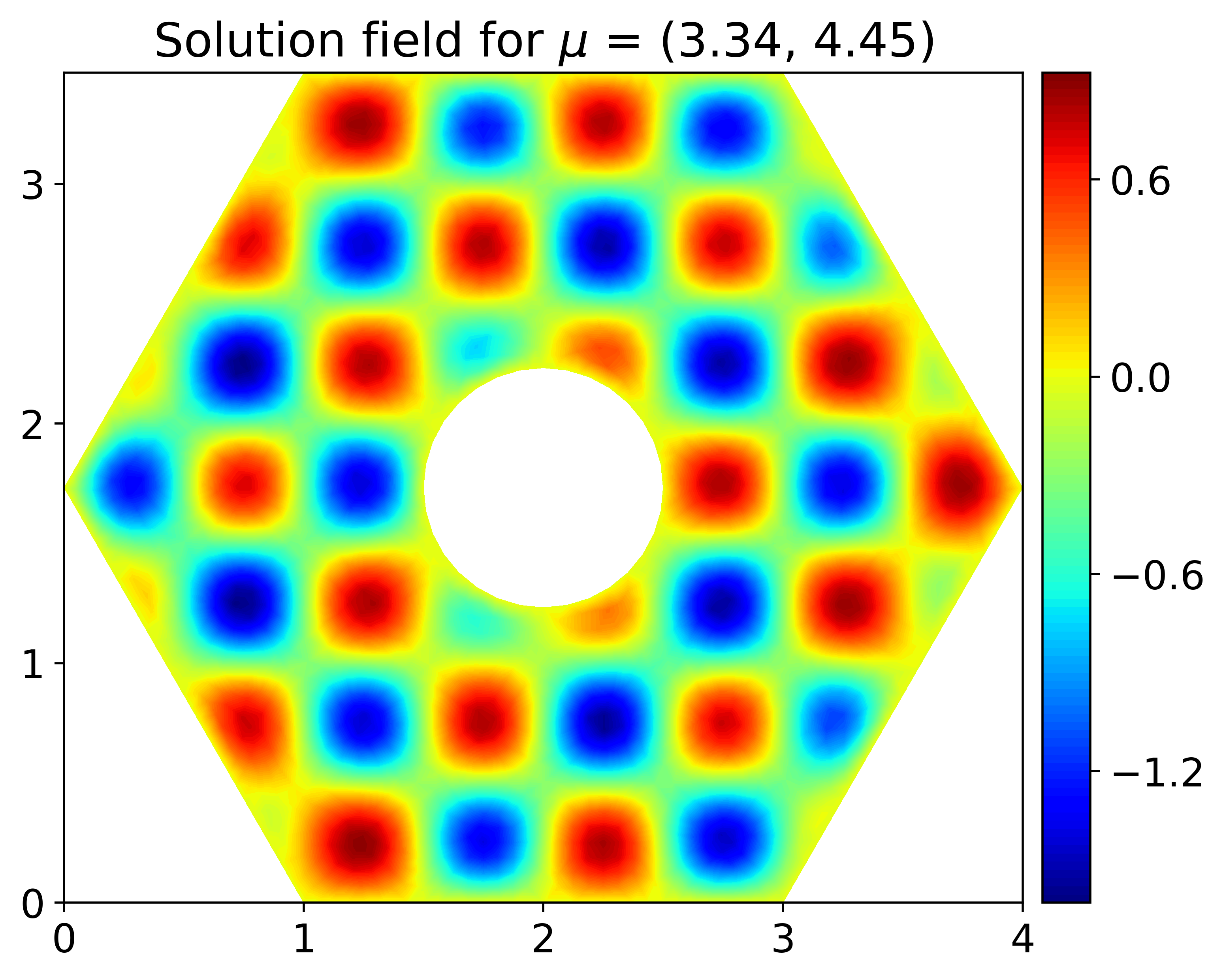}
	\caption{Solution field}
	\label{fig:error_nopooling_poisson_fields_sol}
	\end{subfigure}\hfill
	\begin{subfigure}{0.32\textwidth}
	\includegraphics[width=\textwidth]{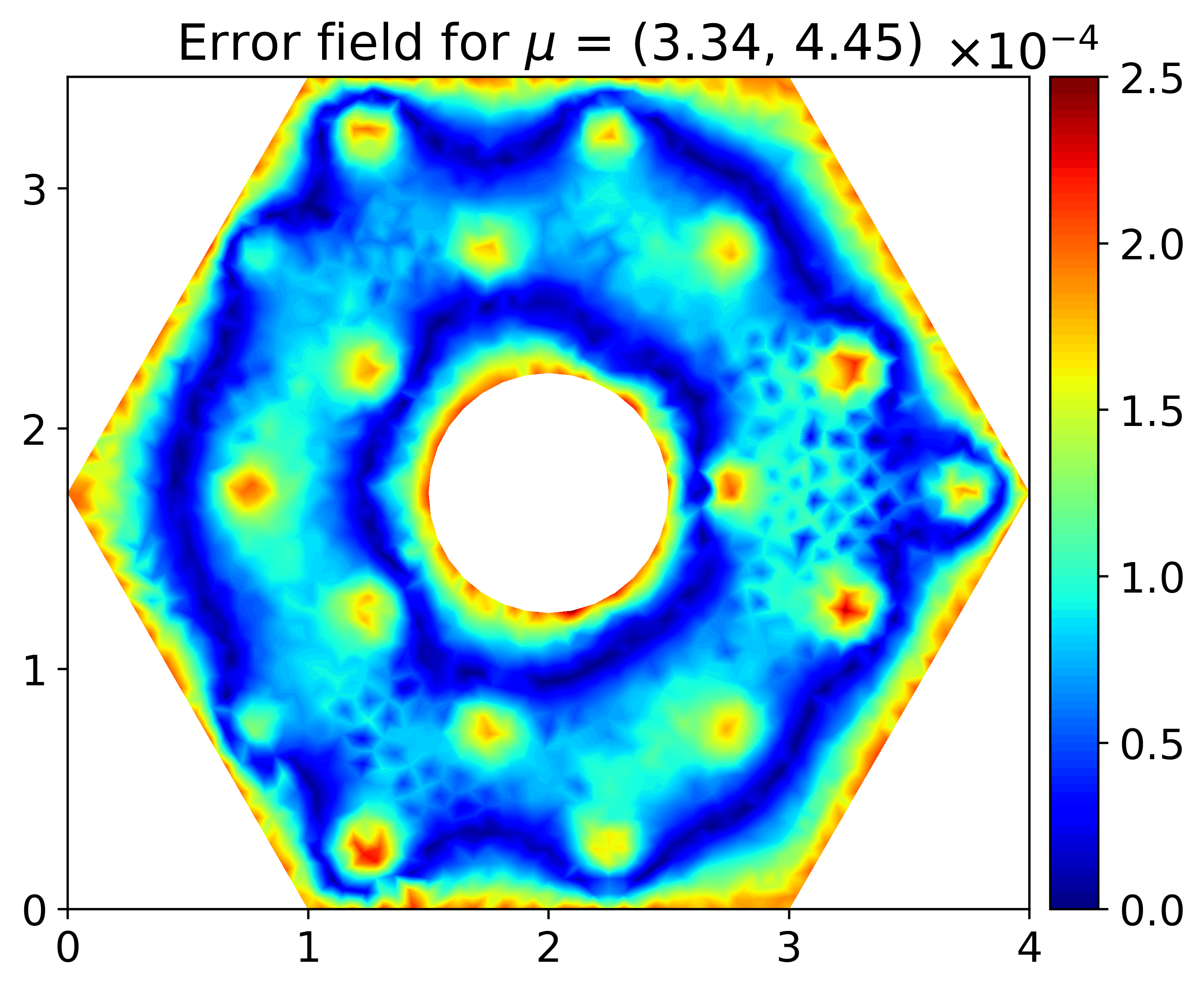}
	\caption{Error field}
	\label{fig:error_nopooling_poisson_fields_e1}
	\end{subfigure}\hfill
	\begin{subfigure}{0.32\textwidth}
	\includegraphics[width=\textwidth]{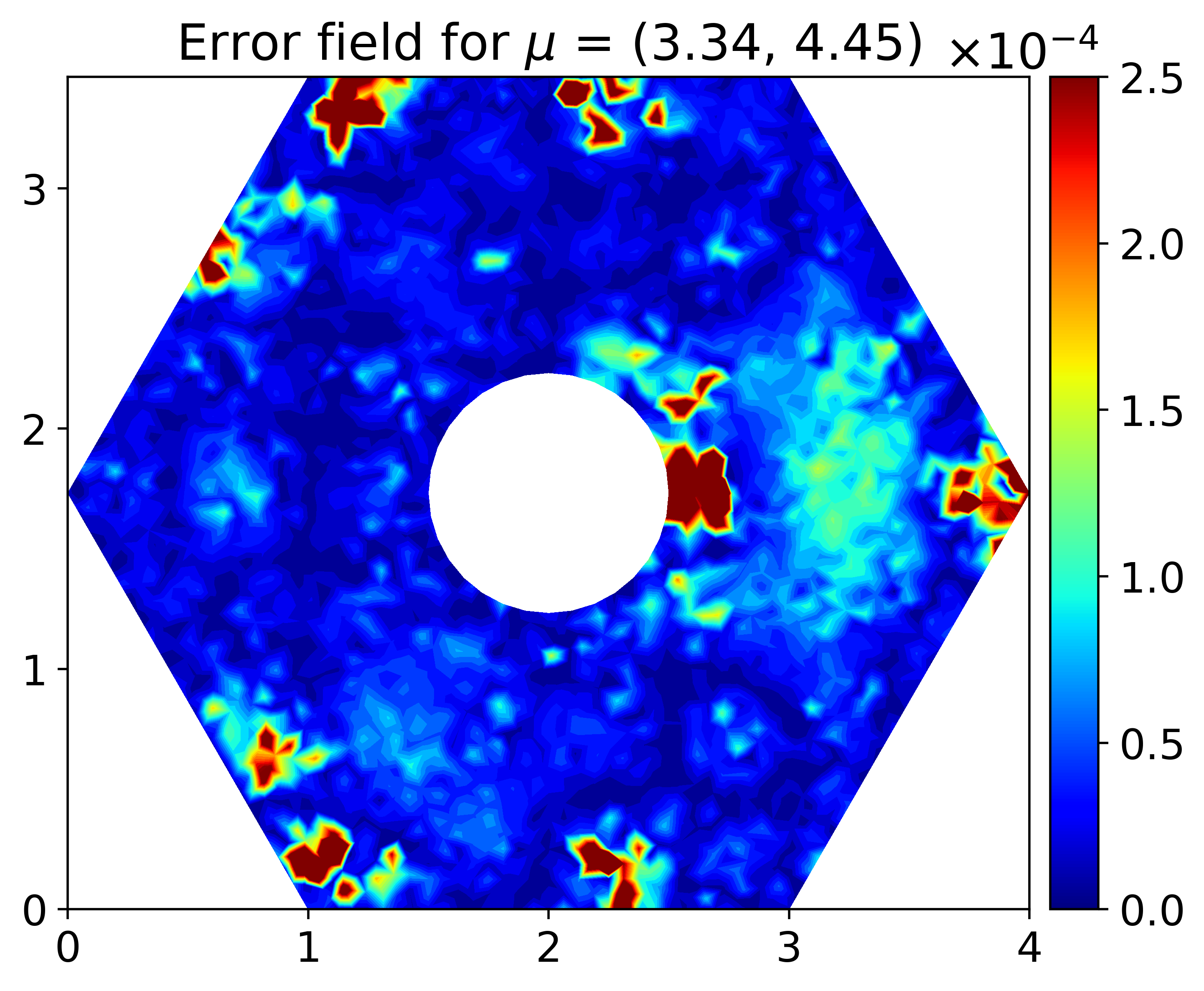}
	\caption{Error field with pooling}
	\label{fig:error_nopooling_poisson_fields_e2}
	\end{subfigure}
	\caption{{Poisson problem for $\boldsymbol{\mu} = (3.34, 4.45)$.}}
	\label{fig:error_nopooling_poisson_fields}
\end{figure}

We show in Figure \ref{fig:error_nopooling_poisson_fields} the high-fidelity solution for Poisson obtained for $\boldsymbol{\mu} = (3.34, 4.45) \in \Xi_\text{te}$, and the GCA-ROM relative errors obtained with and without exploiting pooling procedures. {In both cases, we can accurately reconstruct the corresponding field, with the pooling architecture performing slightly worse due to the high peaks in some regions of the domain, $\varepsilon_{\text{GCA}}(\boldsymbol{\mu}) = 5.0 \times 10^{-3}$ and $\varepsilon_{\text{GCA}}(\boldsymbol{\mu}) = 4.3 \times 10^{-3}$, respectively for pooling and plain versions.} 
This is expected given the availability of only partial information over the solution field during the training stage. 

In Figure \ref{fig:error_nopooling_poisson_fields_e1}, the error is higher near the maxima of the field, {where the dataset $\Xi$ shows the higher variance w.r.t.\ the sampling,} and at the boundary of the domain. This is clearly an undesired effect of the convolution procedure, sharing the message with neighbors even where we impose homogeneous Dirichlet boundary conditions. 
{However, this effect is only present for a few instances of the parameter, while in general the error is spatially distributed in the whole domain, confirming that, regardless the complexity of the solution manifold, GCA-ROM learns the intrinsic behavior and generalize to unseen data.} 
The augmentation of our data-driven architecture with physical information can fix such behavior by guiding the learning task towards physically consistent fields, and will be the subject of forthcoming investigations.
With the pooling approach, see Figure \ref{fig:error_nopooling_poisson_fields_e2}, the accuracy drops especially near the down-sampled regions, where we apply the mask.

\begin{figure}[!ht]
	\centering
	\includegraphics[width=0.4\textwidth]{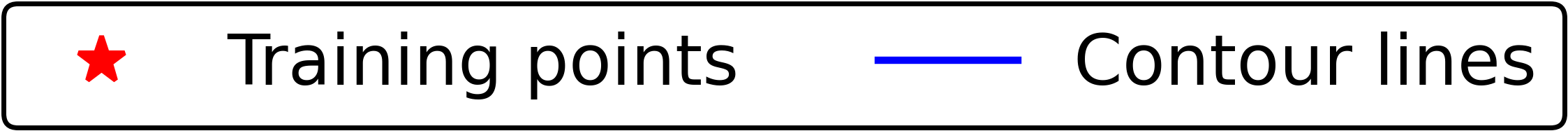}

	\begin{subfigure}{0.49\textwidth}
	\includegraphics[width=\textwidth, clip=true, trim = 20mm 0mm 10mm 0mm]{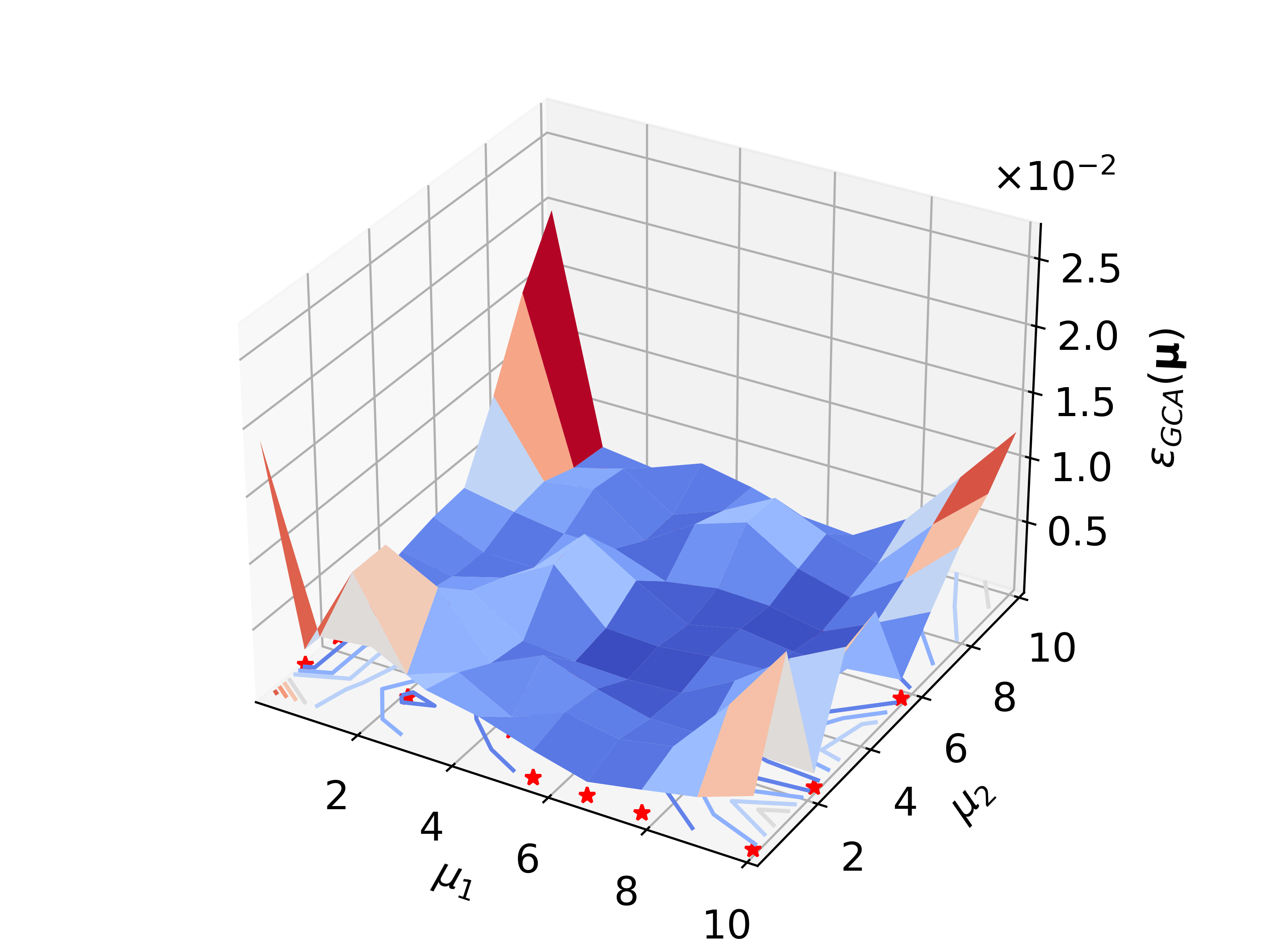}
	\caption{GCA-ROM}
	\label{fig:error_pooling_poisson_fields_nopool}
	\end{subfigure}\hfill
	\begin{subfigure}{0.49\textwidth}
	\includegraphics[width=\textwidth, clip=true, trim = 20mm 0mm 10mm 0mm]{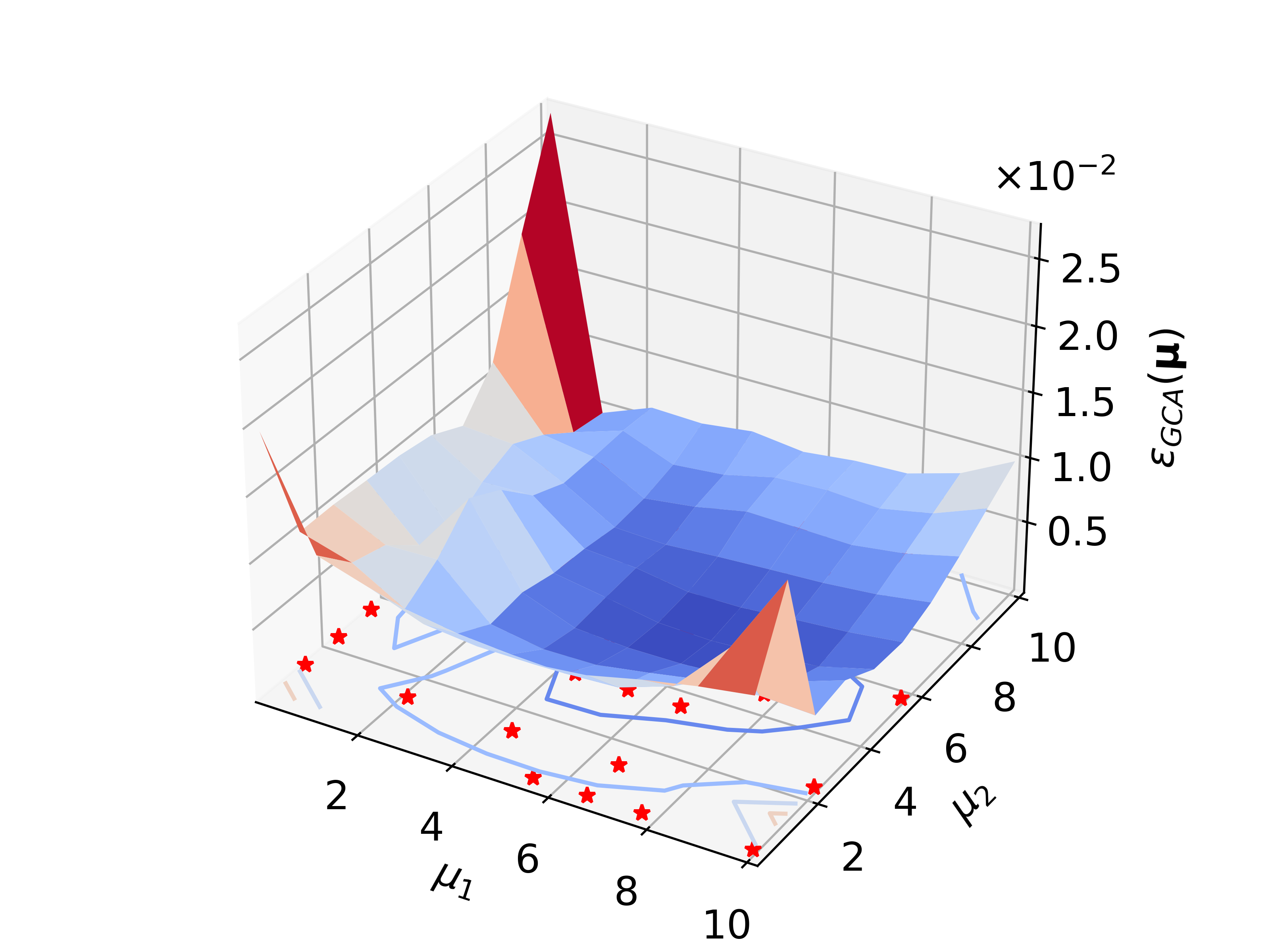}
	\caption{GCA-ROM with pooling}
	\label{fig:error_pooling_poisson_fields_pool}
	\end{subfigure}
	\caption{{GCA-ROM relative errors for the Poisson problem on the dataset $\Xi$, with red markers corresponding to the parameters used in the training set $\Xi_\text{tr}$.}}
	\label{fig:error_pooling_poisson_fields}
\end{figure}

To highlight the generalization capabilities achieved by GCA-ROM, we show in Figure \ref{fig:error_pooling_poisson_fields}, for both approaches, the evolution of the relative error $\varepsilon_{GCA}(\boldsymbol{\mu})$ in the whole dataset $\Xi$. {As expected, the accuracy of the pooling approach, both in terms of mean and maximum errors, marginally deteriorates due to the lack of information, but still produces remarkably close results.}
For the plain architecture in Figure \ref{fig:error_pooling_poisson_fields_nopool}, it can be seen that, exploiting only $r_\text{t} = 30\%$ of the original dataset (reported in the plot with red stars), the relative error over $\Xi_\text{te}$ is below $2 \times 10^{-2}$ and its mean is $\overline{\varepsilon}_{GCA} =  8 \times 10^{-3}$.
The accuracy of the pooling procedure in Figure \ref{fig:error_pooling_poisson_fields_pool} is also surprisingly good, again exploiting only $r_\text{t} = 30\%$ of the dataset, but also retaining only $r_\text{p} = 70\%$ of each high-fidelity solution, i.e.\ we do not have access to these fixed nodes. We are still able to reconstruct the Poisson solutions with maximum relative error below { $2.7 \times 10^{-2}$} and mean $\overline{\varepsilon}_{GCA} = 1 \times 10^{-2}$.

This result shows the potential of graph convolutional architecture in capturing parametric behavior, even when dealing with sparse measurements, and paves the way for the development of robust up- and down-sampling procedures.

\subsection{Advection dominated problem}\label{sec:adv_result}
Standard projection-based reduced order models struggle with slow Kolmogorov decay phenomena, such as those in the advection-dominated regime. Here, we test the performance of our nonlinear and data-driven approach for an advection-diffusion equation with a parametrized transport term. The domain is given by the unstructured mesh defined over the square domain $\Omega = [0, 1]^2 \in \mathbb{R}^2$, depicted in Figure \ref{advection_dominated_domain}, with $N_h = 3967$ nodes.
{
    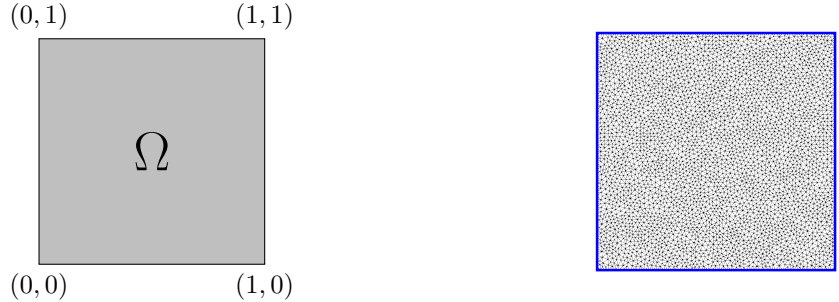
\begin{figure}[!ht]
    \centering
\begin{minipage}{0.4\textwidth}
	\centering
	\begin{tikzpicture}[scale=1.5]
	\draw[fill=lightgray] (-1,-1) rectangle (1,1);
	\node[below] at (-1,-1) {\normalsize{$(0, 0)$}};
	\node[below] at (1,-1) {\normalsize{$(1, 0)$}};
	\node[above] at (-1,1) {\normalsize{$(0, 1)$}};
	\node[above] at (1,1) {\normalsize{$(1, 1)$}};
	\node at (0,0) {\Huge{$\Omega$}};
	\end{tikzpicture}
\end{minipage}\qquad\qquad
\begin{minipage}{0.4\textwidth}
\centering
    \includegraphics[width=0.6\textwidth]{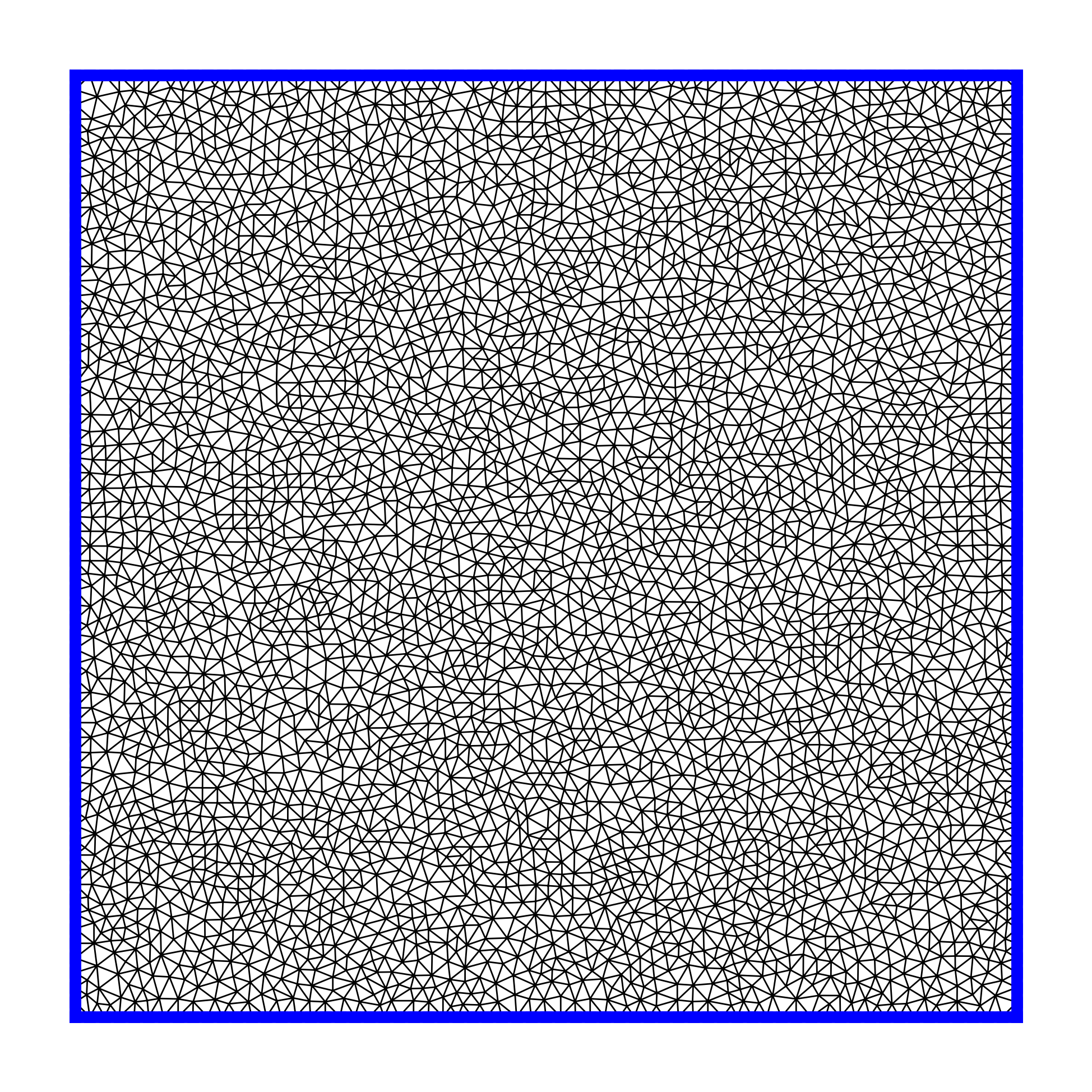}
\end{minipage}
	\caption{{The square domain for the advection-dominated problem and its computational mesh.}}
	\label{advection_dominated_domain}
\end{figure}}

The advection-dominated problem reads as: given $\boldsymbol{\mu}\in\mathbb{P}$, find $u(\boldsymbol{\mu})$ such that
$$
-\alpha(\mu_1)\Delta\,u(\boldsymbol{\mu})+\beta(\mu_2)\cdot\nabla u(\boldsymbol{\mu})=g \quad \text{in }\Omega,
$$
where $\alpha(\mu_1) = 10^{-\mu_1}$ represents the diffusion coefficient defining the Péclet number, $\beta(\mu_2)=(\mu_2,\mu_2)$ is the direction of the transport, and $g=1$ is the source term.
Thus, the model is characterized by the physical multi-parameter $\boldsymbol{\mu} = (\mu_1, \mu_2) \in \mathbb{P} = [0, 6] \times [-1, 1]$. As before, we form the dataset by means of  $N_\text{S} = 100$ snapshots over a uniform grid in $\mathbb{P}$.

\begin{figure}[!ht]
	\centering
	\begin{subfigure}{0.32\textwidth}
	\includegraphics[width=\textwidth]{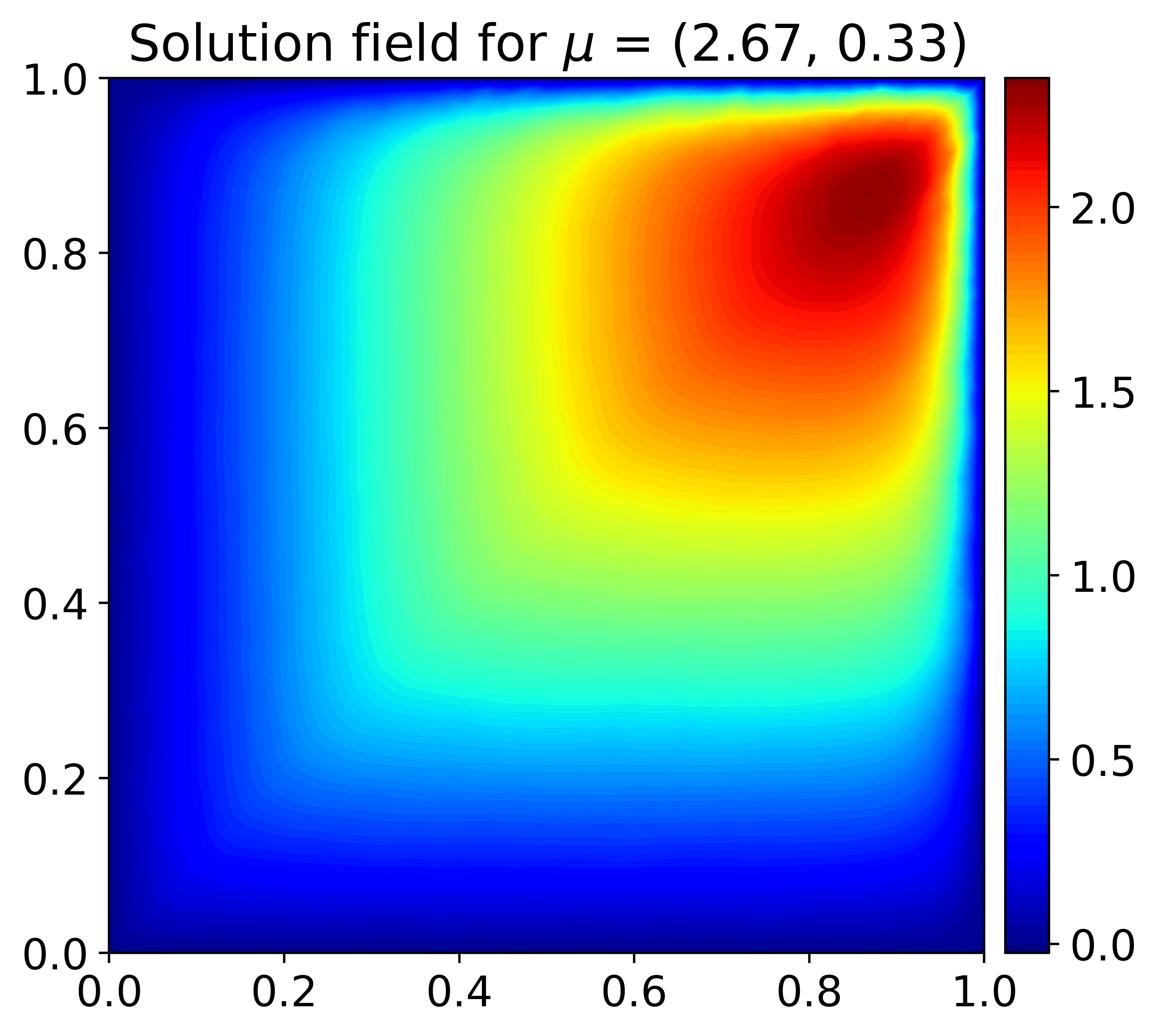}
	\caption{Solution field}
	\label{fig:error_nopooling_advection_fields_sol}
	\end{subfigure}\hfill
	\begin{subfigure}{0.32\textwidth}
	\includegraphics[width=\textwidth]{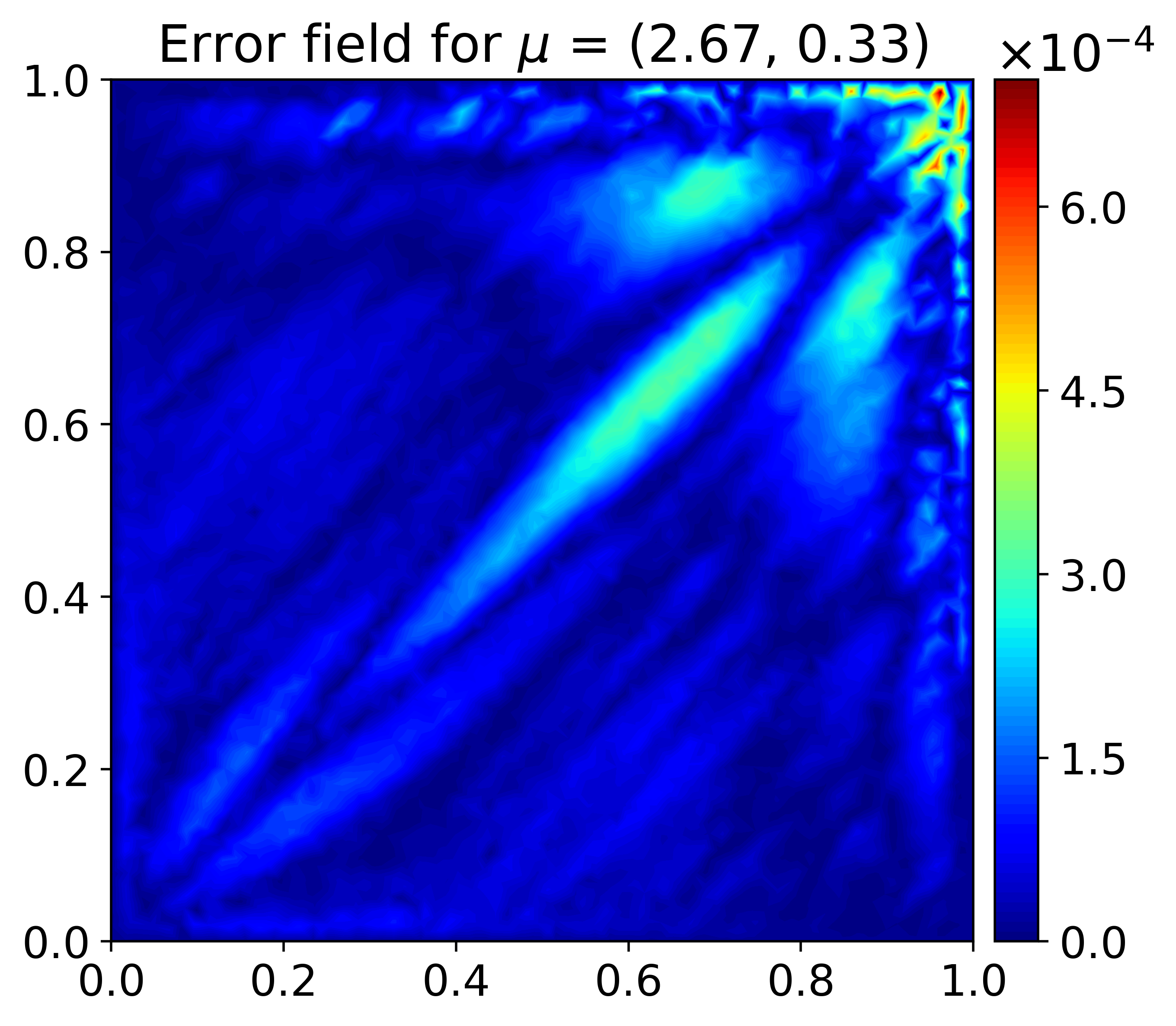}
	\caption{Error field}
	\label{fig:error_nopooling_advection_fields_e1}
	\end{subfigure}\hfill
	\begin{subfigure}{0.32\textwidth}
	\includegraphics[width=\textwidth]{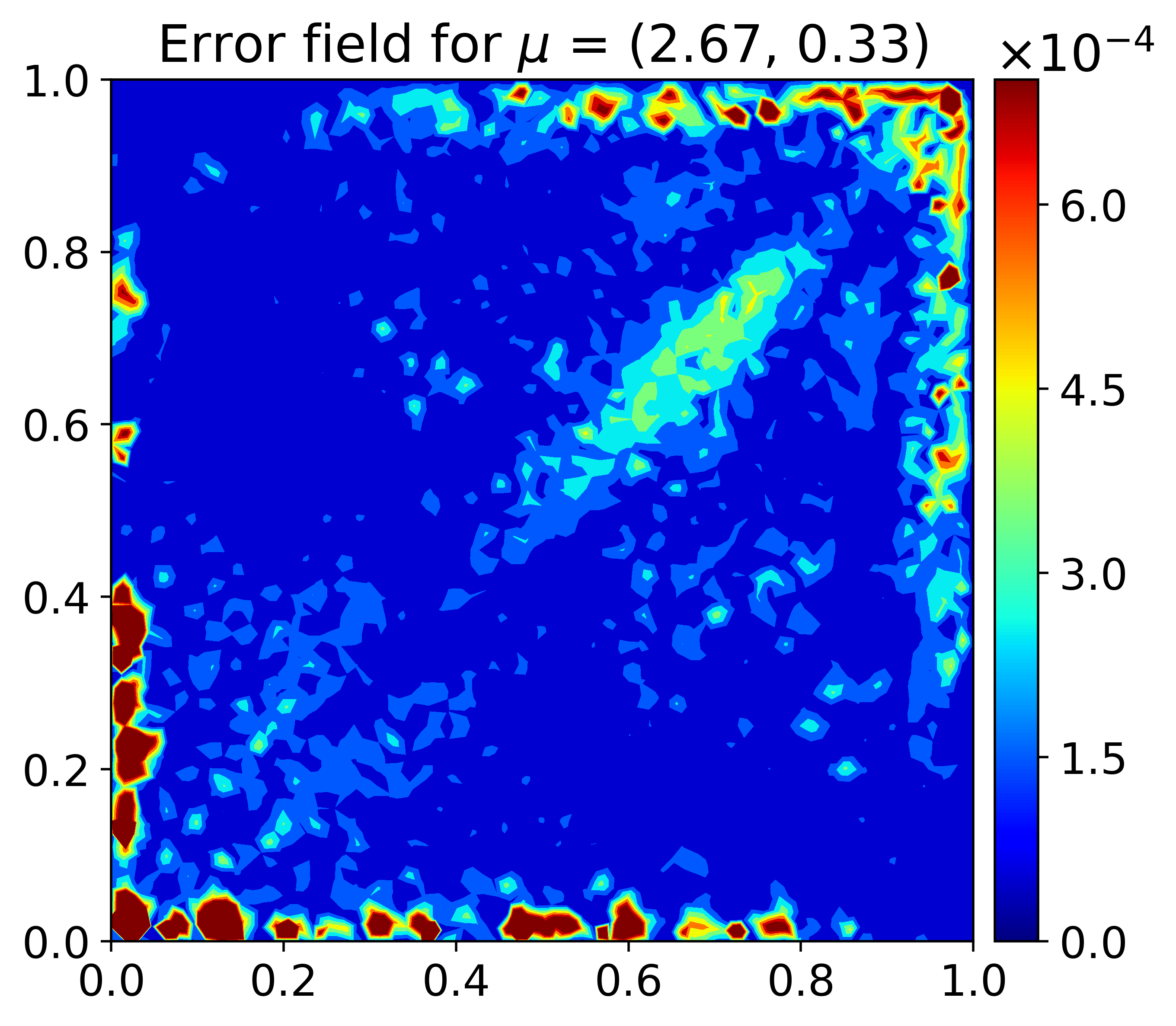}
	\caption{Error field with pooling}
	\label{fig:error_nopooling_advection_fields_e2}
	\end{subfigure}
	\caption{{Advection problem for $\boldsymbol{\mu} = (2.67, 0.33)$.}}
	\label{fig:error_nopooling_advection_fields}
\end{figure}

In Figure \ref{fig:error_nopooling_advection_fields} we plot the solution of the advection dominated problem for $\boldsymbol{\mu} = (2.67, 0.33) \in \Xi_\text{te}$, and the two GCA-ROM relative error fields. Again, the qualitative behavior of the methodology shows accurate reconstructions, with the maximum errors located near the boundary layer's region.
Given the different amount of information, the error in Figure \ref{fig:error_nopooling_advection_fields_e1} is more distributed over the domain, and one order of magnitude smaller than the one in Figure \ref{fig:error_nopooling_advection_fields_e2}.

\begin{figure}[!ht]
	\centering
	\includegraphics[width=0.4\textwidth]{legend.png}

	\begin{subfigure}{0.49\textwidth}
	\includegraphics[width=\textwidth, clip=true, trim = 20mm 0mm 10mm 0mm]{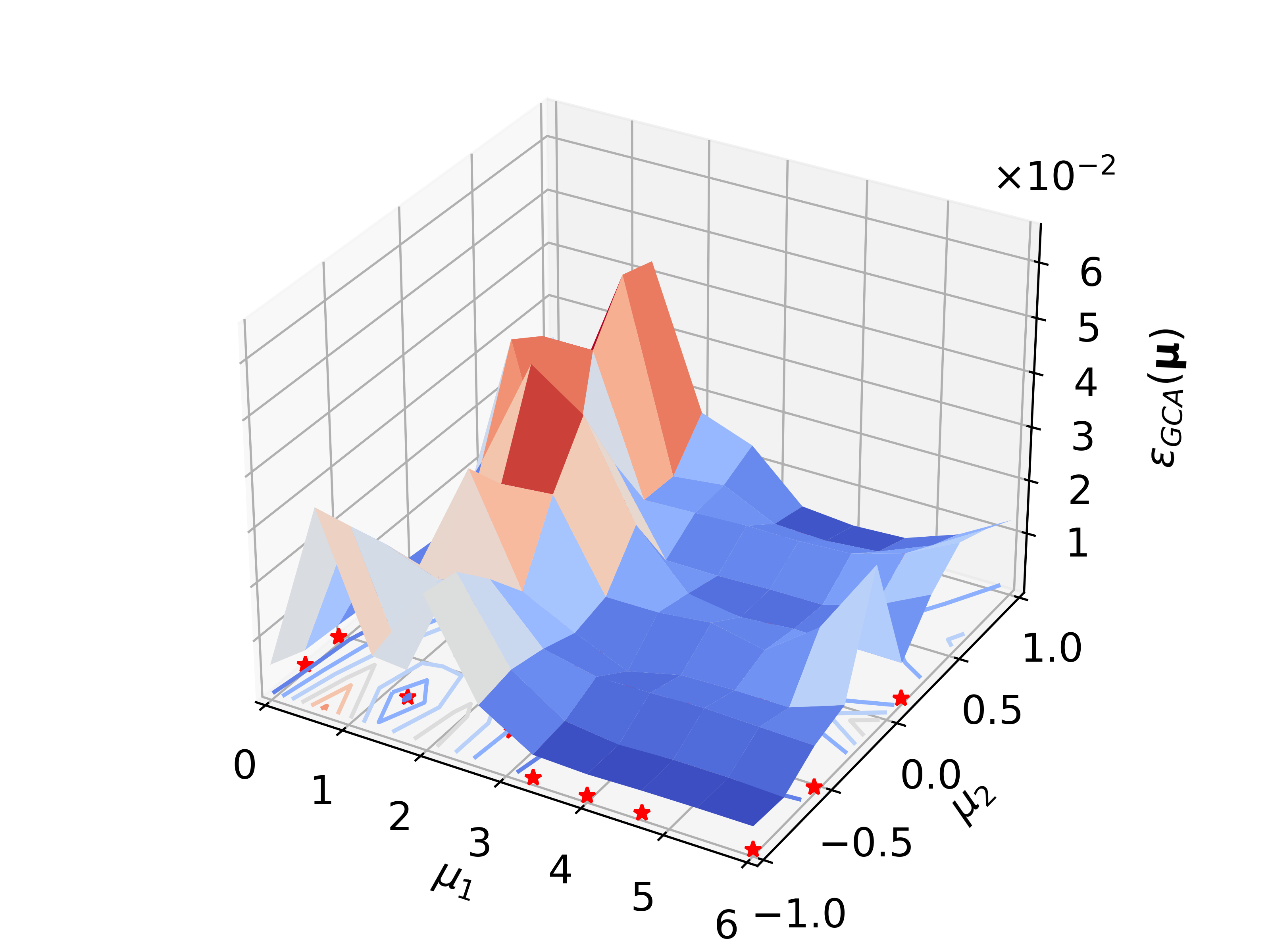}
	\caption{GCA-ROM}
	\label{fig:error_pooling_advection_fields_nopool}
	\end{subfigure}\hfill
	\begin{subfigure}{0.49\textwidth}	
	\includegraphics[width=\textwidth, clip=true, trim = 20mm 0mm 10mm 0mm]{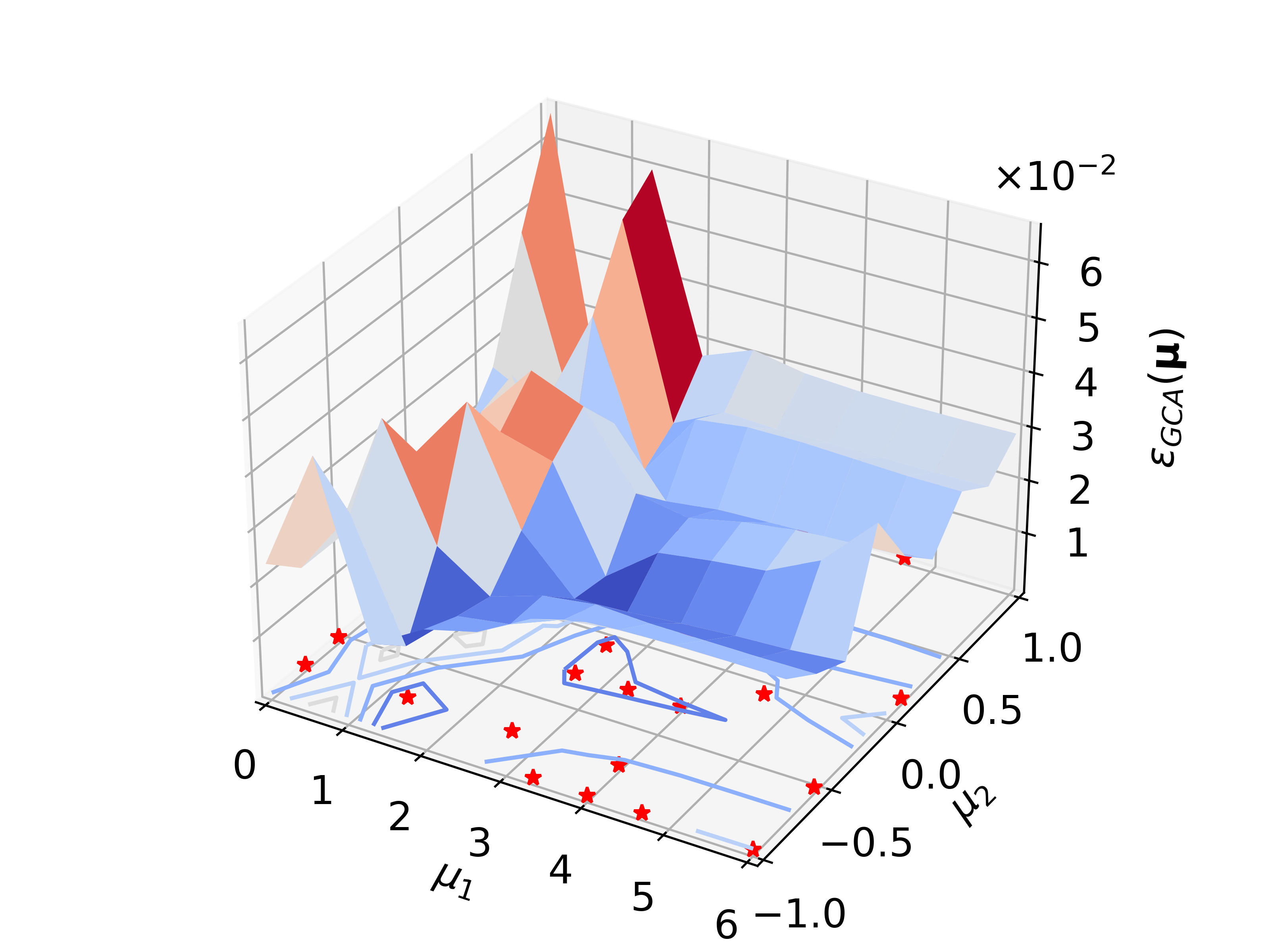}
	\caption{GCA-ROM with pooling}
	\label{fig:error_pooling_advection_fields_pool}
	\end{subfigure}
	\caption{{GCA-ROM relative errors for the Advection problem on the dataset $\Xi$, with red markers corresponding to the parameters used in the training set $\Xi_\text{tr}$.}}
	\label{fig:error_pooling_advection_fields}
\end{figure}

For the accuracy, we plot in Figure \ref{fig:error_pooling_advection_fields} the relative errors $\varepsilon_{GCA}(\boldsymbol{\mu})$ in the plain and down-sampled architectures over the dataset $\Xi$.
In both cases, the training set comprises only $r_\text{t} = 30\%$ of the computed snapshots, but the maximum relative error over $\Xi_\text{te}$ is around $5.2 \times 10^{-2}$, while its mean is $\overline{\varepsilon}_{GCA} =  2.4 \times 10^{-2}$. We stress the fact that to learn this complex behavior, we are only exploiting 30 high-fidelity solutions fixed a-priori by the seed, and testing on the remaining 70.
The pooling strategy gives comparable results by selecting $r_\text{p} = 70\%$ nodes, and thus masking $30\%$ of them, with maximum relative error below $7.2 \times 10^{-2}$ and mean $\overline{\varepsilon}_{GCA} = 3.3 \times 10^{-2}$.

The graph structure is thus still able to predict unseen patterns, even when dealing with a combination of: a low-data regime, scattered mesh information, and complex parametric behavior compromising the linear reducibility of the model. 

\subsection{Graetz flow model}\label{sec:gra_result}
One of the main advantages of the GCA-ROM architecture is the possibility to work directly on the (unstructured) computational domain, introducing geometric priors and defining physically consistent convolutional operations. Now, we consider the case in which the PDE is posed on a geometrically parametrized domain\footnote{The geometric parametrization is performed w.r.t.\ affine mappings, thus producing a set of snapshots defined on the same number of nodes as the original mesh cardinality.}, i.e.\ when the distance between nodes changes w.r.t.\ the parameter. Standard convolutions are agnostic to this behavior.  We study a steady-state Graetz problem on the parametrized rectangular geometry $\Omega(\mu_1) = \Omega_1 \cup \Omega_2(\mu_1) \in \mathbb{R}^2$ depicted in Figure \ref{fig:04_graetz_domain}, with $\Omega_1 = [0, 1] \times [0, 1]$ and $\Omega_2 = [0, \mu_1] \times [0, 1]$, and discretized with $N_h = 5160$ nodes.
{
    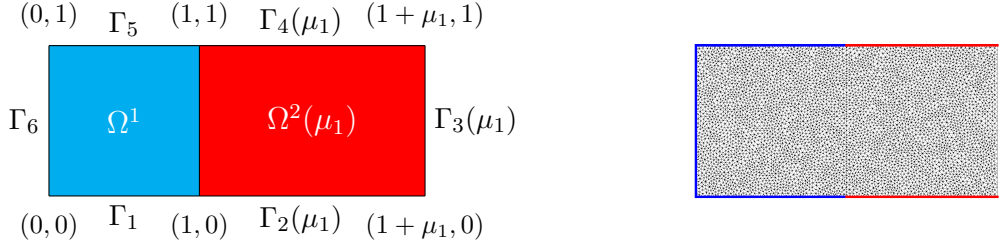
\begin{figure}[!ht]
    \centering
\begin{minipage}{0.4\textwidth}
	\centering
\begin{tikzpicture}[scale=1]
\fill [cyan] (-1,-1) rectangle (1,1);
\fill [red] (1,-1) rectangle (4,1);
\node[below,yshift=-0.1cm] at (-1,-1) {\normalsize{$(0,0)$}};
\node[below,yshift=-0.1cm] at (1,-1) {\normalsize{$(1,0)$}};
\node[below,yshift=-0.1cm] at (4,-1) {\normalsize{$(1+\mu_1,0)$}};
\node[above,yshift=0.1cm] at (-1,1) {\normalsize{$(0,1)$}};
\node[above,yshift=0.1cm] at (1,1) {\normalsize{$(1,1)$}};
\node[above,yshift=0.1cm] at (4,1) {\normalsize{$(1+\mu_1,1)$}};
\draw (-1,-1) -- (1,-1) node[pos=0.5,sloped,below] {\large{$\Gamma_{1}$}};
\draw (1,-1) -- (1,1);
\draw (1,1) -- (-1,1) node[pos=0.5,sloped,above] {\large{$\Gamma_{5}$}};
\draw (-1,1) -- (-1,-1) node[pos=0.5,left] {\large{$\Gamma_{6}$}};
\draw (1,-1) -- (4,-1) node[pos=0.45,sloped,below] {\large{$\Gamma_{2}(\mu_1)$}};
\draw (4,-1) -- (4,1) node[pos=0.5,right] {\large{$\Gamma_{3}(\mu_1)$}};
\draw (4,1) -- (1,1) node[pos=0.55,sloped,above] {\large{$\Gamma_{4}(\mu_1)$}};
\node[white] at (0,0) {\Large{$\Omega^1$}};
\node[white] at (2.5,0) {\Large{$\Omega^2(\mu_1)$}};
\end{tikzpicture}
\end{minipage}\qquad\qquad\qquad
\begin{minipage}{0.4\textwidth}
\centering
    \includegraphics[width=0.75\textwidth]{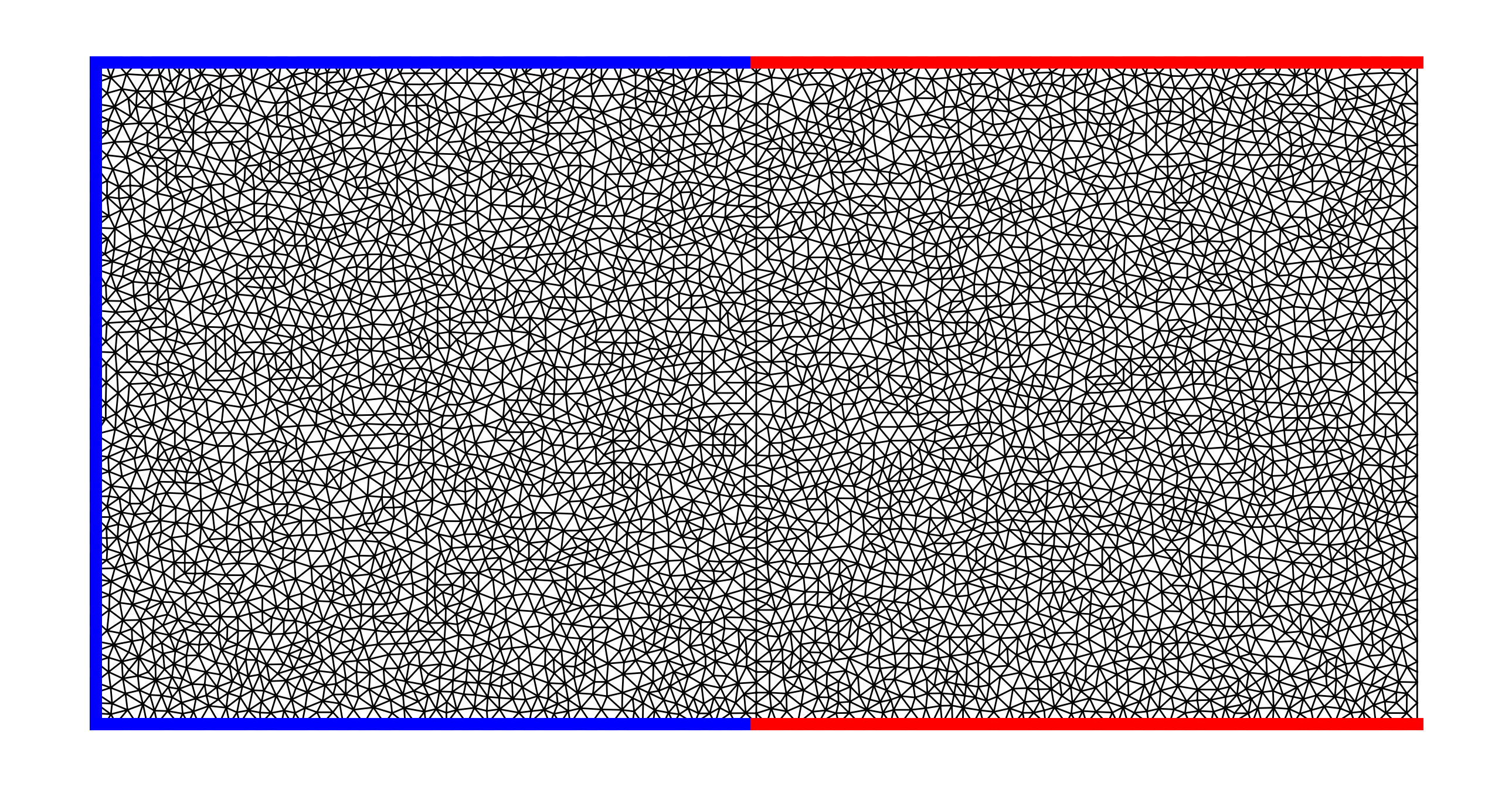}
\end{minipage}
\caption{{The geometrically parametrized domain for the Graetz problem and its computational mesh.}}
\label{fig:04_graetz_domain}
\end{figure}}

The Graetz model describes a thermal field under the combined effects of heat's diffusion and advection along the horizontal direction, with homogeneous and non-homogeneous boundary conditions.
The continuous PDE's formulation of the Graetz problem reads as: given $\boldsymbol{\mu}\in\mathbb{P}$, find $u(\boldsymbol{\mu})$ such that
$$
\begin{cases}
-\mu_2\Delta u(\boldsymbol{\mu})+y(1-y)\partial_{x}u(\boldsymbol{\mu}) = 0 & \text{in } \Omega(\mu_1),\\
u(\boldsymbol{\mu}) = 0 & \text{on } \Gamma_{\text{D}},\\
u(\boldsymbol{\mu}) = 1 & \text{on } \Gamma_{\text{G}}(\mu_1),\\
\frac{\partial u}{\partial\boldsymbol{n}}(\boldsymbol{\mu})=0 & \text{on } \Gamma_{\text{N}}(\mu_1),
\end{cases}
$$
where $\boldsymbol{n}$ denotes the outer normal to the boundary $\Gamma_{\text{N}}$.
We impose homogeneous Dirichlet conditions on $\Gamma_{\text{D}} = \Gamma_{1} \cup \Gamma_{5} \cup \Gamma_{6}$, while the temperature is kept at $u=1$ at the boundaries $\Gamma_{\text{G}} = \Gamma_{2}(\mu_1) \cup \Gamma_{4}(\mu_1)$. Finally, homogeneous Neumann conditions are imposed on $\Gamma_{\text{N}} = \Gamma_{3}(\mu_1)$. 

The model has been parametrized by means of a physical and geometrical multi-parameter $\boldsymbol{\mu} = (\mu_1, \mu_2) \in \mathbb{P} = [1, 3] \times [0.01, 0.1]$, controlling the length $\mu_1$  of the idealized pipes, and the diffusivity coefficient $\mu_2$. The dataset is formed by $N_\text{S} = 200$ snapshots, computed on a uniform grid in $\mathbb{P}$ with respectively 10 and 20 equispaced values for the two parameters.

\begin{figure}[t]
	\centering
	\begin{subfigure}{0.6\textwidth}
	\includegraphics[width=\textwidth]{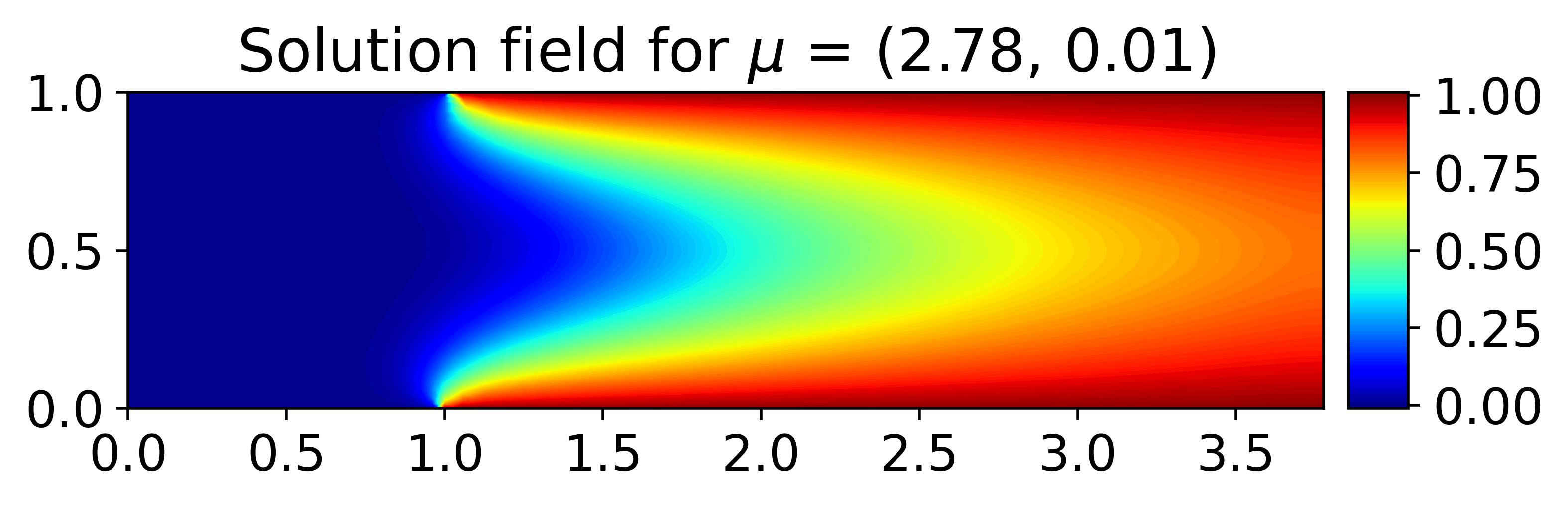}
	\caption{Solution field}
	\label{fig:error_nopooling_graetz_fields_sol}
	\end{subfigure}\\\vspace{0.2cm}

	\begin{subfigure}{0.6\textwidth}
	\includegraphics[width=\textwidth]{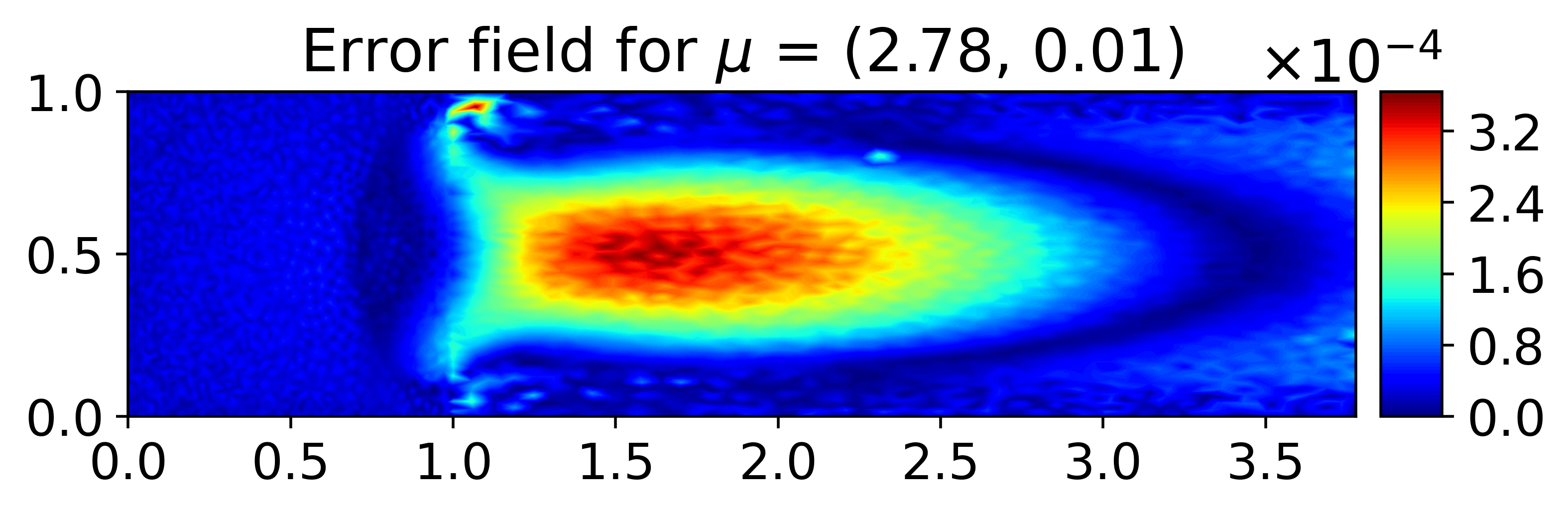}
	\caption{Error field}
	\label{fig:error_nopooling_graetz_fields_e1}
	\end{subfigure}\\\vspace{0.2cm}

	\begin{subfigure}{0.6\textwidth}
	\includegraphics[width=\textwidth]{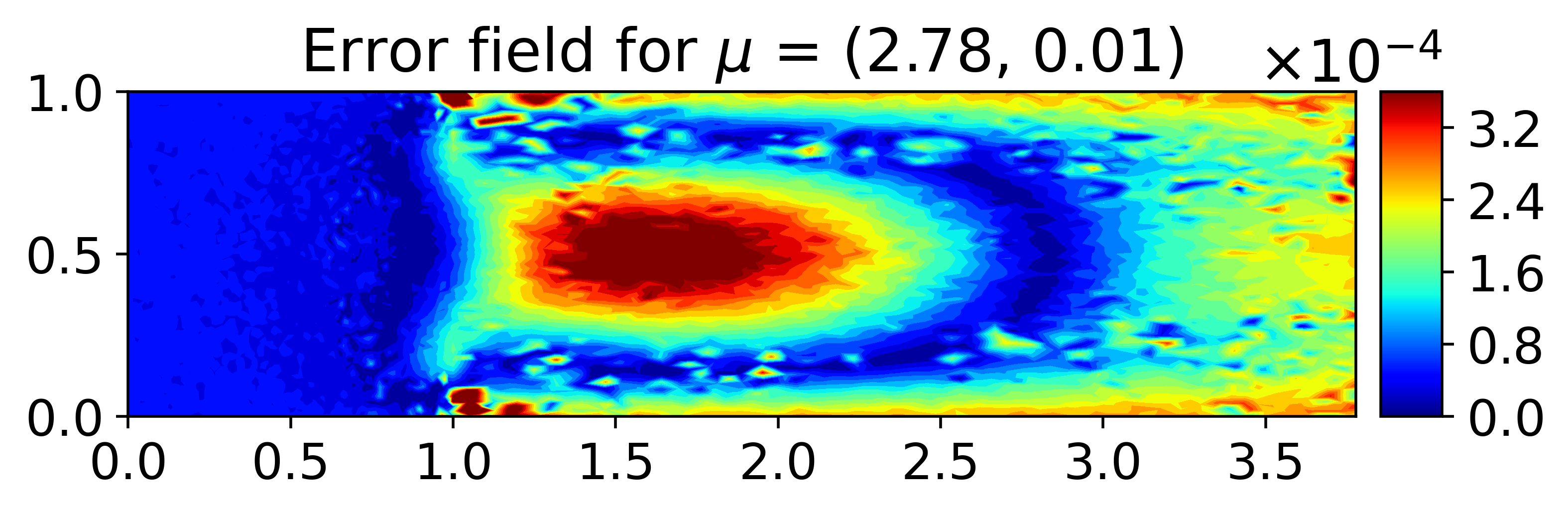}
	\caption{Error field with pooling}
	\label{fig:error_nopooling_graetz_fields_e2}
	\end{subfigure}
	\caption{{Graetz flow model for $\boldsymbol{\mu} = (2.78, 0.01)$.}}
	\label{fig:error_nopooling_graetz_fields}
\end{figure}

In Figure \ref{fig:error_nopooling_graetz_fields} we plot the solution of the Graetz model obtained for $\boldsymbol{\mu} = (2.78, 0.01) \in \Xi_\text{te}$, comparing the relative errors for two GCA-ROM approaches. The error in Figure \ref{fig:error_nopooling_graetz_fields_e1} has a maximum in the center of the stretched subdomain, where the transport phenomena takes place, while the pooling approach in Figure \ref{fig:error_nopooling_graetz_fields_e2} produces higher errors at the interface between the two boundary conditions.

\begin{figure}[!ht]
	\centering
	\includegraphics[width=0.4\textwidth]{legend.png}

	\begin{subfigure}{0.49\textwidth}
	\includegraphics[width=\textwidth, clip=true, trim = 20mm 0mm 10mm 0mm]{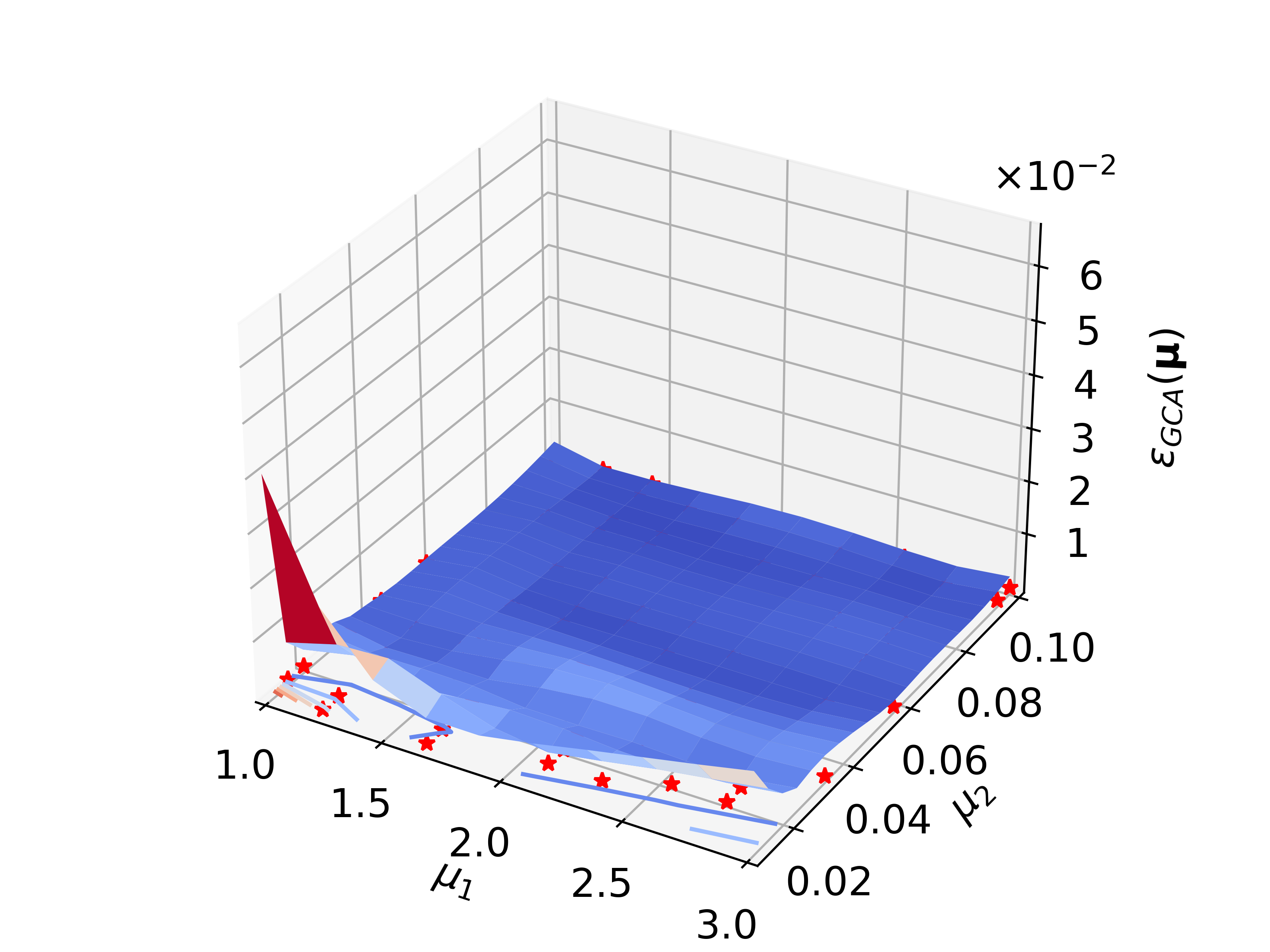}
	\caption{GCA-ROM}
	\label{fig:error_pooling_graetz_fields_nopool}
	\end{subfigure}\hfill
	\begin{subfigure}{0.49\textwidth}	
	\includegraphics[width=\textwidth, clip=true, trim = 20mm 0mm 10mm 0mm]{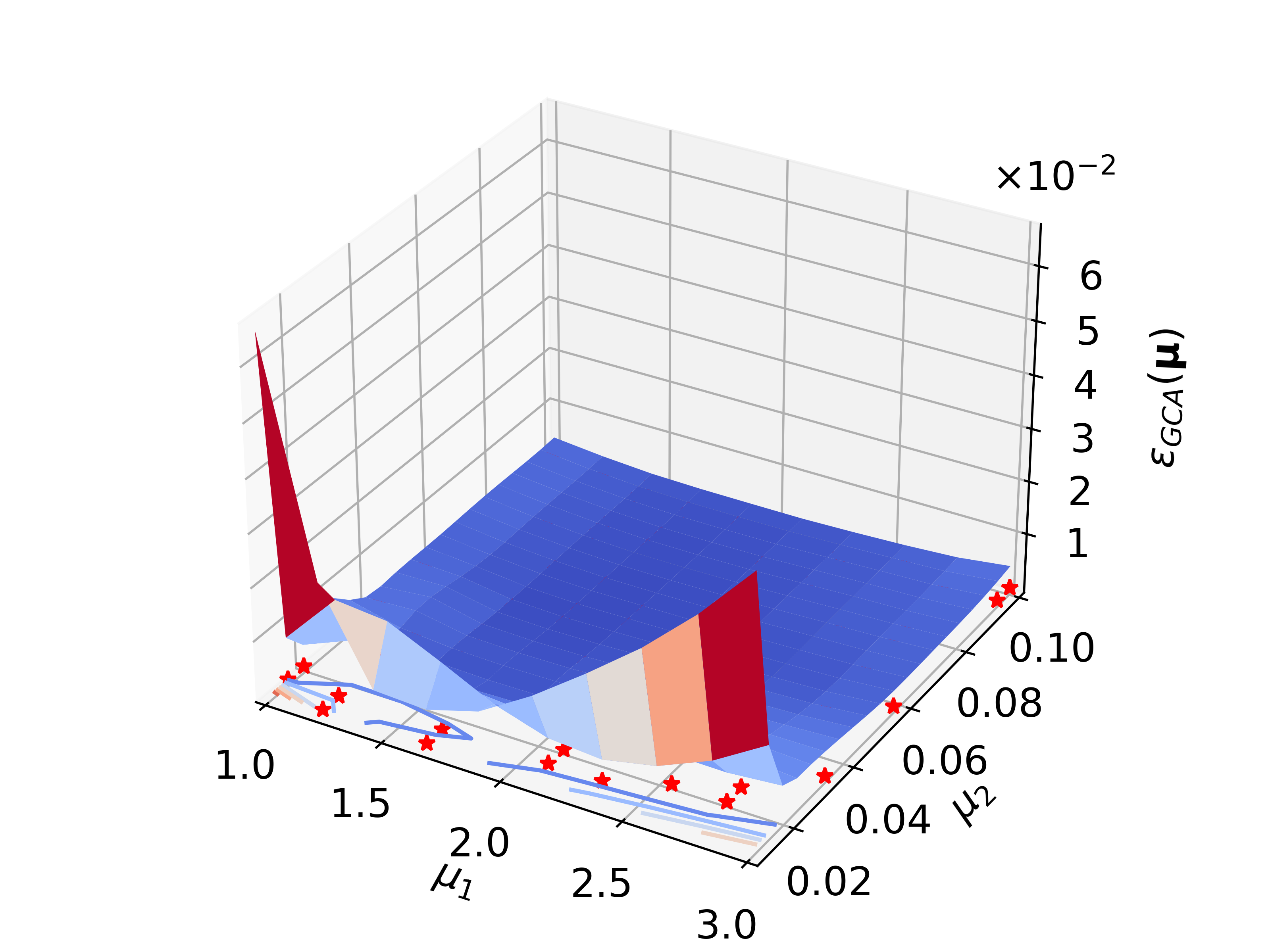}
	\caption{GCA-ROM with pooling}
	\label{fig:error_pooling_graetz_fields_pool}
	\end{subfigure}
	\caption{{GCA-ROM relative errors for the Graetz problem on the dataset $\Xi$, with red markers corresponding to the parameters used in the training set $\Xi_\text{tr}$.}}
	\label{fig:error_pooling_graetz_fields}
\end{figure}

Figure \ref{fig:error_pooling_graetz_fields} shows the investigation of the reconstruction accuracy over the dataset $\Xi$ through the relative errors $\varepsilon_{GCA}(\boldsymbol{\mu})$.
The plain GCA-ROM architecture is still able to learn the low-dimensional evolution in the context of a varying geometry by exploiting only $r_\text{t} = 30\%$ of the dataset information, with the relative error over $\Xi_\text{te}$ always below $6.7 \times 10^{-2}$, with mean $\overline{\varepsilon}_{GCA} =  5.6 \times 10^{-3}$. 
The results for the pooling approach are again comparable, with maximum relative error below $9.3 \times 10^{-2}$ and mean $\overline{\varepsilon}_{GCA} = 7 \times 10^{-3}$, when exploiting $r_\text{p} = 70\%$ of the mesh.
In both cases, we observe high generalization accuracy of the GCA-ROM in capturing the geometric parametrization. This is a result of the introduction of the geometric priors through the GNN architecture endowed with the MoNet convolutional layers. The higher errors are indeed localized for smaller values of the diffusivity parameter, where the phenomenon approaches the advection dominated regime.

\subsection{Navier-Stokes bifurcating system}\label{sec:ns_result}
As the last benchmark, we consider a more complex problem that combines nonlinear terms, parametrized geometry and a vector unknown with non-unique solutions. The so-called Coanda's effect in the Navier-Stokes system is a typical example of fluid dynamic bifurcation problem \cite{PichiArtificialNeuralNetwork2023,PichiDrivingBifurcatingParametrized2022a, KhamlichModelOrderReduction2022, TonicelloNonintrusiveReducedOrder2022}. For such problems, a small variation of the viscosity, in the neighborhood of the bifurcation points, produces a sudden change in the qualitative behavior of the solution, causing the co-existence of different states for the same parameter value. 
Here, we are not interested in recovering all the possible solutions. Rather, we seek to investigate the GCA-ROM performance in the approximation of the transition between two qualitatively different stable states, i.e.\ the symmetric flux and the wall-hugging profile.

We study the steady-state Navier-Stokes system on a parametrized sudden-expansion channel geometry $\Omega(\mu_1)$, depicted in Figure \ref{fig:24_coanda_effect_domain}, discretized with $N_h = 2719$ nodes. 
{
    \begin{figure}[!ht]
    \centering
\begin{minipage}{\textwidth}
	\centering
    \hspace*{-0.3cm}\begin{tikzpicture}[scale=2]
    \draw[fill=blue] (-0.3,0.3) -- (0.8,0.3) -- (0.8,0) -- (5,0) -- (5,0.9) -- (0.8,0.9) -- (0.8,0.6) -- (-0.3,0.6) -- cycle;
    \node[below] at (-0.3,0.3) {\normalsize{(0, 3.75 - $\mu_1$)}};
    \node[above] at (-0.3,0.6) {\normalsize{(0, 3.75 + $\mu_1$)}};
    \node[above] at (5,0.9) {\normalsize{(50, 7.5)}};
    \node[below] at (5,0) {\normalsize{(50, 0)}};
    \node[above] at (0.8,0.9) {\normalsize{(10, 7.5)}};
    \node[below] at (0.8,0) {\normalsize{(10, 0)}};
    \draw[->,>=stealth] (-0.5,0.55) -- (-0.3,0.55);
    \draw[->,>=stealth] (-0.5,0.45) -- (-0.3,0.45);
    \draw[->,>=stealth] (-0.5,0.35) -- (-0.3,0.35);
    \node[below] at (0.6,0.3) {\large{$\Gamma_\text{w}$}};
    \node[above] at (0.6,0.6) {\large{$\Gamma_\text{w}$}};
    \node[below] at (2.9,0) {\large{$\Gamma_\text{w}$}};
    \node[above] at (2.9,0.9) {\large{$\Gamma_\text{w}$}};
    \node[white,right] at (-0.3,0.45) {\large{$\Gamma_\text{in}$}};
    \node[white,left] at (5,0.45) {\large{$\Gamma_\text{out}$}};
    \node[white] at (2.9,0.45) {\Large{$\Omega$}};
    \end{tikzpicture}
    \end{minipage}

\begin{minipage}{\textwidth}
\centering
    \includegraphics[width=0.79\textwidth]{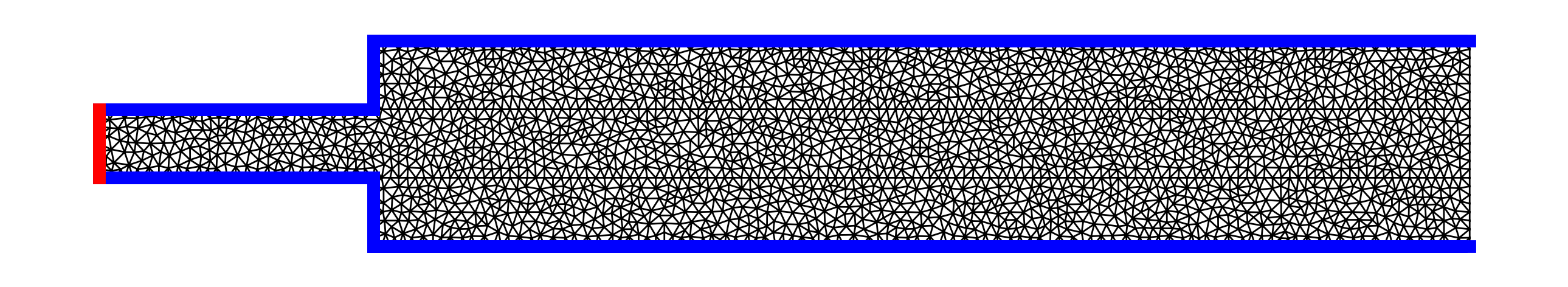}
\end{minipage}
    \caption{{The geometrically parametrized domain for the Coanda effect problem and its computational mesh.}}
    \label{fig:24_coanda_effect_domain}
    \end{figure}}

\noindent The domain is composed by a nozzle, with parametrized width, on which we impose an inlet velocity through a non-homogeneous Dirichlet boundary condition, and a fixed expansion chamber with a Neumann condition at the outlet.

The continuous formulation of the Navier-Stokes problem reads as: given $\boldsymbol{\mu}\in\mathbb{P}$, find the velocity and pressure fields $\textbf{u}(\boldsymbol{\mu}) = (u_x(\boldsymbol{\mu}), u_y(\boldsymbol{\mu}))$ and $p(\boldsymbol{\mu})$, such that
$$
\begin{cases}
-\mu_2\Delta\textbf{u}(\bmu)+(\textbf{u}(\bmu)\cdot\nabla)\,\textbf{u}(\bmu)+\nabla p(\bmu) = 0 & \text{in } \Omega(\mu_1),\\
\nabla\cdot\textbf{u}(\bmu)=0 & \text{in } \Omega(\mu_1),\\
\textbf{u}(\bmu)=\textbf{u}_{\text{in}} & \text{on } \Gamma_\text{in},\\
\textbf{u}(\bmu)=\textbf{0} & \text{on } \Gamma_\text{w}(\mu_1),\\
\mu_2\frac{\partial\textbf{u}}{\partial\textbf{n}}(\bmu)-p(\bmu)\,\textbf{n}=0 & \text{on } \Gamma_\text{out},
\end{cases}
$$
where $\textbf{n}$ is the normal unit vector, $\mu_1$ is the half-width of the inlet channel, $\mu_2$ the kinematic viscosity of the fluid, and  $\textbf{u}_{\text{in}}(\boldsymbol{x}) = \begin{bmatrix} 20(5-y)(y-2.5), \ 0 \end{bmatrix}$ the Poiseuille inflow profile.

The Navier-Stokes benchmark is parametrized by means of the physical and geometric multi-parameter $\boldsymbol{\mu} = (\mu_1, \mu_2) \in \mathbb{P} = [0.5, 2]^2$.
The parametric investigation of this model is physically interesting, and the bifurcation point depends on the physical and geometrical flow's configuration. This means that being able to predict the evolution of the bifurcation point in the parameter space allows to detect the changes in the stability properties of the system. 

Given the nature of the model, with velocity and the pressure fields as unknowns, we show the adaptability of the architecture to vector problems. Thus, we consider a monolithic approach, where the three fields are concatenated as features for each node, rather than a partitioned one, where we recover independently each field\footnote{Notice that this is possible given the data-driven nature of the architecture, i.e.\ we are not exploiting the physics coming from the PDE as in projection-based models.}.

{We report in Table \ref{tab:ns_result} the parametric setting for this benchmark. We highlight that the reason for such a high number of snapshots, $N_\text{S} = 3171$, is that the Finite Element approximation of the bifurcating phenomenon requires a continuation technique to follow the symmetry-breaking branch \cite{pichi_phd}. Despite the increased cardinality of the sampling w.r.t.\ the former test cases, we fixed the training rate to only consider $r_\text{t} = 10\%$ of the original dataset.}

\begin{table}[tbp]
    \centering
    \caption{Configuration and results for the Navier Stokes benchmark.}
    \label{tab:ns_result}
    \begin{tabular}{|c|c|c|c|cc|}
    \hline
    \multicolumn{1}{|c|}{Application} & \multicolumn{1}{c|}{Sampling}  & \multicolumn{2}{c|}{Dataset} & \multicolumn{1}{c}{Component} &  \multicolumn{1}{c|}{$\overline{\epsilon}_{GCA}$}\\
    \hline
    \cellcolor[HTML]{E5E3E3} & $21 \times 151$ & \cellcolor[HTML]{E5E3E3}Train & Test & \cellcolor[HTML]{E5E3E3}Horizontal velocity $u_x$ & \cellcolor[HTML]{E5E3E3}$4.6 \times 10^{-3}$ \\
    \cellcolor[HTML]{E5E3E3}& equispaced & \cellcolor[HTML]{E5E3E3}$r_t = 0.1$ & $r_t = 0.9$ &  Vertical velocity $u_y$ & $2 \times 10^{-2}$ \\
    \cellcolor[HTML]{E5E3E3}\multirow{-3}{*}{Navier-Stokes}& $N_\text{S} = 3171$ & \cellcolor[HTML]{E5E3E3}$N_\text{tr} = 318$ & $N_\text{te} = 2853$ & \cellcolor[HTML]{E5E3E3}Pressure $p$ & \cellcolor[HTML]{E5E3E3}$7.7 \times 10^{-3}$ \\
    \hline
    \end{tabular}
\end{table}

\begin{figure}[!ht]
	\centering
	\subfloat[Horizontal component of the velocity field - $u_x$]{\includegraphics[width=0.49\textwidth]{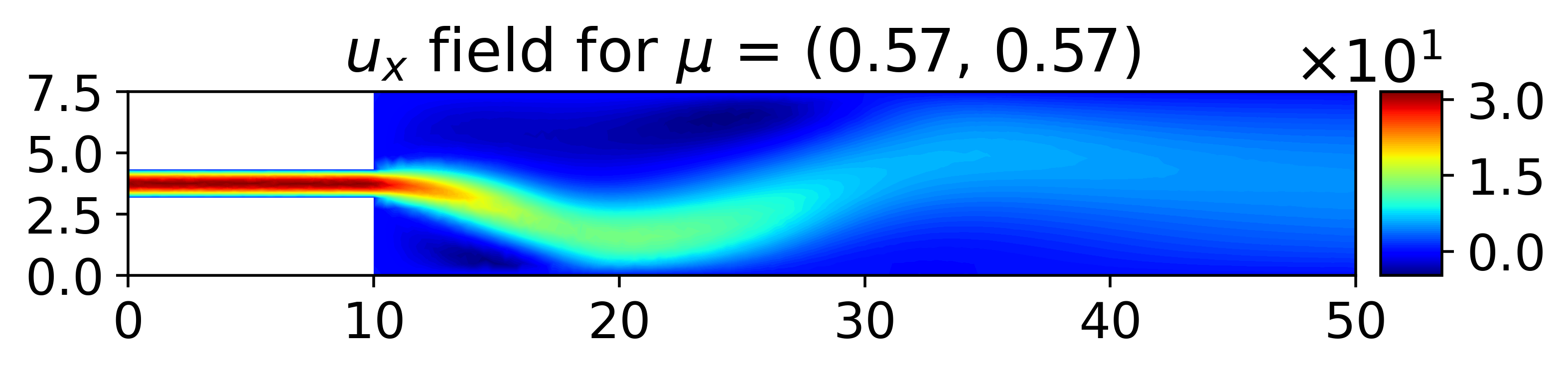}\hfill
	\includegraphics[width=0.49\textwidth]{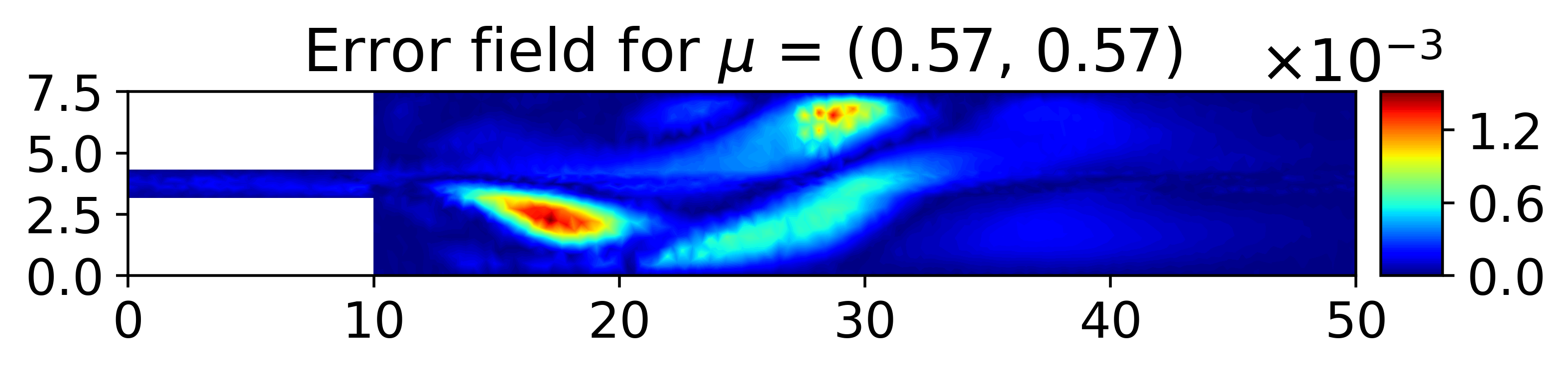}}\\\vspace{0.5cm}
	
	\subfloat[Vertical component of the velocity field - $u_y$]{\includegraphics[width=0.49\textwidth]{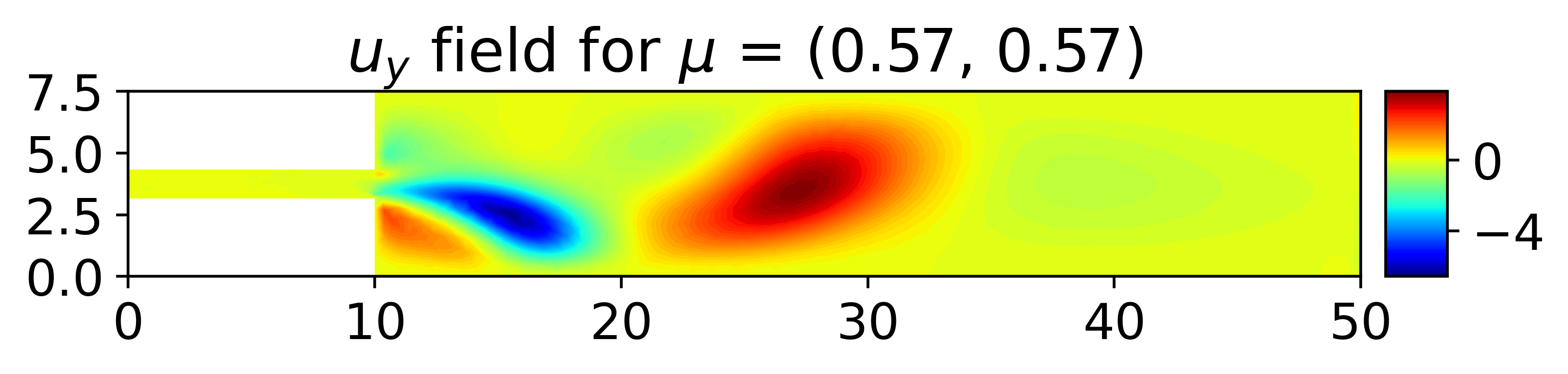}\hfill
	\includegraphics[width=0.49\textwidth]{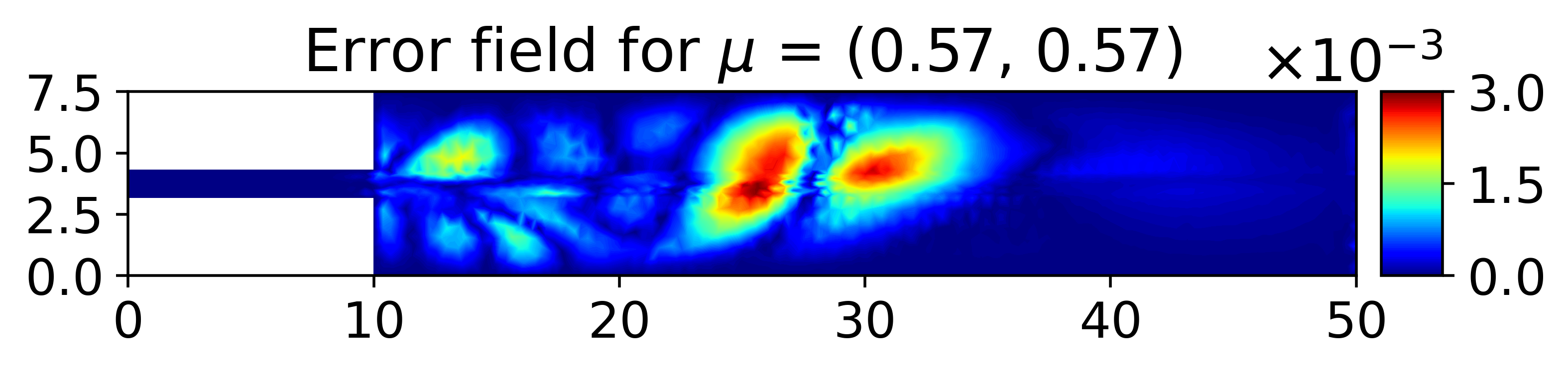}}\\\vspace{0.5cm}
	
	\subfloat[Pressure component - $p$]{\includegraphics[width=0.49\textwidth]{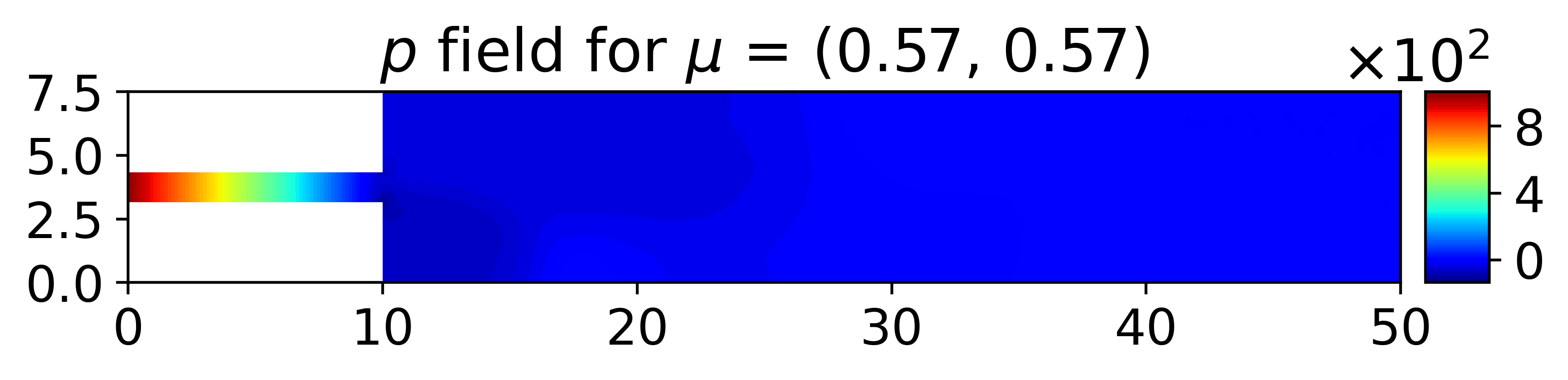}\hfill
	\includegraphics[width=0.49\textwidth]{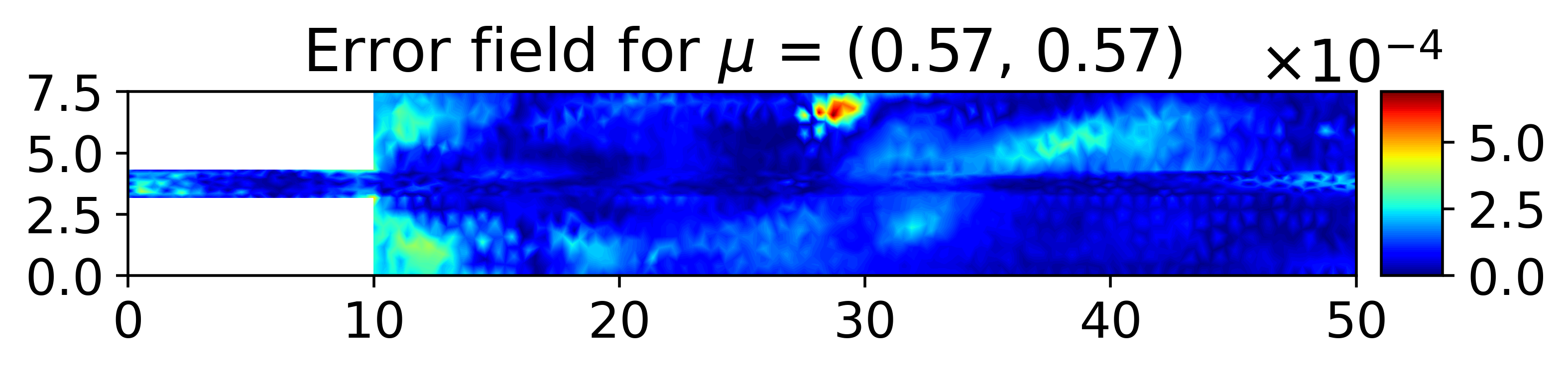}}
	\caption{{Solution and error fields for the components of the Navier-Stokes system for $\boldsymbol{\mu} = (0.57, 0.57)$.}
		\label{fig:error_nopooling_navier_stokes_fields}}
	\end{figure}

In Figure \ref{fig:error_nopooling_navier_stokes_fields} we plot the three components of the Navier-Stokes solution for $\boldsymbol{\mu} = (0.57, 0.57) \in \Xi_\text{te}$, and the GCA-ROM relative errors w.r.t.\ the high-fidelity snapshots. The errors are mainly localized in the regions where the flow suddenly changes, and among the three components, the vertical velocity is the one reconstructed least accurately.
The components recover the wall-hugging behavior with the flow attaching to the bottom boundary, upon which we built the dataset.

\begin{figure}[!ht]
	\centering
    \includegraphics[width=0.4\textwidth]{legend.png}

	\subfloat[Horizontal velocity - $u_x$]{\includegraphics[width=0.33\textwidth,clip=true, trim = 20mm 0mm 10mm 0mm]{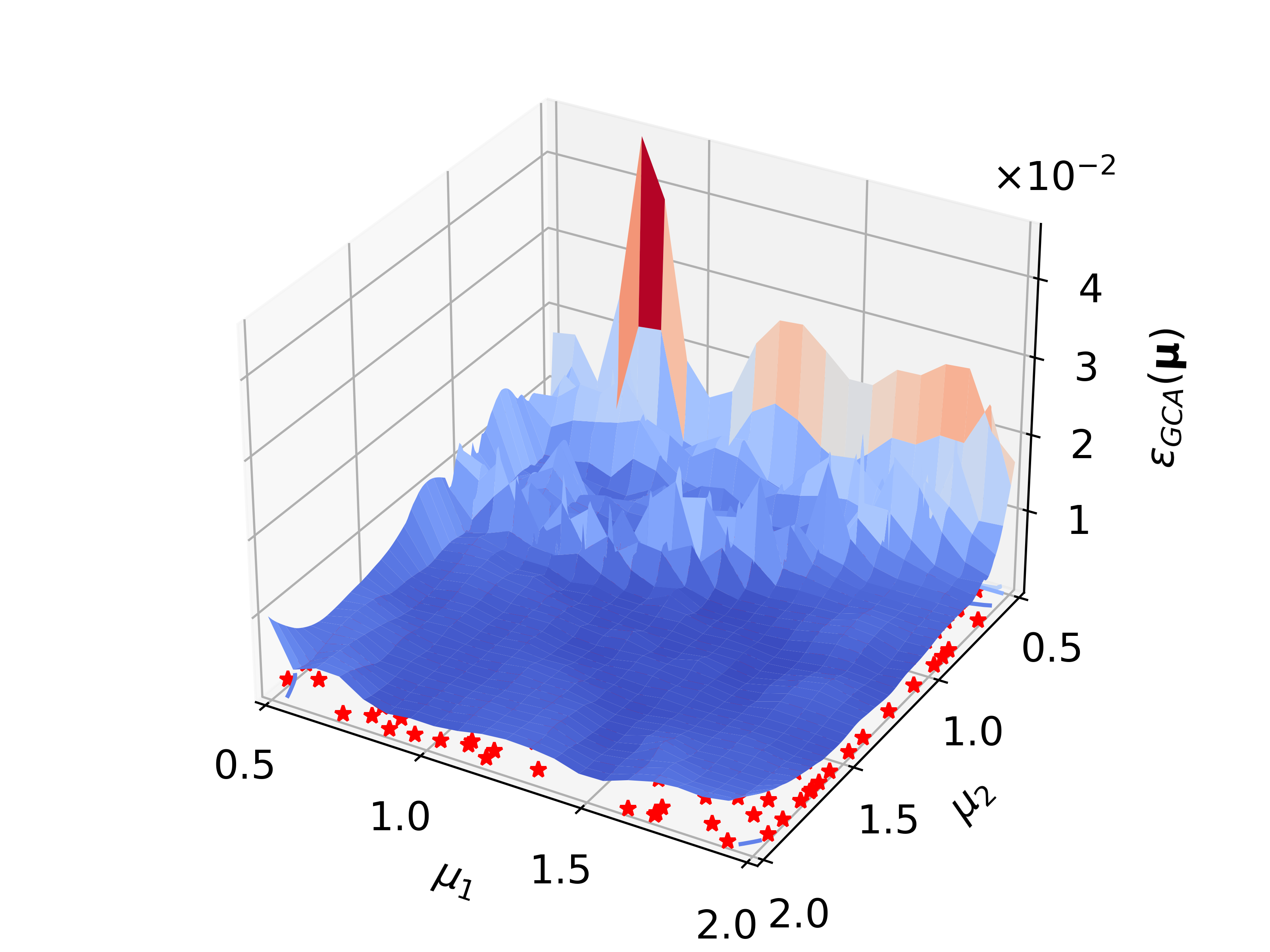}}\hfill
	\subfloat[Vertical velocity - $u_y$]{\includegraphics[width=0.33\textwidth, clip=true, trim = 20mm 0mm 10mm 0mm]{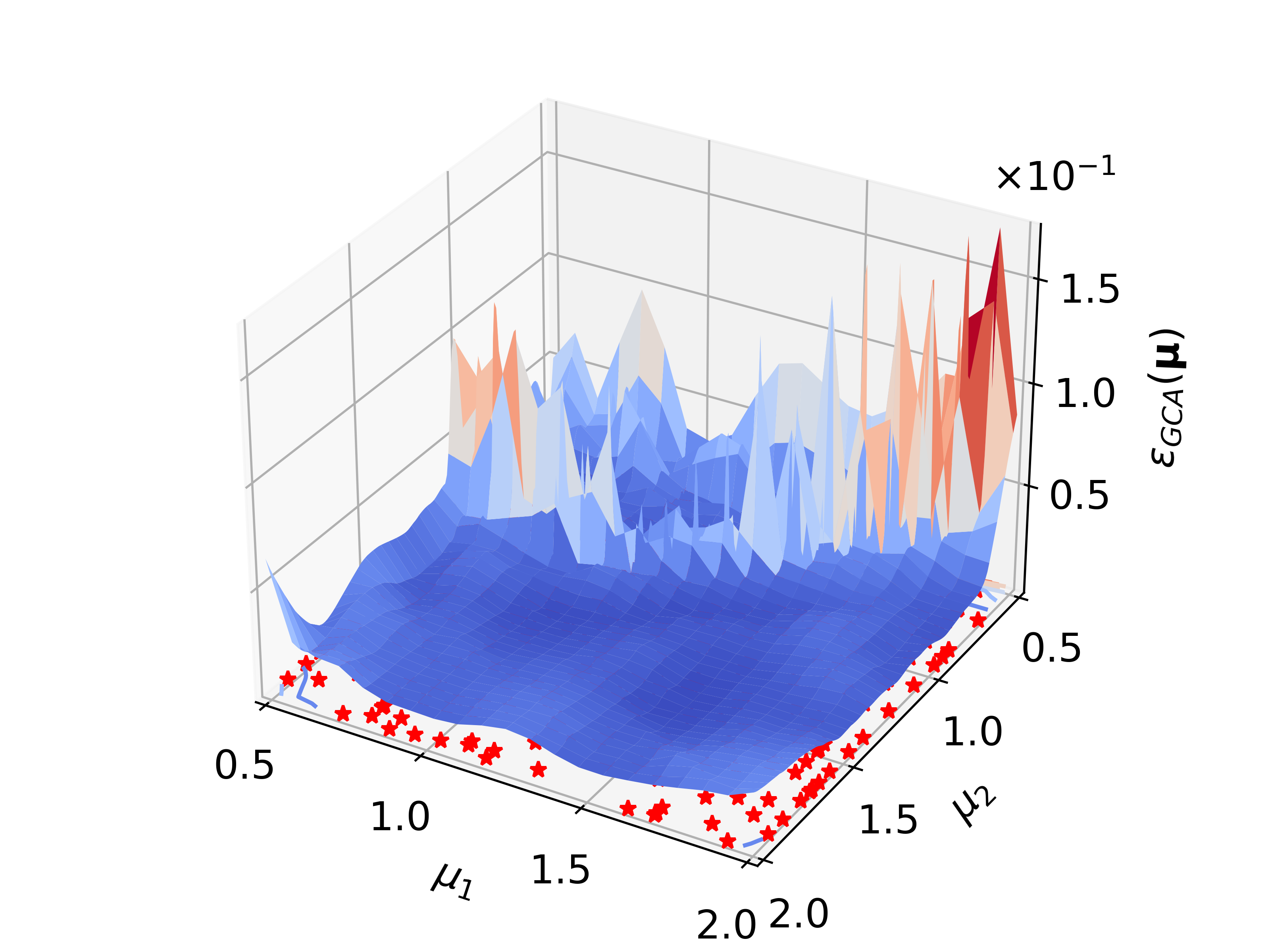}}\hfill
	\subfloat[Pressure - $p$]{\includegraphics[width=0.33\textwidth, clip=true, trim = 20mm 0mm 10mm 0mm]{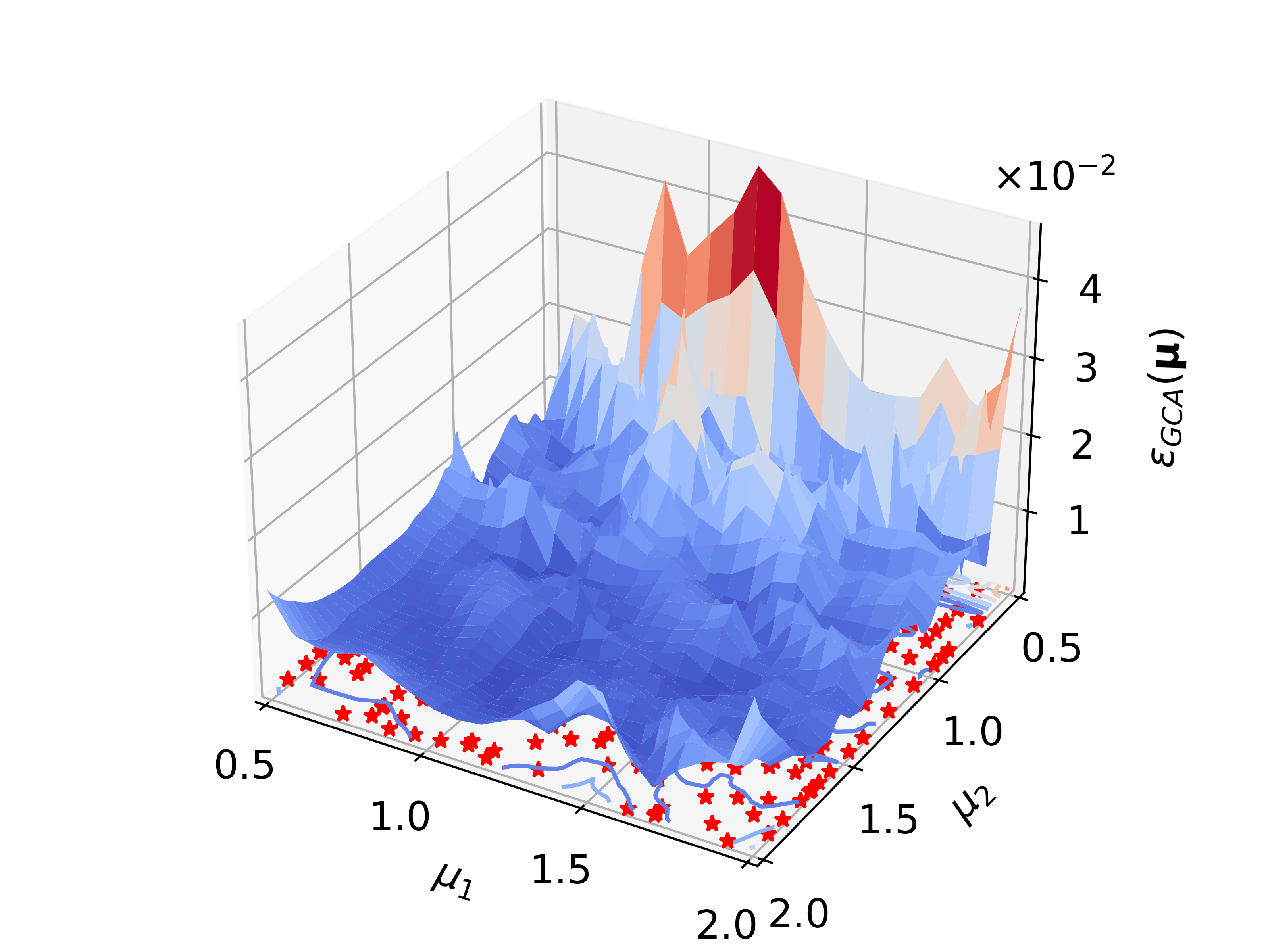}}
	\caption{{GCA-ROM relative errors for the components of the Navier-Stokes test case on the dataset $\Xi$, with red markers corresponding to the parameters used in the training set $\Xi_\text{tr}$.}}
	\label{fig:error_nopooling_navier_stokes_results}
	\end{figure}

\noindent Figure \ref{fig:error_nopooling_navier_stokes_results} depicts the relative errors $\epsilon_{GCA}(\boldsymbol{\mu})$ for the velocity and pressure components over the dataset $\Xi$.
We observe that, even in the vector test case featuring a bifurcating behavior, the GCA-ROM architecture accurately reconstructs all the three components of the Navier-Stokes equations with small mean relative errors, as reported in Table \ref{tab:ns_result}.

On one side, we notice that this could still be considered a low-data regime, especially in the bifurcating context, while on the other this shows that the architecture can benefit when increasing the cardinality of $\Xi_\text{tr}$ without overfitting.
As observed previously, the vertical velocity field is the most difficult to approximate, due to the fact that it is the component responsible for the symmetry breaking phenomenon, and thus the bifurcation.
Higher generalization errors occur again for smaller values of the viscosity. For the vertical velocity, this is maximized for bigger inlet's width, since in that region the geometry still does not admit the wall-hugging profile, while the (bifurcation-agnostic) data-driven approach is influenced by the neighboring training snapshots.

    \subsubsection{Clustering and classification task for bifurcating problems}
    Here, we show a simple yet important exploitation of the low-dimensional manifold encoded in the GCA-ROM architecture's bottleneck. In the bifurcating scenario, from structural to fluid mechanics, deciphering whether a combination of parameters belongs to the critical regime is of importance. Knowing the evolution of the low-dimensional manifold for such problems can help the detection of the bifurcating behavior without performing costly many simulations. Having built a computationally cheap mapping from the parameter to the latent space, the idea is to perform a clustering of the bottleneck's dataset to assign labels for the different states, {and a classification task to predict the state's properties (bifurcation or uniqueness).}
    
    This approach has been explored within linear and non-intrusive ROMs \cite{PichiArtificialNeuralNetwork2023}. Here, we exploit the connection between the reduced coefficients in POD-based techniques and the compressed latent information obtained through the nonlinear encoding. The availability of real-time evaluations of the bottleneck allows to finely sample the parameter space, and obtain an accurate description of its evolution.
    
    We choose the spectral clustering \cite{shi2000normalized}, which exploits eigenvalues' information of the dataset's similarity matrix. Performing the clustering in the similarity space instead of the Euclidean one, this technique is advantageous for problems with irregular and/or unbalanced clusters. In order to distinguish between bifurcating and the uniqueness regime, we only need two clusters, to identify the bifurcation points' curve in $\mathbb{P}$. Given the non-trivial behavior of the bifurcation curve, due to the 
    relation between the viscosity and the inlet's width, we apply $k$-NN {as classification algorithm} for its efficiency in learning nonlinear boundaries.
    {Based on these clustering results, we trained the classifier within a range of $20\%$ to $90\%$ of the labels, and we noticed that the accuracy was always greater than $99.5\%$.}
    
\begin{figure}[t]
    \centering
    \includegraphics[width=0.5\textwidth]{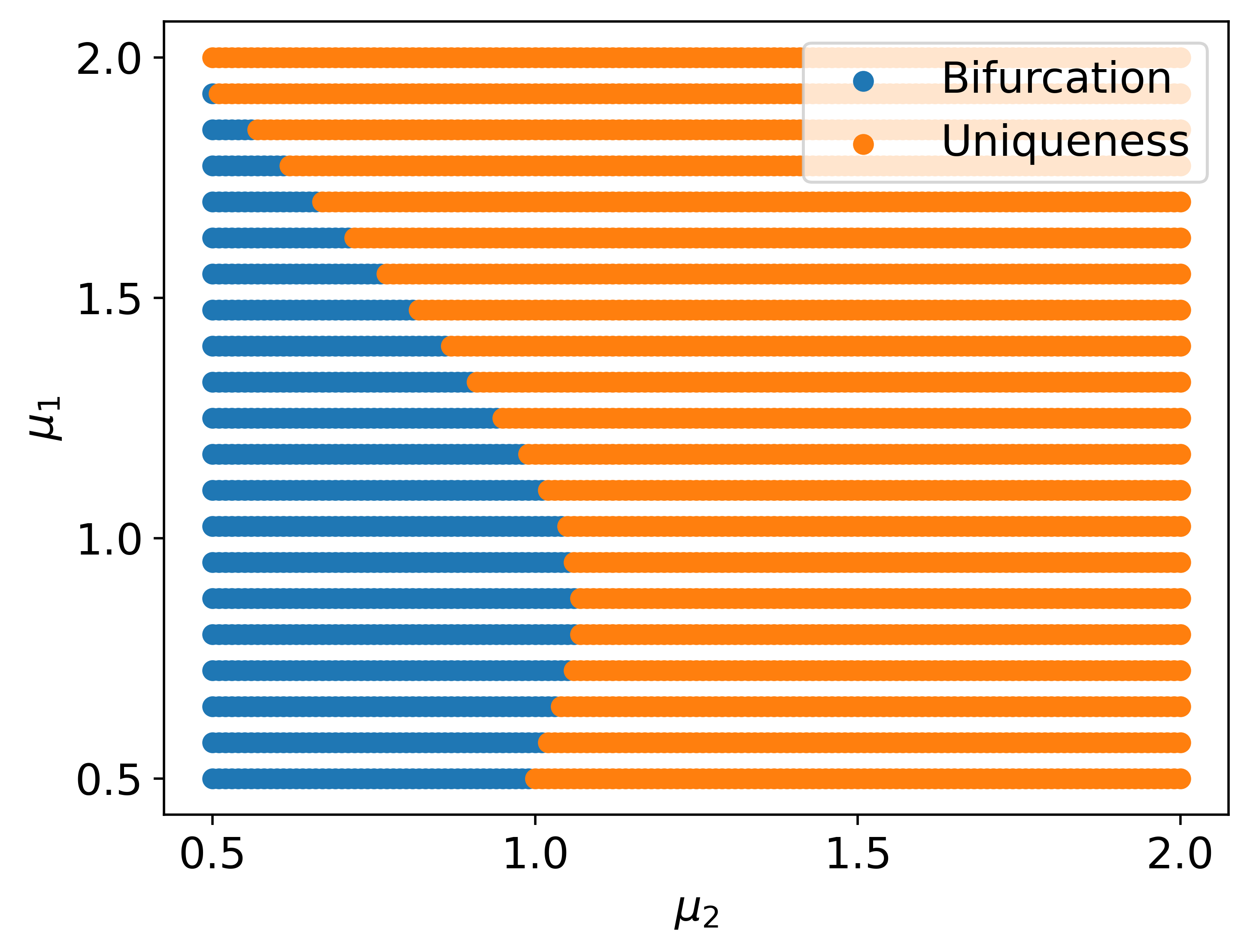}
    \caption{Partitioning of the parameter space in bifurcating and non-bifurcating regimes, obtained with the clustering of the latent manifold based on the GCA-ROM reconstruction of $u_y$.}
        \label{fig:cluster}
    \end{figure}

    In Figure \ref{fig:cluster}, we show the clustered parameter space, clearly identifying the two regimes and detecting the evolution of the bifurcation points. We notice that such a curve was already present in the error analysis in Figure \ref{fig:error_nopooling_navier_stokes_results}, identified by the maximum locations of the GCA-ROM error for the vertical component of the velocity. This means that, breaking the symmetry, the information from $u_y$ can produce more accurate clustering. Finally, we note that the plot in Figure \ref{fig:cluster} is in complete agreement with the results obtained with the Reduced Manifold based Bifurcation diagram (RMB) technique in \cite{PichiArtificialNeuralNetwork2023}.
    
    \subsection{Comparison with linear and nonlinear approaches}
    {In this section, we compare the performance of GCA-ROM with standard linear techniques based on POD, and the original nonlinear DL-ROM architecture.}
    Rather than focusing on the best case scenario, we detail the key features and main drawbacks of these approaches, considering the three scalar benchmarks previously discussed. 
    As metrics, we do not only consider the mean relative error over the testing dataset $\Xi_\text{te}$, but also the offline and online timings, as well as the dimension of the optimization problem.
    
    {
    In Table \ref{table:comparison_error}, we report the accuracy over $\Xi_\text{te}$ of the different methodologies when trained on the same set of snapshots $\Xi_\text{tr}$ with $r_\text{t} = 30\%$. To have a fair comparison, i.e.\ based on the same amount of information, we considered the GCA-ROM architecture without down-sampling, and we fixed the reduced/latent dimension to $n = 15$ for all benchmarks and methodologies.
    
    In particular, we consider the projection error in the POD space, denoted with (POD), the standard POD-Galerkin approach (POD-G), DL-ROM and GCA-ROM.
    As expected, when considering a large enough latent dimension, thanks to the exponential decay of the reduced error, projecting in the POD space gives the best results in terms of approximation accuracy.
    POD-G still produces better performance as compared to ML-based techniques (at least for large enough values of $n$) when dealing with easily reducible problems, such as Poisson or Graetz.
    
    On the contrary, nonlinear techniques show greater capabilities when dealing with advection-dominated phenomena, which are characterized by a slow Kolmogorov $n$-width decay. In all cases, GCA-ROM provided better results than its CNN-based counterpart, even beating the POD-G approach for the Advection test case. This confirms that augmenting the learning procedure with geometrical biases helps the optimization step.
    
    Other non-intrusive POD-based approaches, such as POD-NN or PODI, have been investigated, but the low-data regime we are interested in is a well-known issue for these methodologies. Indeed, both non-intrusive regression and interpolation approaches usually require a large training dataset, without which the learning capabilities are compromised, and the model is unable to generalize.}

    \begin{table}[tbp]    
        \centering
        \caption{Mean relative errors of the scalar benchmarks over $\Xi_\text{te}$ for POD with projection, POD-G, DL-ROM and GCA-ROM techniques, with reduced/latent dimension $n = 15$.}
        \label{table:comparison_error}
        \begin{tabular}{ccccc}
        \hline
        \multicolumn{1}{c}{Application}  & \multicolumn{1}{c}{POD} & \multicolumn{1}{c}{POD-G} & \multicolumn{1}{c}{DL-ROM} & \multicolumn{1}{c}{GCA-ROM} \\ \hline
        \cellcolor[HTML]{E5E3E3}Poisson & \cellcolor[HTML]{E5E3E3}$9.9\times 10^{-5}$ & \cellcolor[HTML]{E5E3E3}$1.0\times 10^{-4}$ & \cellcolor[HTML]{E5E3E3}$1.5\times 10^{-2}$ & \cellcolor[HTML]{E5E3E3}$7.8\times 10^{-3}$ \\
        Advection & $3.1\times 10^{-2}$ & $4.0\times 10^{-2}$ & $5.0\times 10^{-2}$ & $2.4\times 10^{-2}$ \\
        \cellcolor[HTML]{E5E3E3}Graetz & \cellcolor[HTML]{E5E3E3}$2.3\times 10^{-4}$ & \cellcolor[HTML]{E5E3E3}$2.4\times 10^{-4}$ & \cellcolor[HTML]{E5E3E3}$1.4\times 10^{-2}$ & \cellcolor[HTML]{E5E3E3}$6.8\times 10^{-3}$ \\
        \hline
        \end{tabular}
    \end{table}

\begin{figure}[t]
        \centering
        \begin{minipage}{0.4\textwidth}
            \centering
            \includegraphics[width=\textwidth]{legend.png}
    
            \includegraphics[width=\textwidth, clip=true, trim = 20mm 0mm 10mm 0mm]{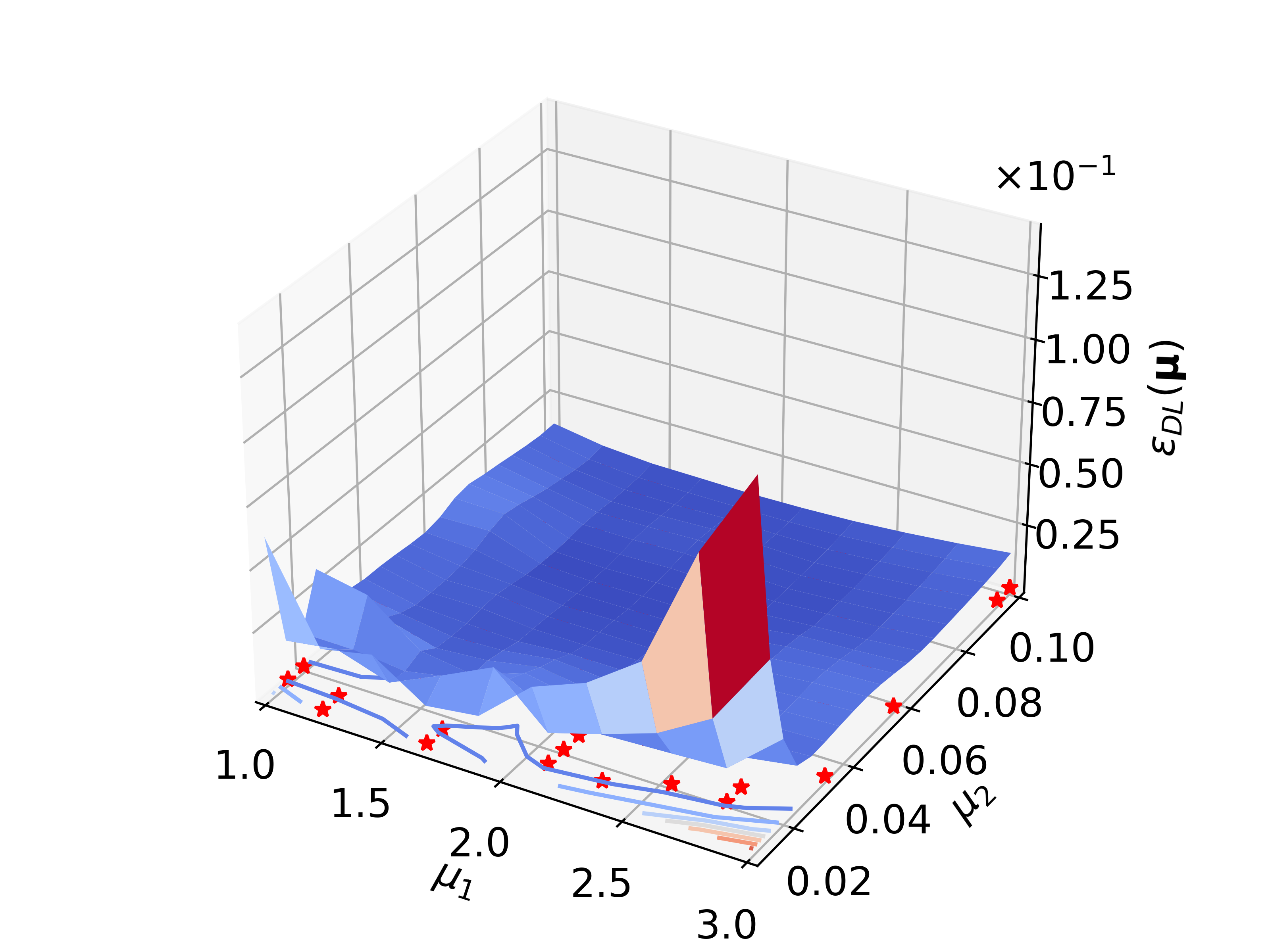}
        \end{minipage}\qquad
        \begin{minipage}{0.4\textwidth}
            \includegraphics[width=\textwidth]{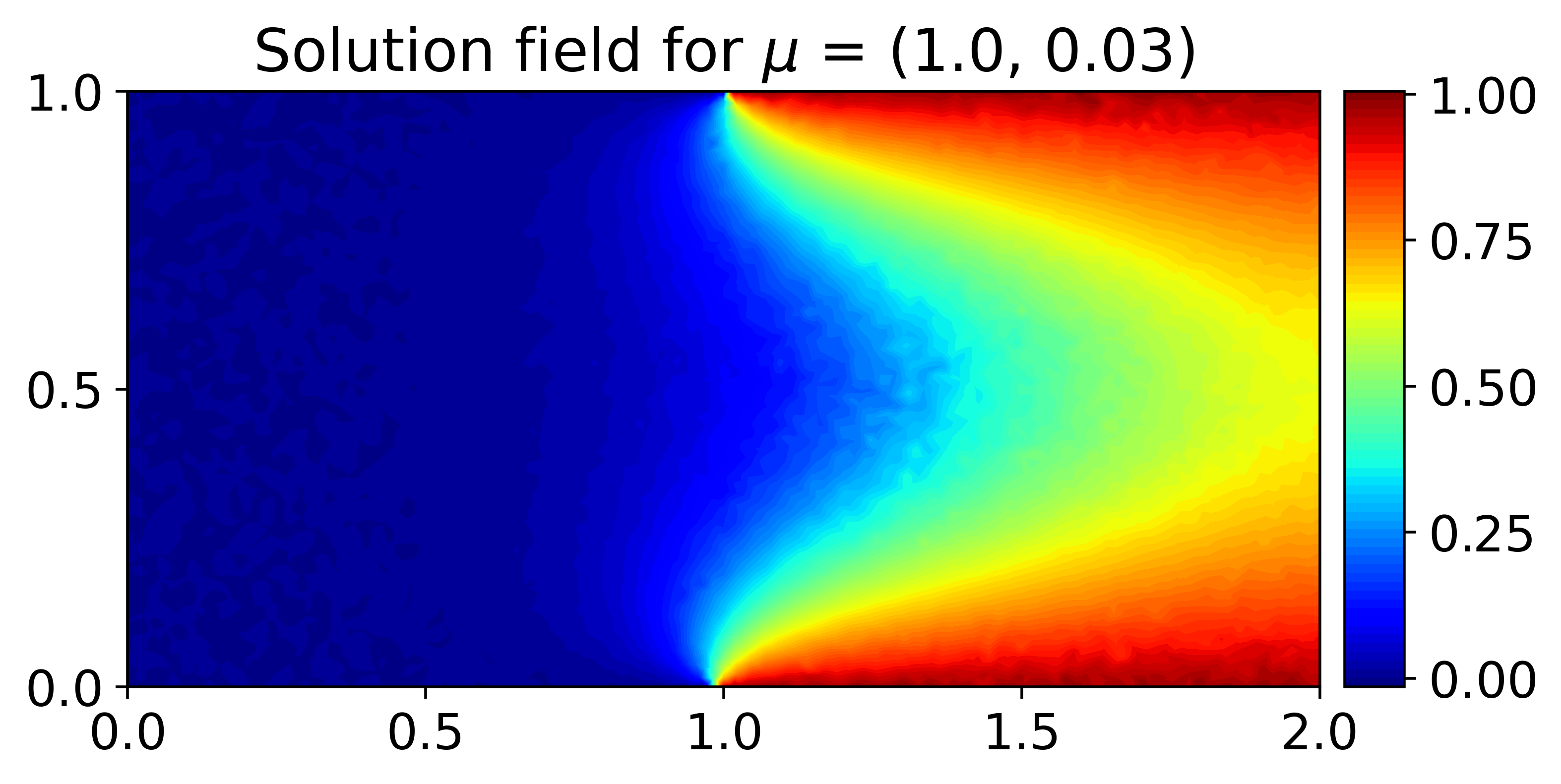}\\
        
            \includegraphics[width=\textwidth]{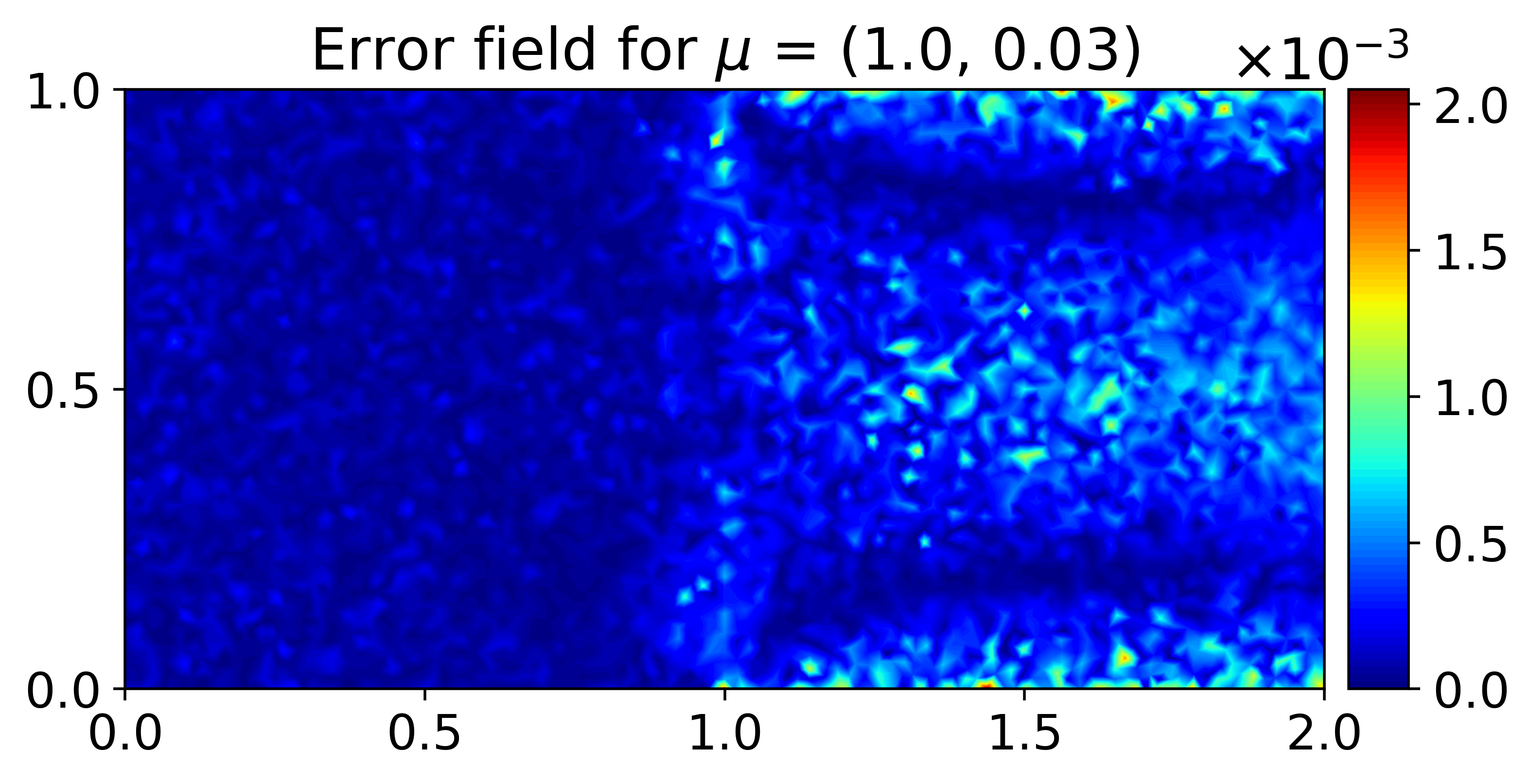}
        \end{minipage}
        \caption{{DL-ROM relative error for the Graetz problem on the dataset $\Xi$, and solution and error fields for $\boldsymbol{\mu} = (1, 0.03)$, left and right respectively.}}
        \label{fig:compareDLrom}
    \end{figure}

    As an example, we show the DL-ROM results for the geometrically parametrized Graetz problem. Relying on an affine map transformation between the original and reference domains, this corresponds to the easiest setting in which to investigate the performance for unstructured grids in advanced applications.
    In Figure \ref{fig:compareDLrom} we plot the relative error of the DL-ROM approach $\varepsilon_{DL}(\boldsymbol{\mu})$, and the solution and error fields of the Graetz model obtained for $\boldsymbol{\mu} = (1, 0.03) \in \Xi_\text{te}$. While the original approach is still able to learn the main features of the model, with possible implications still to be understood \cite{ZhangUnderstandingDeepLearning2021,FrancoApproximationBoundsConvolutional2023}, and a mean $\overline{\varepsilon}_{DL} =  1.6\times 10^{-2}$, we notice that it also produces an unsatisfactory maximum one around $1.4\times 10^{-1}$. This is localized in the area corresponding to the maximum stretching for the domain, where the geometric consistency is even more crucial.
    
    Furthermore, given the non-locality of the convolutional operator in DL-ROM, the reconstruction of the solution field is in general not smooth, as the one depicted in Figure \ref{fig:compareDLrom}. The error shows a speckled pattern affecting the solution in its entire domain, and producing an oscillating reconstruction of the thermal field. Comparing the results with the ones in Figure \ref{fig:error_nopooling_graetz_fields}, the same effect is much less evident, producing a smooth field even considering a highly deformed geometrical configuration and a lower viscosity value.

    As regard the efficiency, we now discuss both the offline and online complexity of the two methodologies.
    While the GCA-ROM approach takes as input the high-fidelity fields already defined on their unstructured geometry, DL-ROM approaches require reshaping the data as a matrix, applying a padding strategy to reach the desired dimensions and apply convolution. This artificial addition of zeroes could potentially compromise the learning procedure, as well as increasing its computational complexity.
    
    In Table \ref{table:comparisondlrom} we compare GCA-ROM and DL-ROM in terms of training and testing times, and number of parameters involved in the optimization process.
    It is worth noting that message passing is a computationally demanding operation for GNNs, but exploiting local convolutions in such context, the sharing of weights reduce significantly the amount of learnable parameters in the network. Depending on the dimension of the filters, we notice a reduction of order three to four. Indeed, the amount of weights exploited by DL-ROM heavily depends on the filter sizes, while for GCA-ROM it is almost the same, i.e.\ the cost of the weighting functions in the Gaussian kernels of MoNet is dominated by the one of the feed-forward layer. For this reason, a more in depth investigation of the up-sampling and down-sampling procedures, could drastically improve the efficiency of the presented methodology. 
    
    All simulations have been performed on a workstation equipped with an Nvidia Quadro GP100 GPU. We report both CPU and GPU times, showing that GCA-ROM is computationally affordable even when only modest hardware is available\footnote{We remark that computational GPU times in the GCA-ROM setting can be further optimized by exploiting new releases of PyTorch Geometric, enabling even faster computations.}. Testing times are including the computation of the reconstructed solutions for all the parameter values in $\Xi$.

    \begin{table}[tbp]
        \centering
        \caption{Comparison between the number of trainable parameters, training, and testing times for DL-ROM and GCA-ROM with different filter sizes.}\label{table:comparisondlrom}
        \begin{tabular}{cccccc}
        \hline
        \multicolumn{1}{c}{Method} & \multicolumn{1}{c}{Device} & \multicolumn{1}{c}{Filters} & \multicolumn{1}{c}{Parameters} & \multicolumn{1}{c}{Training time (s)} & \multicolumn{1}{c}{Testing time (s)} \\ \hline
        \multirow{4}{*}{DL-ROM}  & \multirow{2}{*}{CPU} & 3x3 & \num{8476109} & \num{66518}  & 9.49 \\
        &  & 5x5 & \num{6592461} & \num{62060} & 9.74 \\ 
        & \cellcolor[HTML]{E5E3E3} & \cellcolor[HTML]{E5E3E3}3x3 & \cellcolor[HTML]{E5E3E3}\num{8476109} & \cellcolor[HTML]{E5E3E3}\num{1172}  & \cellcolor[HTML]{E5E3E3}8.93 \\
        & \cellcolor[HTML]{E5E3E3}\multirow{-2}{*}{GPU} & \cellcolor[HTML]{E5E3E3}5x5 & \cellcolor[HTML]{E5E3E3}\num{6592461} & \cellcolor[HTML]{E5E3E3}\num{1323}  & \cellcolor[HTML]{E5E3E3}9.29 \\ 
        \multirow{4}{*}{GCA-ROM} & \multirow{2}{*}{CPU} & 3 & \num{2088682} & \num{3967} & 13.94 \\ 
        &  & 5 & \num{2088694} & \num{6944} & 14.86 \\ 
        & \cellcolor[HTML]{E5E3E3} & \cellcolor[HTML]{E5E3E3}3 & \cellcolor[HTML]{E5E3E3}\num{2088682} & \cellcolor[HTML]{E5E3E3}\num{553} & \cellcolor[HTML]{E5E3E3}15.43 \\ 
        & \cellcolor[HTML]{E5E3E3}\multirow{-2}{*}{GPU} & \cellcolor[HTML]{E5E3E3}5 & \cellcolor[HTML]{E5E3E3}\num{2088694} & \cellcolor[HTML]{E5E3E3}\num{714} & \cellcolor[HTML]{E5E3E3}14.92  \\
        \hline
        \end{tabular}
        \end{table}

\section{Conclusions}

This work showcases the performance of deep learning based on graph neural networks to learn reduced representations of PDEs' solutions. The use of geometric information through networks helps to process unstructured information properly. Introducing a consistent way of treating complex and potentially parametrized geometries, the proposed GCA-ROM architectures improve the generalization property of deep learning approaches, thus paving the way for the analysis of many developments based on graph neural networks.

We tested our novel methodology on a diverse set of parametric problems, exhibiting different characteristics for challenging linear and nonlinear reduced order approaches. In particular, we considered both complex physical behaviors, as advection-dominated problems and bifurcating systems, and unstructured non-rectangular domains, characterized by parameter dependent geometries. In all scenarios the methodology showed great accuracy, even in the low-data regime, exploiting as little as the 30$\%$ of the original dataset.

GCA-ROM was compared with the original DL-ROM approach, highlighting the advantages of adopting graph-based operations in deep learning task related to PDE's solutions. We achieved better accuracy while reducing the computational cost and number of trainable parameters. As a result, the method learns the general pattern of the graphs, thus preventing the autoencoder from overfitting. Given the diffusive nature of the message passing algorithm, based on localized convolutions, GCA-ROM seems to learn smoother fields as compared to standard CNN-based approaches. 

{The application of convolutional layers, particularly of those defined on graphs, shows robust performance while learning across different mesh topologies. In addition, the chosen graph convolution operation, MoNet, incorporates the information regarding the node's position, and the edge attributes as an attention mechanism. These two key attributes enable the method to adapt to any topological configuration efficiently.
However, the current algorithm relies on a fixed size of the layers for the fully connected multilayer perceptron, connecting the encoder/decoder structures to the latent space, thus requiring fixed input/output shapes. A way to solve this issue could be to exploit the pooling operation to reduce the mesh to a fixed number of nodes, and no longer to a percentage $r_\text{p}$ of the original dimension.}

We considered a plain GCA-ROM architecture, and its pooling-based version
comprising also down- and up-sampling operations to coarsen and refine the mesh during learning. We stress that pooling procedures are a valuable tool to extract features from datasets, and effectively reduce the dimensionality without adding further complexity to the network. However, these operations are still in the early stage of development, imposing limitations especially in the decoder structure, where costly interpolations have to be computed. Despite this, such approaches could drastically decrease the size of the optimization problem, and their analysis could create bridges with several methodologies, such as multigrid methods \cite{AntoniettiAgglomerationPolygonalGrids2023}, physics-informed approaches \cite{GaoPhysicsinformedGraphNeural2022,kim2022fast} and time-integration schemes.
Another interesting direction is towards three-dimensional problems, with promising explorations for cardiovascular and musculoskeletal applications \cite{SukMeshNeuralNetworks2022,PegolottiLearningReducedOrderModels2023,KneiflLowdimensionalDatabasedSurrogate2023}. Future work will include high-dimensional parameter space, the enforcement of boundary conditions and the combination of GNNs with different discretization techniques, exploiting even further their versatility.

\section*{Acknowledgements}
Federico Pichi acknowledges the support by ``GO for IT" program within a CRUI fund for the project ``Reduced order method for nonlinear PDEs enhanced by machine learning". Beatriz Moya acknowledges the support of the University of Zaragoza, the bank institution Ibercaja, and the institution CAI (project reference IT1/21). She also thanks the Ministry of Science and Innovation of the Government of Spain (project reference AEI /10.13039/501100011033 with grant number PID2020-113463RB-C31). 

\appendix

\bibliographystyle{abbrv}
\bibliography{manuscript.bib}

\begin{appendix}\label{sec:appendix}
    \section{Architectures}
    
    \addtocontents{toc}{\protect\setcounter{tocdepth}{1}}
    
    \subsection{Architectures and hyperparameters' analysis}\label{sec:arch}
    This section summarizes the GCA-ROM structures used for the numerical results, the hyperparameter investigation to understand the properties of the architecture w.r.t.\ each module, and the corresponding best network configurations for the reproducibility of the results.
    
    We briefly recall the notations used to define the methodology. The dataset is formed by $N_S$ snapshots, i.e.\ the high-fidelity solutions with $d \geq 1$ components (scalar or vector problem) defined over the $N_h$ nodes of the mesh. The training rate $r_\text{t}$ determines the amount of data available for the learning task, while the pooling rate $r_\text{p}$ defines the amount of sensors exploited for the down- and up-sampling procedures. As regard the autoencoder architecture, we have hcp a-priori convolutional layers, and hcd a-posteriori convolutional layers when the pooling module is considered. The connection with the bottleneck of size $n$ is obtained through a 2-layers fully-connected neural network (FC) of size $\text{ffn}$. The parameter map, needed for the online evaluation, is formed by an MLP with $n_l$ neurons in each of the 5-layers. The parameter $\lambda$ is the weight in the loss term balancing the autoencoder and the MLP.
    The dimensionality of the parameter space $\mathbb{P}$ is equal to $P = 2$ for all benchmarks.
    The learning rate  and the weight decay have been fixed throughout all the simulations to the values $l_\text{r} = 10^{-3}$ and $w_\text{d} = 10^{-5}$. In Table \ref{table:architecture} we schematize the GCA-ROM architecture described in Figures \ref{fig:gca_rom_offline} and \ref{fig:gca_rom_online}, and used for the scalar and vector benchmarks. The shaded area corresponds to the pooling module which can be enabled. Otherwise, we consider $r_\text{p} = 100\%$. 
    
\begin{figure}[!ht]
    \centering
    \includegraphics[width=\textwidth]{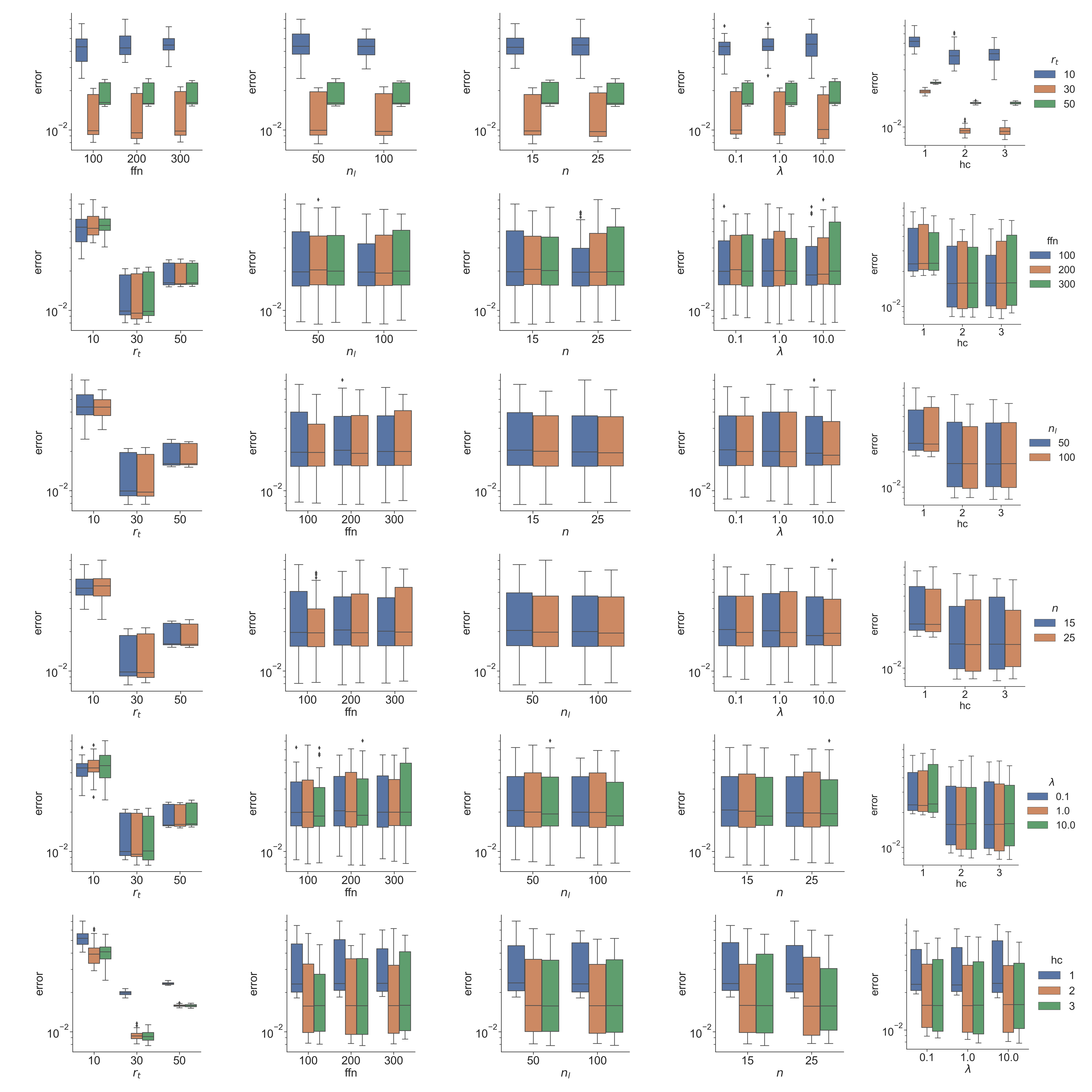}
    \caption{Box-plot of the mean relative error $\varepsilon_{GCA}(\boldsymbol{\mu})$ for the Poisson test case.}
        \label{fig:box_error_nopooling_poisson}
    \end{figure}

\begin{figure}[!ht]
    \centering
    \includegraphics[width=\textwidth]{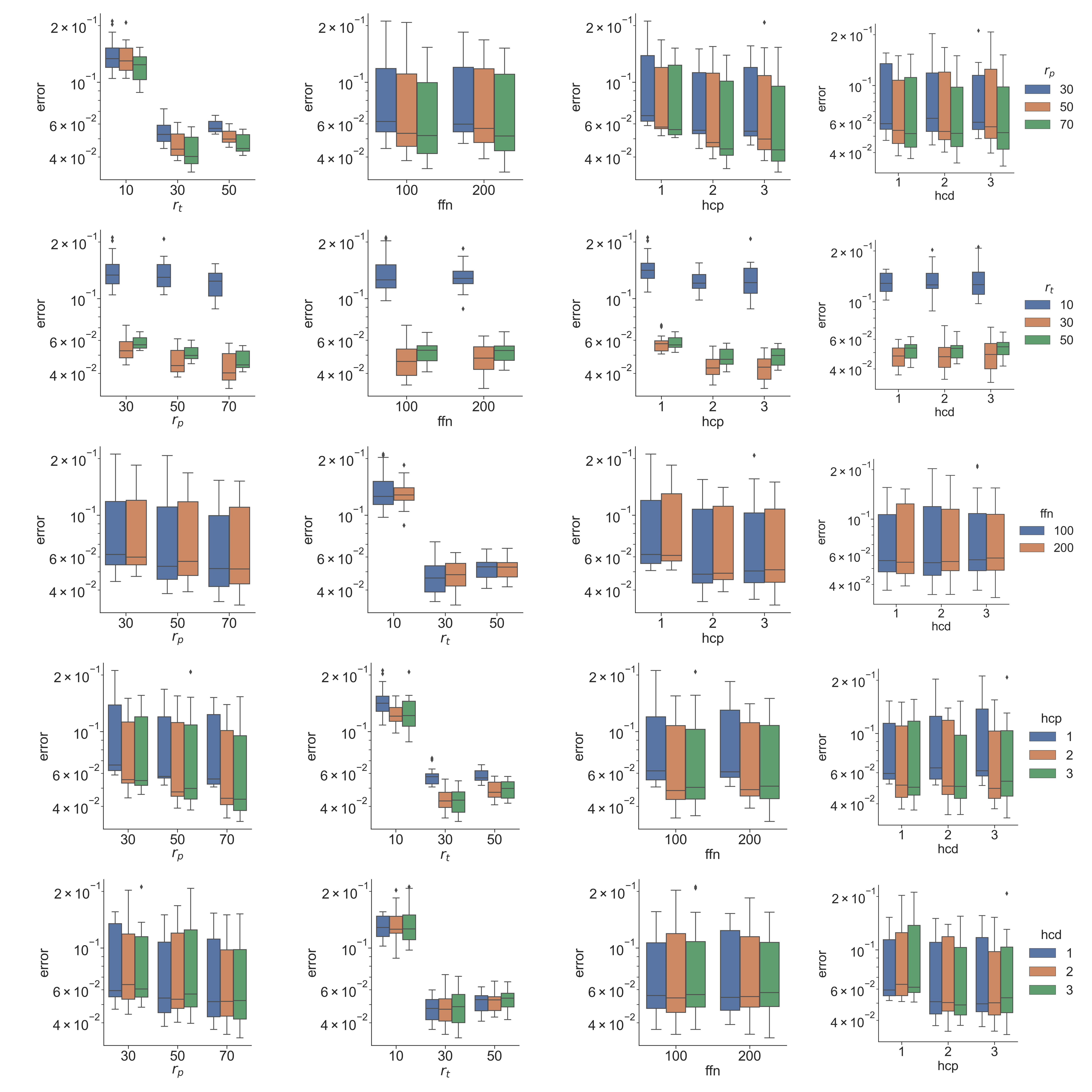}
    \caption{Box-plot of the mean relative error $\varepsilon_{GCA}(\boldsymbol{\mu})$ for the Advection test case with pooling strategy.}
        \label{fig:box_error_pooling_advection}
    \end{figure}

    We show in Figures \ref{fig:box_error_nopooling_poisson} and \ref{fig:box_error_pooling_advection} the hyperparameter investigations obtained with the plain GCA-ROM approach for the Poisson problem, and the GCA-ROM with pooling for the Advection benchmark, respectively. The former illustrates the box-plot of the mean relative errors $\varepsilon_{GCA}(\boldsymbol{\mu})$ for all possible combinations of the following set of hyperparameters:  $r_\text{t} \in [10, 30, 50]$, $r_\text{p} \in [30, 50, 70]$, $\text{ffn} \in [100, 200, 300]$, $n_l \in [50, 100]$, $n \in [15, 25]$, $\lambda \in [0.1, 1, 10]$, and $\text{hc} \in [1, 2, 3]$
    
    In Figures \ref{fig:box_error_pooling_advection} we plot the same quantity, fixing the number of nodes in the latent map $n_l = 50$, the dimension of the bottleneck $n = 25$, and the weight of the loss $\lambda = 1$, while varying: $r_\text{p} \in [30, 50, 70]$, $r_t \in [10, 30, 50]$, $\text{ffn} \in [100, 200]$, and $\text{hcp}, \text{hcd} \in [1, 2, 3]$.
    
    Finally, we report in Table \ref{table:best_config} the best configuration, for the plain and pooling versions of GCA-ROM, resulting from the aforementioned hyperparameter analysis, and which we have used throughout the manuscript. 
    
    \begin{table}[]
        \centering
        \caption{Autoencoder architecture for GCA-ROM approach, with optional pooling module in gray.}\label{table:architecture}
        \begin{tabular}{cccccc}
        \hline
        \multicolumn{1}{c}{Module} & \multicolumn{1}{c}{Layers} & \multicolumn{1}{c}{Input size} &  \multicolumn{1}{c}{Output size} & \multicolumn{1}{c}{Activation} & \multicolumn{1}{c}{Kernel} \\ \hline
        \multirow{5}{*}{Encoder} & $\text{Conv}_{\text{hcp}}$ & $n_b \times N_h \times d$ & $n_b \times N_h \times d$ & ELU & Q \\ 
        & \cellcolor[HTML]{E5E3E3}Pool & \cellcolor[HTML]{E5E3E3}$n_b \times N_h \times d$ & \cellcolor[HTML]{E5E3E3}$n_b \times r_\text{p}N_h \times d$ & \cellcolor[HTML]{E5E3E3} & \cellcolor[HTML]{E5E3E3} \\
        & \cellcolor[HTML]{E5E3E3}$\text{Conv}_{\text{hcd}}$ & \cellcolor[HTML]{E5E3E3}$n_b \times d r_\text{p}N_h$ & \cellcolor[HTML]{E5E3E3}$n_b \times d r_\text{p}N_h$ & \cellcolor[HTML]{E5E3E3}ELU & \cellcolor[HTML]{E5E3E3}Q \\ 
        & \multirow{2}{*}{$\text{FC}$} & $n_b \times d r_\text{p}N_h$ & $n_b \times \text{fnn}$ & ELU & \\
        & & $n_b \times \text{fnn}$ & $n_b \times n$ & Identity & \\ \hline
         \multirow{5}{*}{Decoder} & \multirow{2}{*}{$\text{FC}$}  & $n_b \times n$ & $n_b \times \text{fnn}$ & Identity & \\
         & & $n_b \times \text{fnn}$ & $n_b \times d r_\text{p}N_h$ & ELU & \\
         & \cellcolor[HTML]{E5E3E3}$\text{Conv}_{\text{hcd}}$ & \cellcolor[HTML]{E5E3E3}$n_b \times d r_\text{p}N_h$ & \cellcolor[HTML]{E5E3E3}$n_b \times d r_\text{p}N_h$ & \cellcolor[HTML]{E5E3E3}ELU & \cellcolor[HTML]{E5E3E3}Q \\ 
         & \cellcolor[HTML]{E5E3E3}Unpool & \cellcolor[HTML]{E5E3E3}$n_b \times r_\text{p}N_h \times d$ & \cellcolor[HTML]{E5E3E3}$n_b \times N_h \times d$ & \cellcolor[HTML]{E5E3E3} & \cellcolor[HTML]{E5E3E3} \\
         & $\text{Conv}_{\text{hcp}}$ & $n_b \times N_h \times d$ & $n_b \times N_h \times d$ & ELU & Q \\  \hline
        \multirow{3}{*}{Parameter map} & \multirow{3}{*}{$\text{FC}$} & $n_b \times P$ & $n_b \times n_l$ & $\tanh$ & \\ 
        & & $\vdots$ & $\vdots$ & $\vdots$ & \\
        & & $n_b \times n_l$ & $n_b \times n$ & Identity & \\ \hline
        \end{tabular}
    \end{table}
    
    \begin{table}[]
        \centering 
        \caption{Hyperparameter configuration of the best network for each benchmark.}
        \label{table:best_config}
        \begin{tabular}{cccccccccc}
        \hline
        \multicolumn{1}{c}{Application} & \multicolumn{1}{c}{pooling} & \multicolumn{1}{c}{$r_\text{p}$} &  \multicolumn{1}{c}{$r_\text{t}$} & \multicolumn{1}{c}{$\lambda$} & \multicolumn{1}{c}{btt} & \multicolumn{1}{c}{hcp} & \multicolumn{1}{c}{hcd} & \multicolumn{1}{c}{nd} & \multicolumn{1}{c}{ffn} \\ \hline
        \multirow{2}{*}{Poisson} & \xmark & \xmark & 30 & 10 & 15 & 3 & \xmark & 50 & 100\\
         & \cellcolor[HTML]{E5E3E3}\cmark & \cellcolor[HTML]{E5E3E3}70 & \cellcolor[HTML]{E5E3E3}30 & \cellcolor[HTML]{E5E3E3}1 & \cellcolor[HTML]{E5E3E3}25 & \cellcolor[HTML]{E5E3E3}3 & \cellcolor[HTML]{E5E3E3}3 & \cellcolor[HTML]{E5E3E3}50 & \cellcolor[HTML]{E5E3E3}200 \\ \hline
    
         \multirow{2}{*}{Advection} & \xmark & \xmark & 30 & 10 & 15 & 2 & \xmark & 100 & 200\\
         & \cellcolor[HTML]{E5E3E3}\cmark & \cellcolor[HTML]{E5E3E3}70 & \cellcolor[HTML]{E5E3E3}30 & \cellcolor[HTML]{E5E3E3}1 & \cellcolor[HTML]{E5E3E3}25 & \cellcolor[HTML]{E5E3E3}3 & \cellcolor[HTML]{E5E3E3}3 & \cellcolor[HTML]{E5E3E3}50 & \cellcolor[HTML]{E5E3E3}200 \\ \hline
        
         \multirow{2}{*}{Graetz} & \xmark & \xmark & 30 & 10 & 25 & 2 & \xmark & 50 & 200\\
         & \cellcolor[HTML]{E5E3E3}\cmark & \cellcolor[HTML]{E5E3E3}70 & \cellcolor[HTML]{E5E3E3}30 & \cellcolor[HTML]{E5E3E3}1 & \cellcolor[HTML]{E5E3E3}25 & \cellcolor[HTML]{E5E3E3}2 & \cellcolor[HTML]{E5E3E3}1 & \cellcolor[HTML]{E5E3E3}50 & \cellcolor[HTML]{E5E3E3}100 \\ \hline
    
        Navier-Stokes & \xmark & \xmark & 10 & 1 & 25 & 3 & \xmark & 100 & 200\\ \hline
        \end{tabular}
        \end{table}

    \addtocontents{toc}{\protect\setcounter{tocdepth}{2}}
    \end{appendix}

\end{document}